\documentclass[11pt,english,reqno]{amsart}

\usepackage{layout}
\usepackage{amssymb}
\usepackage{amsfonts}
\usepackage{amsmath,amscd,amsthm}
\usepackage{newlfont}
\usepackage{mathrsfs}
\usepackage[dvips]{graphicx}
\usepackage{epsf}
\usepackage{color}
\usepackage{hhline}
\usepackage{rotating}
\usepackage{verbatim}
\usepackage{enumitem}
\usepackage{nccmath}
\usepackage[normalem]{ulem}
\usepackage{soul}

\everymath{\displaystyle}

\usepackage[footnotesize]{subfigure}
\usepackage{caption}
\newcommand{\E}{\mathbb E}

\renewcommand{\P}{\mathbb P}\newcommand{\R}{\mathbb R}

\newcommand{\indic}{1\negthickspace\text{I}}
\newcommand{\indicat}{{\text{1}}\negthickspace\text{I}}

\newcommand{\cA}{\cal{A}}\newcommand{\cC}{\cal{C}}\newcommand{\cD}{\cal{D}}\newcommand{\cE}{\cal{E}}
\newcommand{\cF}{\cal{F}}\newcommand{\cG}{\cal{G}}
\newcommand{\cN}{\cal{N}}
\newcommand{\cP}{\cal{P}}\newcommand{\cQ}{\cal{Q}}\newcommand{\cR}{\cal{R}}\newcommand{\cS}{\cal{S}}\newcommand{\cT}{\cal{T}}
\newcommand{\cU}{\cal{U}}

\newcommand{\eps}{\varepsilon}

\newcommand{\dd}{\mathrm{d}}

\newcommand{\pass}{\vspace{0.3cm}\noindent}
\newcommand{\noi}{\noindent}

\DeclareMathOperator{\conv}{conv}

\theoremstyle{plain}
\newtheorem{theo}{Theorem}[section]
\newtheorem{prop}[theo]{Proposition}
\newtheorem{lem}[theo]{Lemma}

\newtheorem*{theora}{Theorem A}
\newtheorem*{theorb}{Theorem B}

\theoremstyle{definition}

\definecolor{note}{rgb}{1,0,0}

\definecolor{modif}{rgb}{0,0,1}

\author[P. Calka, Y. Demichel and N. Enriquez]{Pierre Calka\textsuperscript{1}, Yann Demichel\textsuperscript{2} and Nathana\"el Enriquez\textsuperscript{3}}

\address{\textsuperscript{1} Laboratoire de Math\'ematiques Rapha\"el Salem, UMR 6085, Universit\'e de Rouen, avenue de l'Universit\'e,
Technop\^ole du Madrillet, 76801 Saint-Etienne-du-Rouvray, France.}
\email{pierre.calka@univ-rouen.fr}%
\address{\textsuperscript{2} Laboratoire MODAL'X, EA 3454, Universit\'e Paris Ouest Nanterre La D\'efense, 200 avenue de la R\'e\-pu\-bli\-que, 92001 Nanterre, France.}
\email{yann.demichel@u-paris10.fr}%
\address{\textsuperscript{3} Laboratoire MODAL'X, EA 3454, Universit\'e Paris Ouest Nanterre La D\'efense, 200 avenue de la R\'e\-pu\-bli\-que, 92001 Nanterre, France.}
\email{nathanael.enriquez@u-paris10.fr}%

\address{This work was partially supported by the French ANR grant PRESAGE (ANR-11-BS02-003) and the French research group GeoSto (CNRS-GDR3477).}

\title[The Poisson-Voronoi cell around an isolated nucleus]{The Poisson-Voronoi cell around an isolated nucleus}

\begin{document}

\begin{abstract}
Consider a planar random point process made of the union of a point (the origin) and of a Poisson point process with a uniform intensity \textit{outside}
a deterministic set surrounding the origin. When the intensity goes to infinity, we show that the Voronoi cell associated with the origin converges \textit{from above} to a deterministic convex set. We describe this set and give the asymptotics of the expectation of its defect area, defect perimeter and number of vertices. On the way, two intermediary questions are treated. First, we describe the mean characteristics of the Poisson-Voronoi cell conditioned on containing a fixed convex body around the origin and secondly, we show that the nucleus of such cell converges to the Steiner point of the convex body. As in R\'enyi and Sulanke's seminal papers on random convex hulls, the regularity of the convex body has crucial importance. We deal with both the smooth and polygonal cases.
Techniques are based notably on accurate estimates of the area of the Voronoi flower and of the support function of the cell containing the origin as well as on an Efron-type relation.
\end{abstract}



\maketitle

\fussy
\section{Introduction}\label{sec:intro}

\subsection{Main issue and related questions}

\pass

\noi One of the questions we address in this paper is

\pass\begin{center}
\begin{minipage}{13cm}
\textsc{Question 1:} \textit{Given a fixed domain $\cD$ containing the origin $o$ in its interior, what is the geometry of the cell containing $o$ in a Voronoi tessellation generated by the union of $o$ with a Poisson point process whose intensity goes to infinity outside $\cD$ and equals $0$ inside?}
\end{minipage}
\end{center}

\begin{figure}[h!]
\centering
\begin{center}
\includegraphics[trim= 5.5cm 12cm 5cm 7cm, clip, scale=0.39]{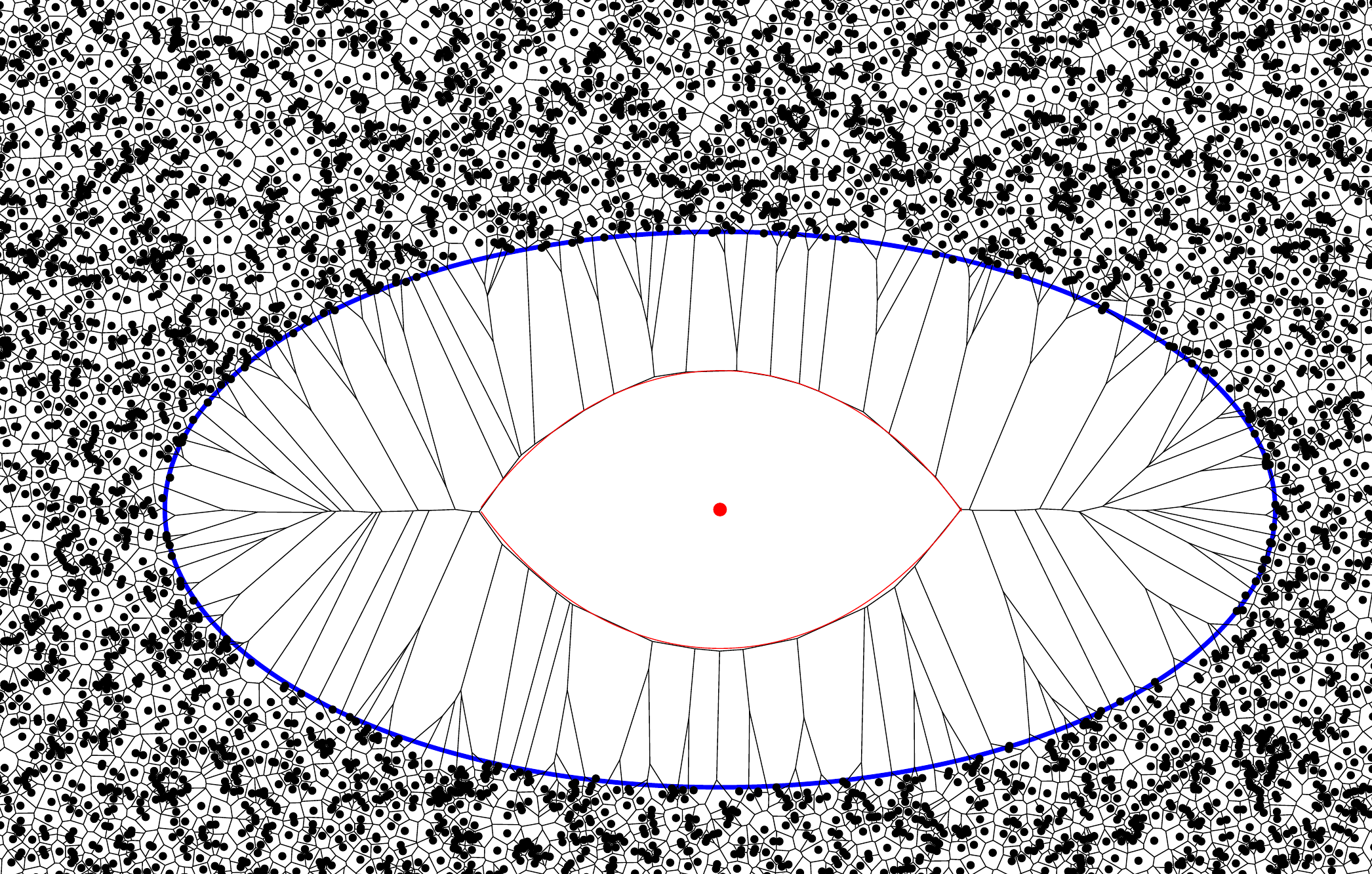}
\hspace*{1cm}
\includegraphics[trim= 3.5cm 2.5cm 3cm 3cm, clip, scale=0.215]{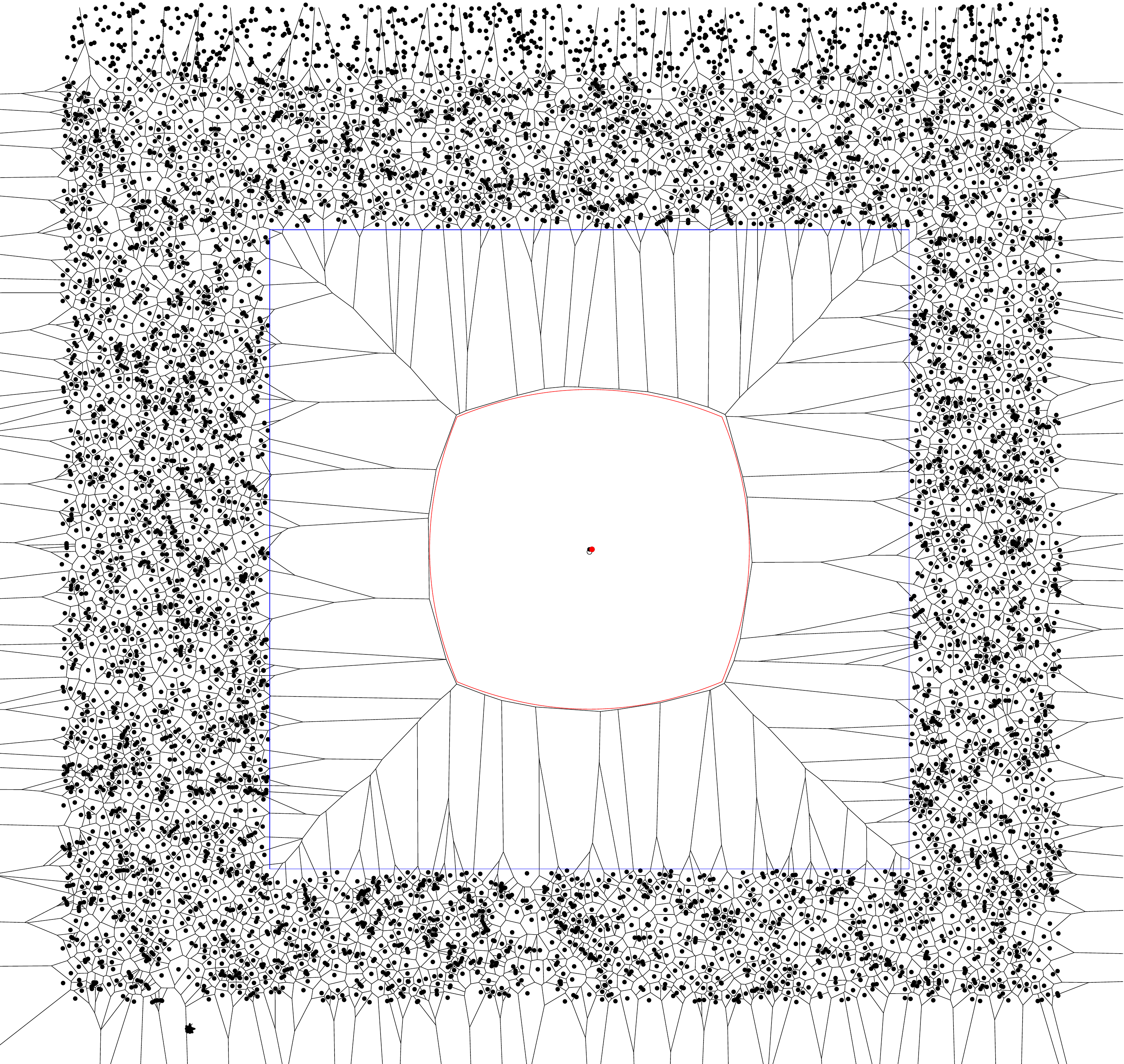}
\end{center}
\caption{Voronoi cells generated by the origin (red) and a Poisson point process intensity $20000$ outside an ellipse (left) and outside a square (right).}
\label{fig:question1}
\end{figure}

\pass This cell converges \textit{from outside} to the convex set of points which are closer to the origin than the boundary of $\cD$ (see Figure \ref{fig:question1}). We are interested in the asymptotic means of the defect area, defect perimeter and number of vertices of this Voronoi cell.

\pass It turns out that the following dual question is the first and key step in answering this problem.

\pass
\begin{center}
\begin{minipage}{13cm}
\textsc{Question 2:} \textit{Given a fixed convex body $K$ containing $o$ in its interior, what is the geometry of the zero-cell of a Voronoi tessellation generated by the union of $o$ with a Poisson point process conditioned on this zero-cell to contain $K$ as its intensity goes to infinity?}
\end{minipage}
\end{center}

\begin{figure}[h!]
\centering
\begin{center}
\includegraphics[trim= 3.5cm 11.5cm 3cm 6cm, clip, scale=0.2]{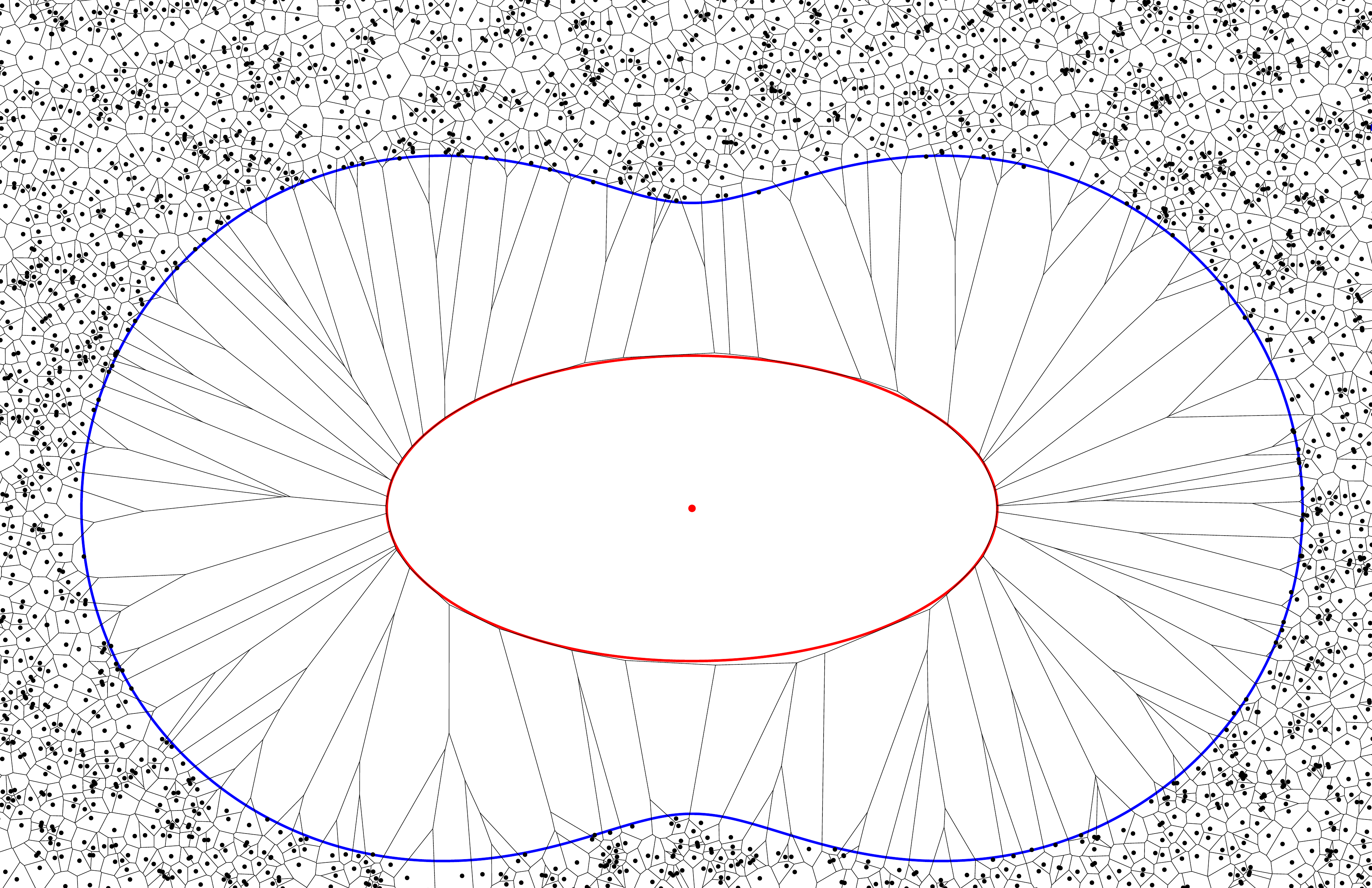}
\hspace*{1cm}
\includegraphics[trim= 3.5cm 2.5cm 3cm 3cm, clip, scale=0.205]{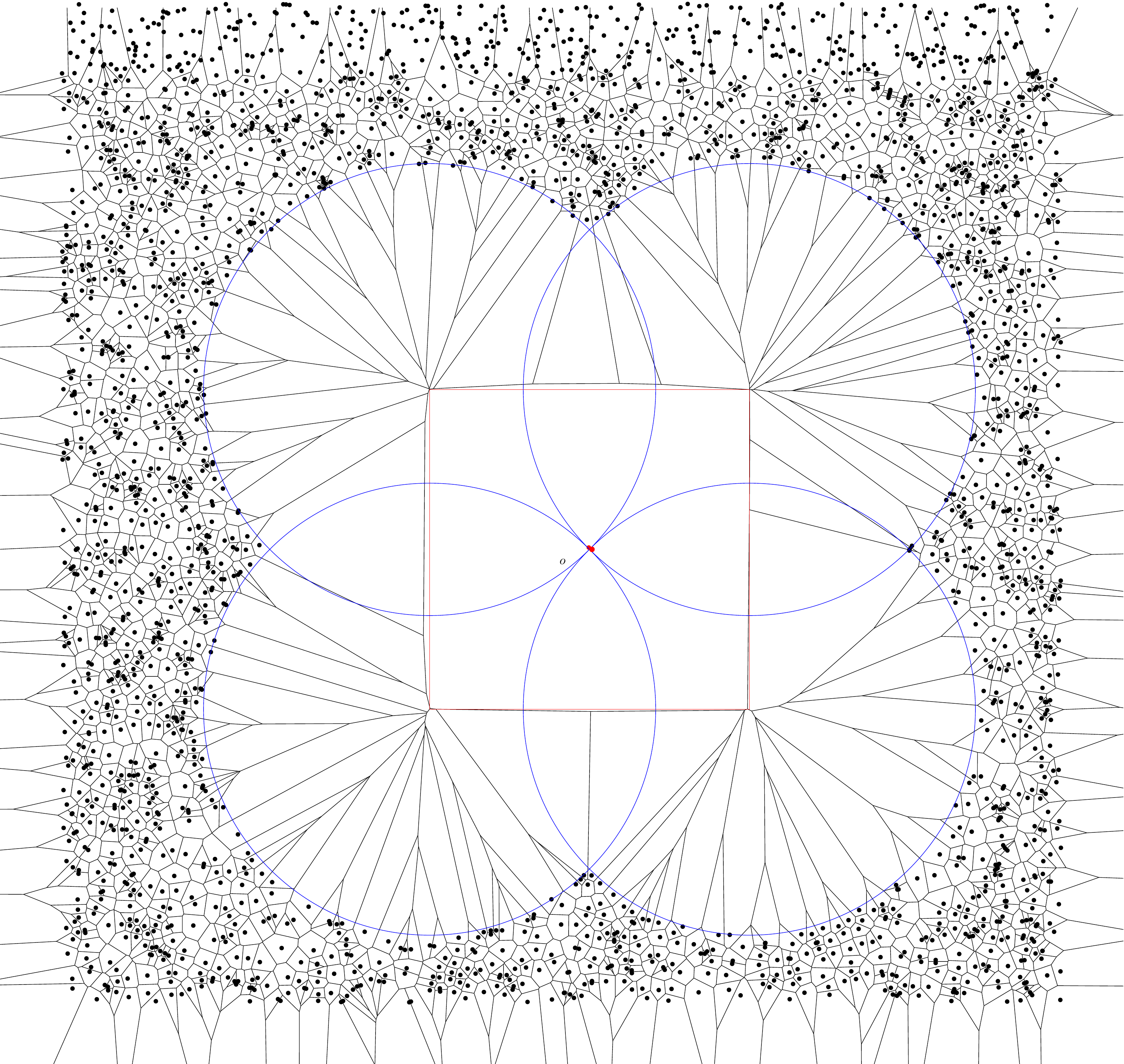}
\end{center}
\caption{Voronoi cells generated by the origin (red) and a Poisson point process of intensity $10000$ outside the Voronoi flower of an ellipse  (left) and outside the Voronoi flower of a square (right).}
\label{fig:question2}
\end{figure}

\pass Question $2$ is a special case of Question $1$ since the conditioning means that the Poisson point process has uniform intensity outside twice the Voronoi flower of $K$ with respect to $o$ (see Figure \ref{fig:question2}). Therefore, the answer to Question $2$ can be deduced from the answer to Question $1$ when the domain $\cD$ is a Voronoi flower.

\noi In order to treat Question $1$ in the case of a general domain $\cD$, we have to prove that, up to higher orders, the quantities we consider coincide with what we find for Question $2$ when $\cD$ is replaced by the maximal Voronoi flower contained in $\cD$.

\pass As in the question of the approximation of a convex body $K$ \textit{from inside} by the convex hull of a large uniform sample of points, see e.g. the two seminal papers by R\'enyi and Sulanke \cite{resu63,resu64}, answering Question $2$ is radically different when $K$ is smooth and when $K$ is polygonal. Our results on asymptotic mean values will reflect this dichotomy.

\pass Another purpose of this paper is to answer an intrinsic version of Question $2$. Indeed, Question $2$ involves an arbitrary additional point, namely the origin $o$, and in this respect this question is not intrinsic in $K$ whereas the following is.

\pass
\begin{center}
\begin{minipage}{13cm}
\textsc{Question 3:} \textit{Given a fixed convex body $K$ and a Poisson point process conditioned on its associated Voronoi tessellation to not intersect $K$, what is the geometry of the cell containing $K$ as the intensity goes to infinity?}
\end{minipage}
\end{center}

\pass The conditioning here corresponds to empty the Voronoi flower associated with the nucleus of the cell containing $K$. We will show that, as the intensity goes to infinity, this nucleus concentrates around a point whose associated flower has the smallest area. This point is known as the Steiner point of $K$.
The answers to Question 2 can then be used by choosing the Steiner point as the origin.

\pass Questions 1 and 2 fall naturally within the general literature on the asymptotic description of large cells from random tessellations. The breakthrough paper \cite{hrs04} proves and extends the famous conjecture stated by D.~G.~Kendall in the 40's and which asserts that large cells from a stationary and isotropic Poisson line tessellation are close to the circular shape.  Thereafter, the work \cite{cs05} investigates the mean defect area and mean number of vertices of the typical Poisson-Voronoi cell and of the zero-cell of a stationary and isotropic Poisson line tessellation  conditioned on containing a disk of radius $r$ when $r\to\infty$. More recently, \cite{hs14} provides an estimate of the Hausdorff distance between $K$ and its random polyhedral approximation in the slightly different model of a zero-cell from a stationary Poisson hyperplane tessellation in any dimension. To the best of our knowledge, there has not been any attempt to prove exact formulae for the means of the main geometric characteristics of a zero-cell of a Poisson-Voronoi tessellation or Poisson hyperplane tessellation around $K$. 

\pass The present work extends the results of the aforementioned paper \cite{cs05} when the disk is replaced by a general convex body $K$. The method used in \cite{cs05} relies on the application of an inversion map and of the results from \cite{resu63,resu64} on the asymptotics for random convex hulls in a disk. It is at first sight specific to the case of the disk, nevertheless we will sketch in Section \ref{sec:rsapproach} a general method for extending it for a general smooth body. It seems hopeless to try to extend it to the case of a polygon. Actually, the growth rates that we obtain when $K$ is a polygon do not even coincide with the rates from \cite{resu63,resu64}. In this paper, we have chosen an intrinsic technique which is at the same time more natural and common to the two cases of a smooth convex body and of a polygon. More precisely, our method involves notably a precise understanding of the geometry of the Voronoi flower and its interpretation as a pedal curve, a constrained version of the Blaschke-Petkantschin change of variables formula, a precise analysis of the process of the support points and a revisited Efron's identity. In particular, all these arguments turn out to be more intricate when $K$ is a polygon.

\pass Since Question $2$ is the cornerstone of all the three questions, we first focus on it.

\subsection{Answer to Question $2$: statement of key results}\label{sec:keyresults}

\pass

\noi The Euclidean plane $\R^2$ with origin $o$ is endowed with its natural scalar product $\langle \cdot,\cdot\rangle$ and the Euclidean norm $\|\cdot\|$. For any $x\in\R^2$ and $r>0$, we denote by $B_r(x)$ the closed ball centered at $x$ and of radius $r$. For two distinct points $a$ and $b$ in $\R^2$, we denote by  $(a,b)$ (resp. $[a,b]$)  the unique line containing the two points (resp. the segment between the two points).

\pass For any locally finite point set $\chi$, we construct its associated Voronoi tessellation as the collection of all cells 
\begin{equation*}
\left\{y\in\R^2:\|y-x\|\le \|y-x'\| \mbox{ for all } x'\in \chi\right\}, \,\, x\in\chi.
\end{equation*}
The set $\chi$ is called the set of nuclei associated with the tessellation. In the rest of the paper, we consider Voronoi tessellations generated by random sets of nuclei.

\pass For any $\lambda>0$, let $\cP_{\lambda}$ be a homogeneous Poisson point process on $\R^2$ with intensity $\lambda$. In particular, all processes $\cP_{\lambda}$, $\lambda>0$, may be coupled on the same space by taking $\cP_{\lambda}$ as the projection on $\R^2$ of $\cP\cap (\R^2\times (0,\lambda))$ where $\cP$ is a homogeneous Poisson point process of intensity $1$ in $\R^2\times \R$.

\pass The so-called Poisson-Voronoi tessellation is the Voronoi tessellation generated by $\cP_{\lambda}$, see e.g. \cite{mol94,okabe,cal10}. Its statistical properties are described by its typical cell which represents roughly a {\it randomly chosen} cell among the set of all cells and turns out to be distributed as the Voronoi cell of the origin $o$ with the set of nuclei $\cP_{\lambda}\cup\{o\}$ (see \cite{SW}, Theorem 3.3.5).

\pass Let $K$ be a convex body containing the origin in its interior and let $K_\lambda$ be the random convex body which is equal to the Voronoi cell
of the nucleus $o$ associated with the set of nuclei $(\cP_{\lambda}\cap (\R^2\setminus 2\cF_o(K)))\cup\{o\}$, where $\cF_o(K)$ is the Voronoi flower of $K$ with respect to $o$ (see \eqref{def:vorflow} for a precise definition). Alternatively, $K_\lambda$ is equal in distribution to the Voronoi cell  of the nucleus $o$ associated with the set of nuclei $\cP_{\lambda}\cup\{o\}$ conditional on the event that it contains $K$. In this paper, we are interested in describing the asymptotics of several characteristics of $K_{\lambda}$. More precisely, our main results provide limiting expectations up to proper rescalings of its area $\cA(K_{\lambda})$, perimeter $\cU(K_{\lambda})$ and number of vertices $\cN(K_{\lambda})$ when $K$ has a smooth boundary or is a polygon.

\pass In the sequel, we set $f(u)\underset{u}{\sim} g(u)$ (resp. $f(u)\underset{u}{=}\mathrm{O}(g(u))$) when the ratio $\tfrac{f(u)}{g(u)} \to 1$ (resp. $\tfrac{f(u)}{g(u)}$ is bounded from above) when the variable $u\to\infty$ or $u\to0$ according to the situation.

\pass We assume now that $K$ is a smooth convex body containing $o$ in its interior, that is $\partial K$ is of class $\cC^2$ with bounded positive curvature. For $s\in\partial K$ we denote by $r_s$ and $n_s$ respectively the radius of curvature and the outer unit normal vector of $\partial K$ at point $s$.

\begin{theo}[Smooth case]\label{theosmooth}
The defect area, defect perimeter and number of vertices of $K_{\lambda}$ have respectively the following asymptotics when the intensity $\lambda\to\infty:$
\begin{enumerate}
\item[$(i)$] $\E(\cA(K_{\lambda}))-\cA(K) \underset{\lambda\to\infty}{\sim}
    \lambda^{-\frac{2}{3}} 2^{-2}3^{-\frac1{3}}\Gamma\Big(\mfrac{2}{3}\Big)\displaystyle\int_{\partial K} r_s^{\frac1{3}}\left\langle s,n_s\right\rangle^{-\frac{2}{3}}\dd s$ \\
\item[$(ii)$] $\E(\cU(K_{\lambda}))-\cU(K)\underset{\lambda\to\infty}{\sim}
    \lambda^{-\frac{2}{3}} 3^{-\frac{4}{3}}\Gamma\Big(\mfrac{2}{3}\Big)\displaystyle\int_{\partial K} r_s^{-\frac{2}{3}}\langle s,n_s\rangle^{-\frac{2}{3}}\dd s$ \\
\item[$(iii)$] $\E(\cN(K_{\lambda}))\underset{\lambda\to\infty}{\sim}
    \lambda^{\frac{1}{3}}2^2 3^{-\frac{4}{3}}\Gamma\Big(\mfrac{2}{3}\Big)\displaystyle\int_{\partial K} r_s^{-\frac{2}{3}}\left\langle s,n_s\right\rangle^{\frac1{3}}\dd s.$
\end{enumerate}
\end{theo}

\pass This theorem is reminiscent of the famous results obtained by R\'enyi and Sulanke \cite{resu63,resu64} in the study of the approximation of a convex body $K$ by the convex hull $K^{\lambda}$ of $\cP_\lambda\cap K$. Their results can be summarized as follows.

\begin{theora} [Smooth case, A. R\'enyi and R. Sulanke, \cite{resu63,resu64}]\label{theo:rssmooth}
Let $K$ be a smooth convex body. Then, the defect area, defect perimeter and number of vertices of $K^{\lambda}$ have respectively the following asymptotics when the intensity $\lambda\to\infty:$
\begin{enumerate}
\item[$(i)$] $\cA(K)-\E(\cA(K^{\lambda}))\underset{\lambda\to\infty}{\sim}
     \lambda^{-\frac{2}{3}}2^{\frac{4}{3}}3^{-\frac{4}{3}}\Gamma\Big(\mfrac{2}{3}\Big)\displaystyle \int_{\partial K}r_s^{-\frac1{3}} \dd s$ \\
\item[$(ii)$] $\cU(K)-\E(\cU(K^{\lambda}))\underset{\lambda\to\infty}{\sim}
\lambda^{-\frac{2}{3}}2^{-\frac{2}{3}}3^{-\frac{1}{3}}\Gamma\Big(\mfrac{2}{3}\Big)\displaystyle\int_{\partial K}r_s^{-\frac{4}{3}} \dd s$ \\
\item[$(iii)$] $\E(\cN(K^{\lambda}))\underset{\lambda\to\infty}{\sim}
 \lambda^{\frac1{3}}2^{\frac{4}{3}}3^{-\frac{4}{3}}\Gamma\Big(\mfrac{2}{3}\Big) \displaystyle \int_{\partial K}r_s^{-\frac1{3}} \dd s.$
\end{enumerate}
\end{theora}

\pass Let us notice that in both theorems, the exponents of $\lambda$ coincide but the geometric quantities involved in the constants differ.
In particular, these quantities in R\'enyi and Sulanke's theorem above are intrinsic whereas they depend not only on $K$ but also on the origin through
the variable $s$ present in the integral terms of our results. Moreover, Efron's identity \cite{Efron} for random convex hulls connects the mean area of $K^{\lambda}$ to its mean number of vertices, which explains that the integrals over $\partial K$ in results (i) and (iii) of R\'enyi and Sulanke's theorem are the same. We will show that we have a similar pattern between points (ii) and (iii) of Theorem \ref{theosmooth}, as emphasized by Proposition \ref{prop:efron}. A deeper explanation of the strong interplay between Theorem \ref{theosmooth} and R\'enyi and Sulanke results will be given in Section \ref{sec:rsapproach}.

\pass Now let $K$ be a convex polygon with $n_K\ge 3$ consecutive vertices in anticlockwise order denoted by $a_1,\ldots,a_{n_K}$ and set $a_{n_K+1}=a_1$. Then, denoting by $o_i$ the orthogonal projection of $o$ onto the line $(a_i,a_{i+1})$
 and by $\alpha_i$ the interior angle at vertex $a_i$, we obtain similar results for $K_{\lambda}$.

\begin{theo}[Polygonal case]\label{theopolygon}
The defect area, defect perimeter and number of vertices of $K_{\lambda}$ have respectively the following asymptotics when the intensity $\lambda\to\infty:$
\begin{enumerate}
\item[$(i)$]  $\E(\cA(K_\lambda))-\cA(K)\underset{\lambda\to\infty}{\sim}
   \displaystyle \lambda^{-\frac1{2}}2^{-\frac{9}{2}}\pi^{\frac{3}{2}} \sum_{i=1}^{n_K} {\|o_i\|}^{-\frac1{2}}\|a_{i+1}-a_i\|^{\frac{3}{2}} $ \\
\item[$(ii)$] $\E(\cU(K_{\lambda}))-\cU(K)\underset{\lambda\to\infty}{\sim}
   (\lambda^{-1}\log \lambda)\cdot 2^{-1} 3^{-1}\sum_{i=1}^{n_K} \|o_i\|^{-1}  $ \\
\item[$(iii)$] $\E(\cN(K_\lambda))\underset{\lambda\to\infty}{\sim}(\log\lambda)\cdot 2\cdot 3^{-1}n_K.$
\end{enumerate}
\end{theo}

\pass A. R\'enyi and R. Sulanke have calculated the asymptotic mean number of vertices in \cite{resu63}. They calculated the mean area and mean perimeter only in the case where $K$ is the square. When $K$ is any convex polygon, the asymptotic mean area is a consequence of Efron's identity \cite{Efron} while the asymptotic mean perimeter is due to C. Buchta \cite{Buchta84} and is an explicit function of the angles $\alpha_i$.

\begin{theorb}[Polygonal case, A. R\'enyi and R. Sulanke, \cite{resu63,resu64}, C. Buchta, \cite{Buchta84}]\label{theopolygon-rs}
Let $K$ be a convex polygon. Then, the defect area, defect perimeter and number of vertices of $K^{\lambda}$ have respectively the following asymptotics when the intensity $\lambda\to\infty:$
\begin{enumerate}
\item[$(i)$] $\cA(K)-\E(\cA(K^{\lambda}))\underset{\lambda\to\infty}{\sim}(\lambda^{-1}\log \lambda)\cdot 2\cdot 3^{-1}n_K $ \\
\item[$(ii)$] $\cU(K)-\E(\cU(K^{\lambda}))\underset{\lambda\to\infty}{\sim}\lambda^{-\frac{1}{2}}\sum_{i=1}^{n_K}\psi(\alpha_i)$ \\
\item[$(iii)$] $\E(\cN(K^{\lambda}))\underset{\lambda\to\infty}{\sim}(\log\lambda)\cdot2\cdot 3^{-1}n_K.$
\end{enumerate}
where the function $\psi$ is explicit $($see Satz $1$ in \cite{Buchta84}$)$.
\end{theorb}

\pass We observe that contrary to the smooth case, the respective growth rates of the mean defect area and mean defect perimeter in Theorem \ref{theopolygon} and in R\'enyi and Sulanke's results do not coincide. Again, the constants in points (i) and (ii) of Theorem \ref{theopolygon} depend on the position of the origin inside $K$. Surprisingly, the limiting mean number of vertices in point (iii) of Theorem \ref{theopolygon} does not and even coincides with R\'enyi and Sulanke corresponding result. To the best of our knowledge, there is no easy explanation of this feature.

\pass The paper is structured as follows. We start by introducing in Section \ref{sec:keytools} the strategy and key tools for proving the announced asymptotic results. 
The proofs of Theorems \ref{theosmooth} and \ref{theopolygon} are presented in Section \ref{sec:smooth} and Section \ref{sec:polygon} for the smooth and polygonal cases respectively. 
Section \ref{sec:steiner} is devoted to Question $3$ and in particular to the convergence of the nucleus of the Voronoi cell containing $K$ to the Steiner point of $K$. In Section \ref{sec:tele}, we show how the results of Section \ref{sec:smooth} and Section \ref{sec:polygon} can be applied to answer Question $1$. Finally, we gather extensions to our work and open questions in Section \ref{sec:conclusion}.

\section{Strategy and key tools}\label{sec:keytools}

\pass In this section, we rewrite in a tractable way the three expectations which appear in Theorems \ref{theosmooth} and \ref{theopolygon}, i.e. we aim at getting the three relations \eqref{eq:airesmooth}, \eqref{eq:croftonperimeter} and \eqref{eq:efronvertices}. We also emphasize the basic ideas and guidelines of the proofs from Sections \ref{sec:smooth} and \ref{sec:polygon}.

\pass Let us introduce the Voronoi flower of a compact set $L$ with respect to a point $x\in\R^2$ as the set defined by
\begin{equation}\label{def:vorflow}
\cF_x(L)= \bigcup_{s\in L} B_{\frac1{2}\|s-x\|}\left(\tfrac1{2}(s+x)\right).
\end{equation}

\noi We notice in particular that $\cF_x(L)=\cF_x(\conv(L))$ where $\conv(\cdot)$ denotes the convex hull. The basic equivalence
\begin{equation*}\label{eq:loipoissfleur}
x\in K_\lambda\setminus K \Longleftrightarrow \cP_{\lambda}\cap 2(\cF_o(K\cup \{x\})\setminus \cF_o(K))=\emptyset
\end{equation*}
and the equality
\begin{equation}\label{eq:emptycap}
\P(\cP_{\lambda}\cap 2(\cF_o(K\cup \{x\})\setminus \cF_o(K)) = \emptyset) = \exp(-4\lambda\cA(\cF_o(K \cup \{x\})\setminus\cF_o(K)))
\end{equation}
imply that
\begin{align}\label{eq:airesmooth}
\E(\cA(K_{\lambda}))-\cA(K) = \int_{\R^2\setminus K} \exp(-4\lambda\cA(\cF_o(K \cup \{x\})\setminus\cF_o(K)))\dd x.
\end{align}

\noi Thus, the basic problem consists in providing accurate estimates for the extra area of the flower of $K$ when adding to $K$ a single point $x$ outside
of it. We will need to treat separately the case where $K$ has a smooth boundary (Lemma \ref{lem:geom}) and the case where $K$ is a convex polygon (Lemma \ref{lem:geompoly1}). In particular, Lemma \ref{lem:geom} will be proved in two different ways. One of the proofs is based on the Cauchy-Crofton formula \eqref{eq:formuleairefleur} involving the so-called support function of $K$.

\pass For every $\theta\in [0,2\pi)$, let us denote by $(u_{\theta},v_{\theta})$ the orthonormal basis in direction $\theta$, i.e. $u_{\theta}=(\cos\theta,\sin\theta)$ and $v_{\theta}=(-\sin\theta,\cos\theta)$. The support function of $K$ with respect to a point $x\in\R^2$ (see e.g. \cite{Schneider}, section 1.7) is the function defined for $z\in\R^2$ by
\begin{equation}\label{def:suppfunc}
p_x(K,z)=  \sup_{y\in K}\langle y-x,z\rangle.
\end{equation}
Observe that $p_x(K,\cdot)$ is homogeneous of degree $1$. We will denote by $p_x(K,\theta)$ the quantity $p_x(K,u_\theta)$ and we will use indifferently both notations in the sequel, depending on the context. In particular, the distance from $x$ to the boundary of $\cF_x(K)$ in direction $u_{\theta}$ is precisely $p_{x}(K,\theta)$, which implies in turn that
\begin{equation}\label{eq:formuleairefleur}
\cA(\cF_x(K))=\frac1{2}\int_0^{2\pi}p^2_x(K,\theta)\dd \theta.
\end{equation}

\pass The support function also makes it possible to rewrite the defect perimeter as an integral using the well-known Cauchy-Crofton formula
\begin{equation}\label{eq:croftonperimeter}
\E(\cU(K_{\lambda}))-\cU(K)=\int_0^{2\pi} \E(p_o(K_{\lambda},\theta)-p_o(K,\theta))\dd\theta.
\end{equation}

\noi Therefore, in order to deal with this expectation we aim to determine the distribution of the point which achieves the support function into a fixed direction (see Propositions \ref{prop:supportsmooth} and \ref{prop:supportpoly}). The strategy will consist in using nontrivial changes of variable according to the case where $K$ is smooth or not (see Lemma \ref{lem:jacob}). The main difficulty will be in the computation of its Jacobian, the determination of the domain of integration for it and finally its integration.

\pass In the next proposition, we prove a relation in the same spirit as the well-known Efron's relation for convex hulls of random inputs, see e.g. \cite{Efron}, which connects the mean number of sides of $K_\lambda$ either to the mean defect area of the flower or to the mean defect support function.

\begin{prop}\label{prop:efron}
\noi
\begin{enumerate}
\item[$(i)$] For every $\lambda>0$, the following identity holds
\begin{equation}
\E(\cN(K_{\lambda})) = 4\lambda (\E(\cA(\cF_o(K_{\lambda}))-\cA(\cF_o(K))).
\end{equation}
\item[$(ii)$] Moreover, when $\lambda\to\infty$,
\begin{equation}\label{eq:efronvertices}
\E(\cN(K_{\lambda})) \underset{\lambda\to\infty}{\sim} 4\lambda\int_0^{2\pi}p_o(K,\theta)\E(p_o(K_{\lambda},\theta)-p_o(K,\theta))\dd\theta.
\end{equation}
\end{enumerate}
\end{prop}

\pass{\bf Proof of Proposition \ref{prop:efron}.}
We recall that $\cN(K_{\lambda})$ is the number of neighbors of $o$, i.e. the set of all  $x\in\cP_{\lambda}\setminus 2\cF_o(K)$ such that the bisecting line of the segment $[o,x]$ has a non-empty intersection with the boundary of $K_\lambda$. Moreover, for any $x\in\cP_{\lambda}\setminus 2\cF_o(K)$, $x$ is a neighbor of $o$ if and only if $\mfrac{1}{2}x\in \cF_o(V_x)\setminus\cF_o(K)$, where $V_x$ is the Voronoi cell of the origin associated with the set of nuclei $(\cP_{\lambda}\setminus 2\cF_o(K))\setminus\{x\}$. Consequently, thanks to Mecke-Slivnyak's formula (see Corollary $3.2.3$ in \cite{SW}) and Fubini theorem, we obtain
\begin{align*}
\E(\cN(K_{\lambda})) & = \E\bigg(\sum_{x\in \cP_{\lambda}\setminus 2\cF_o(K)}{\emph{\indicat}}_{\{\frac{x}{2}\in \cF_o(V_x)\setminus \cF_o(K)\}}\bigg) \\
& = \lambda\int_{\R^2\setminus 2\cF_o(K)} \P(x\in 2(\cF_o(K_\lambda)\setminus\cF_o(K)))\dd x \\
& = 4\lambda\E(\cA(\cF_o(K_{\lambda}))-\cA(\cF_o(K))).
\end{align*}

\noi Now, using \eqref{eq:formuleairefleur}, we obtain
\begin{align*}\label{sophie}
\E(\cN(K_{\lambda})) & = 2\lambda\int_0^{2\pi}\E(p_o^2(K_{\lambda},\theta)-p_o^2(K,\theta))\dd\theta\nonumber\\
& = 2\lambda\int_0^{2\pi}\E\left((p_o(K_{\lambda},\theta)+p_o(K,\theta))(p_o(K_{\lambda},\theta)-p_o(K,\theta))\right)\dd\theta.
\end{align*}
Let us denote by $\dd_H(K_\lambda,K)$ the Hausdorff distance between $K_\lambda$ and $K$. We get from the equality above and from the inequality $p_o(K_\lambda,\theta)\le p_o(K,\theta)+\dd_H(K_\lambda,K)$ that
\begin{equation*}\label{eq:hausdist}
 0 \le \E(\cN(K_{\lambda}))-4\lambda\int_0^{2\pi}p_o(K,\theta)\E(p_o(K_{\lambda},\theta)-p_o(K,\theta))\dd\theta\le 4\pi\lambda\E(\dd_H^2(K_\lambda,K)).
\end{equation*}

\noi Here we will use the fact that the methods developed in \cite{hs14} in the case of a zero-cell of a stationary Poisson hyperplane tessellation may be used in our setting  to show that the distance $\dd_H(K_\lambda,K)$ decreases to zero almost surely. This means in particular that $\E(\dd_H^2(K_\lambda,K))\to 0$ when $\lambda\to\infty$, which completes the proof. $\square$

\pass Let us notice that Proposition \ref{prop:efron} can be extended to higher dimension when the number of sides is replaced by the number of facets.

\section{Answer to Question $2$: proof for the smooth case}\label{sec:smooth}

\noi In this section, $K$ is a smooth convex body containing the origin in its interior. Every $x\in\R^2\setminus K$ can be written as $x=s+hn_s=s_h$ with $s\in\partial K$ and $h>0$. We denote by $\Delta\cF_{s_h}$ the set $\cF_o(K\cup \{s_h\})\setminus \cF_o(K)$. Because of \eqref{eq:emptycap} and \eqref{eq:airesmooth}, a key ingredient for proving Theorem \ref{theosmooth} is the estimate of the area of the increase $\Delta \cF_{s_h}$ of a Voronoi flower of $K$ induced by the addition of a point outside $K$. To the best of our knowledge this estimate is new despite the natural aspect of the question. Our method is based on original considerations on the curvature of the boundary of the Voronoi flower of $K$ and this is done in Subsection \ref{subsec:fleursmooth}. Subsections \ref{subsec:surfacesmooth}, \ref{sub:perimetersmooth} and \ref{subsec:verticessmooth} are then devoted to the asymptotic mean area of $K_\lambda$, the asymptotic mean support function and perimeter of $K_\lambda$ and the asymptotic intensity and mean number of vertices of $K_\lambda$ respectively.

\subsection{Increase of the area of the Voronoi flower}\label{subsec:fleursmooth}

\pass

\noi The next lemma provides the exact calculation of the limiting rescaled defect area of the Voronoi flower as well as a lower-bound.

\begin{lem}\label{lem:geom}
Let us assume that $K$ is a smooth convex body containing the origin $o$ in its interior.
\begin{enumerate}
\item[$(i)$] For every $s\in\partial K$, we get
\begin{equation*}\label{estimeevolfleur}
\cA(\Delta\cF_{s_h})\underset{h\to 0}{\sim} h^{\frac{3}{2}} 2^{\frac{5}{2}}3^{-1}r_s^{-\frac1{2}}\langle s, n_s\rangle.
\end{equation*}
\item[$(ii)$] Moreover, there exists $C>0$ such that, for every $h>0$ and $s\in\partial K$,
\begin{equation*}\label{eq:minorantvol}
h^{-\frac{3}{2}}\cA(\Delta\cF_{s_h}) \ge C >0 .
\end{equation*}
\end{enumerate}
\end{lem}

\pass{\bf Proof of Lemma \ref{lem:geom}.}
Actually we will provide two different proofs of assertion (i). The first proof is possibly the most natural one but specific to the planar case. The second proof is more analytical and can be extended to any dimension.

\pass{\it First proof of $(i)$}. Our first method is based on a precise geometric description of the increase of the Voronoi flower $\cF_o(K)$ and on the
observation that the boundary of $\cF_o(K)$ is nothing but the pedal curve of $\partial K$, i.e. the set of orthogonal projections of $o$ onto the tangent lines of $\partial K$ (see \cite{yates52}, p. 160).

\begin{figure}[!h]
\begin{center}
\includegraphics[scale=0.75]{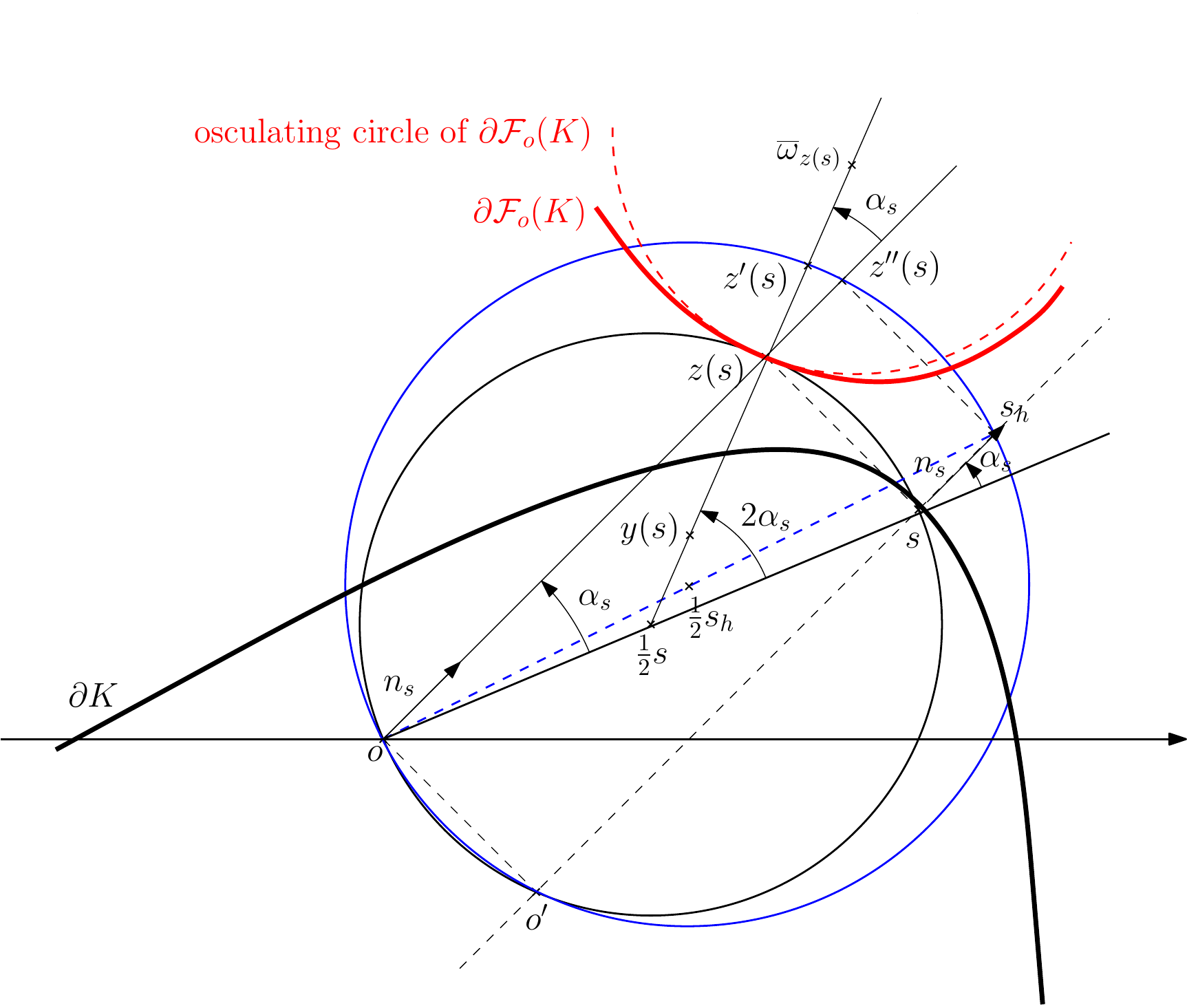}
\vspace{-0.2cm}
\caption{\label{fig:geolem-proof1} The influence on the Voronoi flower of $K$ of adding an extra point $s_h$ in a neighborhood of its boundary viewed from the origin $o$.}
\end{center}
\end{figure}

\pass Let us introduce some notation (see Figure \ref{fig:geolem-proof1}):
\begin{enumerate}[label=-,leftmargin=2\parindent]
\item The point $z(s)$ is the unique point of intersection of the pedal curve $\partial\cF_o(K)$ with the ball $B_{\frac1{2}\|s\|}(\tfrac1{2}s)$
(the uniqueness of $z(s)$ comes from the fact that $K$ is a smooth convex body),
\item The point $z'(s)$ is the intersection of the half-line $\tfrac1{2}s +\R_+(z(s)-\tfrac1{2}s)$ with the circle $\partial B_{\frac1{2}\|s_{h}\|}(\tfrac1{2}s_h)$,
\item The point $z''(s)$ is the intersection of the line $o+\R z(s)$ with the circle $\partial B_{\frac1{2}\|s_{h}\|}(\tfrac1{2}s_h)$,
\item The point $\overline{\omega}_{z(s)}$ is the center of curvature of $\partial\cF_o(K)$ at $z(s)$,
\item The point $y(s)$ is the intersection of the half-line $\tfrac1{2}z'(s)+\R_+(\tfrac1{2}s-z'(s))$ with the circle $\partial B_{\frac1{2}\|s_h\|}(z'(s))$. 
\end{enumerate}

\pass Let us emphasize a few geometric observations:
\begin{enumerate}[label=-,leftmargin=2\parindent]
\item The points $\tfrac1{2}s$, $y(s)$, $z(s)$, $z'(s)$ and $\overline{\omega}_{z(s)}$ are aligned since the pedal curve and the ball
$B_{\frac1{2}\|s\|}(\tfrac1{2}s)$ have the same normal vector at $z(s)$.
\item The line containing $o$, $z(s)$ and $z''(s)$ has the same direction as $n_s$. Indeed, let $o'$ be the intersection point of
$\partial B_{\frac1{2}\|s\|}(\tfrac1{2}s)$ and $\partial  B_{\frac1{2}\|s_h\|}(\tfrac1{2}s_h)$ different from $o$.
The three triangles $oo's$, $oz''(s)s_h$ and $o'oz(s)$ are inscribed in half-balls so that they are right triangles and $oo's_hz''(s)$ is a rectangle).
\item The points $z(s)$ and $z''(s)$ are the respective symmetric points of $s$ and $s_h$ with respect to the line with direction
$n_s$ containing $\tfrac1{2}s$ and $\tfrac1{2}s_h$.
\item The angle $\alpha_s$ between the two half-lines $\R_+ n_s$ and $\tfrac1{2}s+\R_+(z(s)-\tfrac1{2}s)$ satisfies
\begin{equation}\label{eq:alphas}
\cos\alpha_s=\left\langle \frac{s}{\|s\|},n_s\right\rangle
\end{equation}
by symmetry with respect to the same line.
\end{enumerate}

\pass Consequently, we have the equalities $\|z(s)-z''(s)\|=  \|s_h-s\| = h$ and
\begin{equation}\label{eq:hh'}
\|z'(s)-z(s)\|\underset{h\to 0}{=} h\cos\alpha_s+\mathrm{O}(h^2)
\end{equation}

\noi We denote by $\rho_{z(s)}=\|z(s)-{\overline{\omega}}_{z(s)}\|$ the radius of curvature of the pedal curve $\partial \cF_o(K)$
at the point $z(s)$. The strategy is to approximate the set $\Delta\cF_{s_h}$ by the simpler set $\Delta_{s_h}$ defined by
\begin{equation*}
\Delta_{s_h}
=  \left\{\begin{array}{ll}
      B_{\frac1{2}\|s_h\|}(y(s))\cap B_{\rho_{z(s)}}({\overline{\omega}}_{z(s)}) & \hbox{if the curvature of $\partial\cF_o(K)$ at $z(s)$ is negative} \\
      B_{\frac1{2}\|s_h\|}(y(s))\setminus B_{\rho_{z(s)}}({\overline{\omega}}_{z(s)}) & \hbox{if the curvature of $\partial\cF_o(K)$ at $z(s)$ is positive}
          \end{array}\right.
\end{equation*}
where the two balls are both centered on the line containing $\tfrac1{2}s$ and $z(s)$.

\pass In the sequel we only deal with the case if the curvature of $\partial\cF_o(K)$ at $z(s)$ is negative, the other case will be treated similarly. Let $\beta_{z(s)}$ (resp. $\gamma_{z'(s)}$) be the aperture of the arc of circle $\partial B_{\rho_{z(s)}}(\overline{\omega}_{z(s)}) \cap \Delta_{s_h}$ (resp. $\partial B_{\frac1{2}\|s_{h}\|}(y(s)) \cap \Delta_{s_h}$). We assume without loss of generality that $\beta_{z(s)}$ has same sign as $\rho_{z(s)}$ and $\gamma_{z'(s)}$ is positive. Notice that, since $\|z'(s)-z(s)\|$ is proportional to $h$, we get
\begin{equation}\label{eq:anglesouverture}
\beta_{z(s)}\underset{h\to 0}{=}\mathrm{O}(h^{\frac1{2}})\,\text{ and }\,\gamma_{z'(s)}\underset{h\to 0}{=}\mathrm{O}(h^{\frac1{2}}).
\end{equation}

\pass First we will evaluate the error between the two areas $\cA(\Delta\cF_{s_h})$ and $\cA(\Delta_{s_h})$.
Observe that this error is due to two different contributions: the replacement of $\cF_o(K)$ by $B_{\rho_{z(s)}}(\overline{\omega}_{z(s)})$ and
the replacement of $B_{\frac1{2}\|s_{h}\|}(\tfrac1{2}s_{h})$ by $B_{\frac1{2}\|s_{h}\|}(y(s))$. Let us examine each one separately.

\pass On one hand, the region $\delta_{s_h}^{(1)}$ between $\partial\cF_o(K)$ and its osculating circle at $z(s)$ and lying inside the ball
$B_{\frac1{2}\|s_h\|}(\tfrac1{2}s_h)$ has an area of order $\mathrm{O}(h^2)$. Indeed, let us consider the half-line with origin $\overline{\omega}_{z(s)}$
and making the angle $\theta$ with the half-line $\overline{\omega}_{z(s)}+\R_+(z(s)-\overline{\omega}_{z(s)})$. A Taylor expansion shows that the distance
between the intersection points of that half-line with $\partial\cF_o(K)$ and $\partial B_{\rho_{z(s)}}(\overline{\omega}_{z(s)})$ respectively is
of order $\mathrm{O}(\theta^3)$. Moreover, because of \eqref{eq:anglesouverture}, the point $\overline{\omega}_{z(s)}$ sees the intersection
$B_{\rho_{z(s)}}(\overline{\omega}_{z(s)})\cap B_{\frac1{2}\|s_h\|}(\tfrac1{2}s_h)$ with an angle of order $\mathrm{O}(h^{\frac1{2}})$.
Consequently, integrating $\theta^3$ for angles between $0$ and $h^{\frac1{2}}$, we obtain that the area of $\delta_{s_h}^{(1)}$ is of order $\mathrm{O}(h^2)$.

\pass On the other hand, the region $\delta_{s_h}^{(2)}$ between the two balls $B_{\frac1{2}\|s_h\|}(\tfrac1{2}s_h)$ and $B_{\frac1{2}\|s_h\|}(y(s))$
and lying inside the ball $B_{\rho_{z(s)}}(\overline{\omega}_{z(s)})$ has also an area of order $\mathrm{O}(h^2)$.
Indeed, these balls cross at point $z'(s)$ with an angle of order $\mathrm{O}(h)$, i.e. the order of magnitude of the distance between $y(s)$ and $\tfrac1{2}s_h$. By \eqref{eq:anglesouverture}, the intersection  $B_{\rho_{z(s)}}(\overline{\omega}_{z(s)})\cap B_{\frac1{2}\|s_h\|}(y(s))$
has an aperture of order $h^{\frac1{2}}$. Consequently, the region $\delta_{s_h}^{(2)}$ has therefore an area of the same order
as the area of an isocele triangle of height $h^{\frac1{2}}$ and aperture $h$, i.e. of order $\mathrm{O}((h^{\frac1{2}})^2 h) =\mathrm{O}(h^2)$.

\pass Consequently, since $|\cA(\Delta\cF_{s_h})-\cA(\Delta_{s_h})|\le \cA(\delta_{s_h}^{(1)})+\cA(\delta_{s_h}^{(2)})$,
we obtain the following approximation
\begin{equation}\label{eq:voli}
\cA(\Delta\cF_{s_h})\underset{h\to 0}{=}\cA(\Delta_{s_h})+ \mathrm{O}(h^2).
\end{equation}

\pass Now it remains to evaluate the area $\cA(\Delta_{s_h})$. To do this, we need to provide precise estimates of the two angles $\beta_{z(s)}$ and $\gamma_{z'(s)}$.

\pass We can write the half-diameter of $\Delta_{s_h}$ as
\begin{equation} \label{eq:proofheight}
\rho_{z(s)} \sin\Big(\mfrac1{2}\beta_{z(s)}\Big)=\frac1{2}\|s_{h}\| \sin\Big(\mfrac1{2}\gamma_{z'(s)}\Big)
\end{equation}
and its width
\begin{equation} \label{eq:proofthick}
\|z'(s)-z(s)\|= \frac1{2}\|s_{h}\|\left(1-\cos\Big(\mfrac1{2}\gamma_{z'(s)}\Big)\right)-\rho_{z(s)}\left(1-\cos\Big(\mfrac1{2}\beta_{z(s)}\Big)\right),
\end{equation}
so that its area is finally given by
\begin{equation} \label{eq:proofvol}
\cA(\Delta_{s_h})=\frac1{8}\|s_{h}\|^2\big(\gamma_{z'(s)}-\sin(\gamma_{z'(s)})\big)-\frac1{2}\rho_{z(s)}^2\big(\beta_{z(s)}-\sin(\beta_{z(s)})\big).
\end{equation}

\pass Let us notice that \eqref{eq:anglesouverture} and \eqref{eq:proofheight} imply that
\begin{equation}\label{eq:betagamma}
\gamma_{z'(s)}\underset{h\to 0}{=}2\rho_{z(s)}\beta_{z(s)}\|s_h\|^{-1}+\mathrm{O}(h^3).
\end{equation}

\pass Inserting equalities \eqref{eq:hh'} and \eqref{eq:betagamma} in \eqref{eq:proofthick} yields
\begin{equation}\label{eq:betacarre}
\beta_{z(s)}^2\underset{h\to 0}{=}\frac{h\cos \alpha_s}{\frac1{4}\rho_{z(s)}^2\|s_h\|^{-1}-\frac1{8}\rho_{z(s)}}+\mathrm{O}(h^2).
\end{equation}

\pass We then need to calculate the radius of curvature $\rho_{z(s)}$ of the pedal curve $\partial \cF_o(K)$ at point $z(s)$.
It is known (see \cite{steede}) that it is given by
\begin{equation}\label{eq:curvpodaire}
\rho_{z(s)} = \frac{\|s\|^2}{2\|s\| - r_s\cos\alpha_s}.
\end{equation}

\pass Keeping in mind that $\beta_{z(s)}$ and $\rho_{z(s)}$ have the same sign, and using
\begin{equation}\label{eq:estimsh}
\|s_h\|\underset{h\to 0}{=}\|s\|+\mathrm{O}(h),
\end{equation}
we deduce from \eqref{eq:curvpodaire} and \eqref{eq:betacarre} that
\begin{equation}\label{eq:proofalphaOmega}
\beta_{z(s)}\underset{h\to 0}{=} 2^{\frac{3}{2}} r_s^{-\frac1{2}} \left(2\|s\|-r_s\cos\alpha_s\right)\|s\|^{-1}h^{\frac1{2}} + \mathrm{O}(h^{\frac{3}{2}}).
\end{equation}

\pass We deduce from \eqref{eq:betagamma} and \eqref{eq:proofalphaOmega} that
\begin{equation}\label{eq:proofalphaG}
\gamma_{z'(s)} \underset{h\to 0}{=} 2^{\frac{5}{2}}r_s^{-\frac1{2}}h^{\frac1{2}}+\mathrm{O}(h^{\frac{3}{2}}).
\end{equation}

\pass Thus, inserting \eqref{eq:curvpodaire}, \eqref{eq:estimsh}, \eqref{eq:proofalphaOmega} and \eqref{eq:proofalphaG} into \eqref{eq:proofvol} we obtain
\begin{equation}\label{eq:proofvoli}
\cA(\Delta_{s_h})\underset{h\to 0}{=} 2^{\frac{5}{2}}3^{-1}r_s^{-\frac1{2}}\|s\|(\cos\alpha_s) h^{\frac{3}{2}}+ \mathrm{O}(h^{\frac{5}{2}}).
\end{equation}
Finally, we obtain (i) by combining  \eqref{eq:voli}, \eqref{eq:proofvoli} and \eqref{eq:alphas}.

\pass

\pass{\it Second proof of $(ii)$.}
The second proof goes along the following lines.
We wish to use the simpler case where the origin $o$ coincides with the center of curvature $\omega_s$ of $\partial K$ at point $s$. In other words, our aim is to show that the area $\cA(\Delta\cF_{s_h})$ can be calculated in function of $\cA(\Delta\widetilde{\cF}_{s_h})$ where
\begin{equation*}
\Delta\widetilde{\cF}_{s_h}=\cF_{\omega_s}(K \cup \{s+h n_s\})-\cF_{\omega_s}(K).
\end{equation*}

\pass We use the equality \eqref{eq:formuleairefleur} and the relation, for every $x\in\R^2$,
\begin{equation}\label{eq:suppchangeref}
p_x(K,\theta) -p_o(K,\theta) = -\langle x, u_{\theta}\rangle
\end{equation}
to obtain
\begin{equation*}
\cA(\cF_o(K)) =\cA(\cF_{\omega_s}(K))-\frac1{2}\int_0^{2\pi}\langle \omega_s, u_{\theta}\rangle^2\dd \theta+\int_0^{2\pi}p_{\omega_s}(K,\theta)\langle \omega_s, u_{\theta}\rangle\dd \theta.
\end{equation*}

\pass Applying this formula to both $K$ and $K\cup\{s_h\}$ yields
\begin{align}\label{eq:chaslesvol}
\cA(\Delta\cF_{s_h})=\cA(\Delta\widetilde{\cF}_{s_h}) + \int_0^{2\pi} \Delta p_{\omega_s}(\theta)\langle \omega_s, u_{\theta}\rangle\dd\theta
\end{align}
where
\begin{align*}\label{eq:chaslesvol2}
\Delta p_{\omega_s}(\theta)=p_{\omega_s}(K \cup \{s_h\},\theta)-p_{\omega_s}(K, \theta).
\end{align*}

\begin{figure}[!h]
\begin{center}
\includegraphics[scale=0.75]{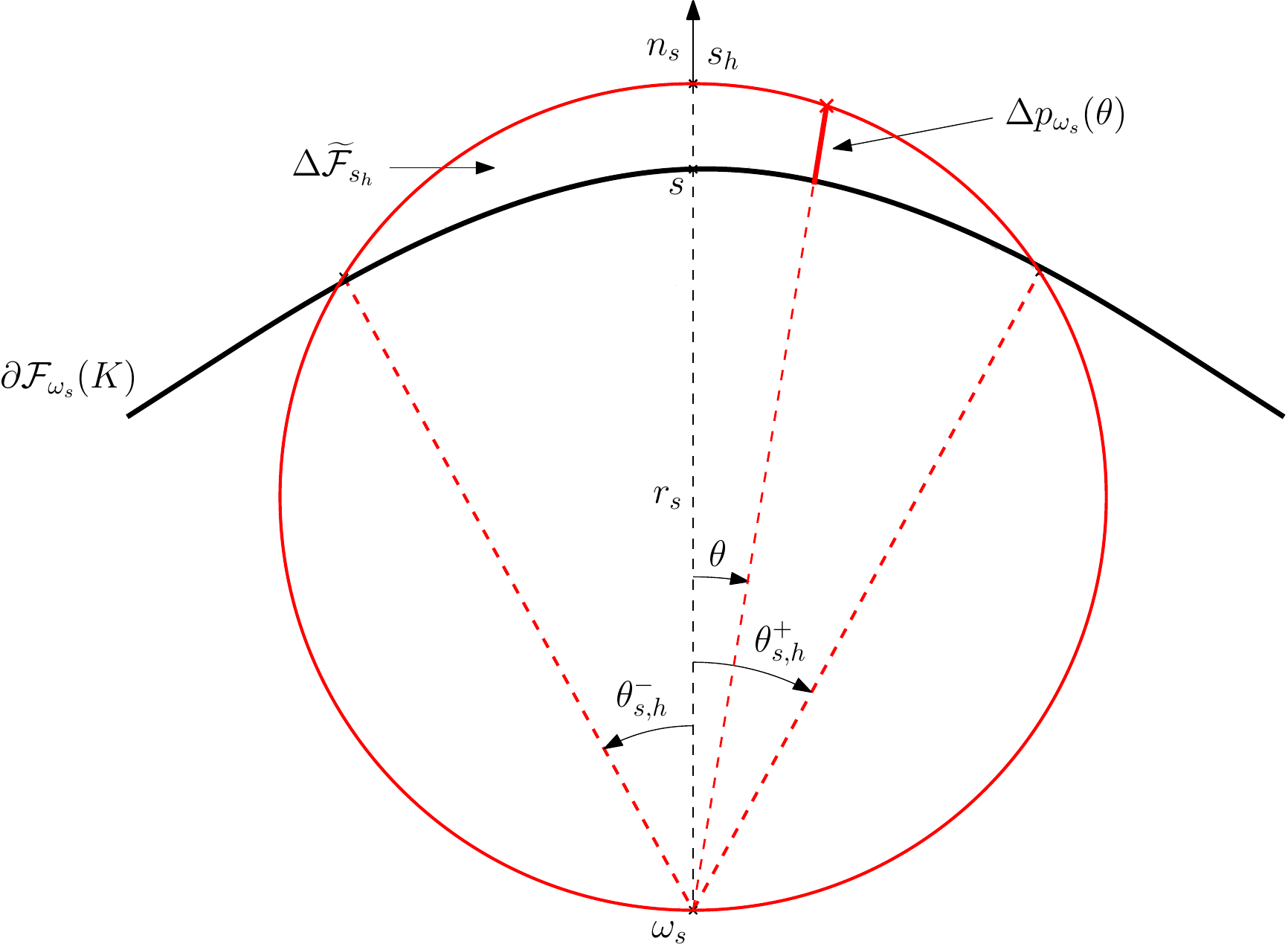}
\vspace{-0.2cm}
\caption{\label{fig:geolem-proof2} Flowers viewed from the center of curvature $\omega_s$ of $\partial K$ at point $s$.}
\end{center}
\end{figure}

\noi We treat separately the two terms of the right-hand side of \eqref{eq:chaslesvol} starting with $\cA(\Delta\widetilde{\cF}_{s_h})$.
Let us define $s(\theta)$ as the point belonging to $\partial K$ such that $\theta \in(-\pi,\pi]$ is the angle between the two half-lines
$\omega_s+\R_+(s-\omega_s)$ and $\omega_s+\R_+(s(\theta)-\omega_s)$ (see Figure \ref{fig:geolem-proof2}). Denote by $-\theta_{s,h}^-$ and $\theta_{s,h}^+$
the two angles such that $p_{\omega_s}(K\cup\{s_h\},\theta) = p_{\omega_s}(K,\theta)$ if and only if $\theta\notin[-\theta_{s,h}^-,\theta_{s,h}^+]$. Then we can write
\begin{equation}\label{eq:voldifflow1}
\cA(\Delta\widetilde{\cF}_{s_h})= \int_{-\theta_{s,h}^-}^{\theta_{s,h}^+}\int_{p_{\omega_s}(K,\theta)}^{(r_s+h)\cos\theta} r \dd r \dd\theta.
\end{equation}

\noi If $\partial K$ was a perfect circle of radius $r_s$ in the neighborhood of $s$ we would have $\theta_{s,h}^-=\theta_{s,h}^+=\theta_{s,h}$ with $\cos(\theta_{s,h})=\mfrac{r_s}{r_s+h}$. In the general case, we can sandwich $\partial K$ between two circles of radii $r_s+Ch$ and $r_s-Ch$ for a certain $C>0$. Consequently, angles $\theta_{s,h}^-$, $\theta_{s,h}^+$ and $\theta_{s,h}$ can all be written as
\begin{equation}\label{eq:angleaperture}
\theta_{s,h}^+\underset{h\to 0}{=}\theta_{s,h}^-+\mathrm{O}(h^{\frac{3}{2}})\underset{h\to 0}{=}\theta_{s,h}+\mathrm{O}(h^{\frac{3}{2}})
\underset{h\to 0}{=}2^{\frac1{2}}r_s^{-\frac{1}{2}}h^{\frac1{2}}+\mathrm{O}(h^{\frac{3}{2}}).
\end{equation}
Moreover, since $\theta_{s,h}^+$ is bounded from above by its value obtained when $\partial K$ is replaced by an outer circle of radius $r_s+Ch$, we have
\begin{equation}\label{eq:dominationangleouverture}
\theta_{s,h}^+\le \arccos\Big(1-\frac{h}{r_s+(C+1)h}\Big)\le 2^{\frac1{2}}r_s^{-\frac1{2}}h^{\frac1{2}}.
\end{equation}
Similarly, the defect of circularity implies that
\begin{equation}\label{eq:defectcircularity}
p_{\omega_s}(K,\theta) \underset{h\to 0}{=} r_s+\mathrm{O}(h^{\frac{3}{2}}).
\end{equation}

\pass Therefore, inserting \eqref{eq:angleaperture} and \eqref{eq:defectcircularity} in \eqref{eq:voldifflow1} yields
\begin{equation*}\label{eq:voldifflow2}
\cA(\Delta\widetilde{\cF}_{s_h})\underset{h\to 0}{=}2^{\frac{5}{2}}3^{-1}  r_s^{\frac1{2}}h^{\frac{3}{2}} + \mathrm{O}(h^2).
\end{equation*}

\pass We treat now the integral term in the right-hand side of \eqref{eq:chaslesvol}. For every $\theta\in[-\theta_{s,h}^-,\theta_{s,h}^+]$, we have
\begin{equation*}
\langle \omega_s, u_{\theta}\rangle \underset{h\to 0}{=} \langle \omega_s, n_s\rangle + \mathrm{O}(h).
\end{equation*}

\pass Moreover
\begin{equation*}
p_{\omega_s}(K\cup\{s_h\},\theta)= (r_s+h)\cos\theta.
\end{equation*}

\pass Therefore, we can write successively:
\begin{align*}
\int_0^{2\pi} \Delta p_{\omega_s}(\theta)\langle \omega_s, u_{\theta}\rangle\dd\theta
& \hspace{0.15cm} = \int_{-\theta_{s,h}^-}^{\theta_{s,h}^+} \Delta p_{\omega_s}(\theta)\langle \omega_s, u_{\theta}\rangle\dd\theta\\
& \underset{h\to 0}{=} \langle \omega_s, n_s\rangle\int_{-\theta_{s,h}^-}^{\theta_{s,h}^+} \Delta p_{\omega_s}(\theta)\dd\theta + \mathrm{O}(h^2)\\
& \underset{h\to 0}{=} \langle \omega_s, n_s\rangle\int_{-\theta_{s,h}^-}^{\theta_{s,h}^+} ((r_s+h)\cos\theta-r_s)\dd\theta + \mathrm{O}(h^2) \\
& \underset{h\to 0}{=} \langle \omega_s, n_s\rangle ((r_s+h)(\sin(\theta_{s,h}^+)+\sin(\theta_{s,h}^-))-r_s(\theta_{s,h}^++\theta_{s,h}^-)) + \mathrm{O}(h^2) \\
& \underset{h\to 0}{=} 2^{\frac{5}{2}}3^{-1}\,  r_s^{-\frac1{2}} \langle \omega_s, n_s\rangle h^{\frac{3}{2}}+ \mathrm{O}(h^2).
\end{align*}

\pass Finally, using $s=\omega_s+r_s n_s$, we get
\begin{align*}
\cA(\Delta\cF_{s_h}) & \underset{h\to 0}{=} 2^{\frac{5}{2}}3^{-1} r_s^{\frac1{2}}(1+r_s^{-1}\langle \omega_s, n_s\rangle)h^{\frac{3}{2}} + \mathrm{O}(h^{2})
\underset{h\to 0}{=} 2^{\frac{5}{2}}3^{-1}  r_s^{-\frac1{2}}\langle s,n_s\rangle h^{\frac{3}{2}} + \mathrm{O}(h^{2})
\end{align*}
that gives again the desired result.

\pass{\it Proof of $(ii)$.}
Thanks to (i) we can fix $\eps>0$ such that $\cA(\Delta\cF_{s_h})h^{-\frac{3}{2}}$ is bounded from below by a constant $C>0$
for all $h\in (0,\eps)$. When $h\geq\eps$, we notice that the region $\Delta\cF_{s_h}$ contains a disk of radius proportional to $h$,
which means that there exists $C'>0$ such that $\cA(\Delta\cF_{s_h})\geq C' h^2\geq C' \eps^{\frac1{2}}h^{\frac{3}{2}}$. $\square$

\subsection{Proof of Theorem \ref{theosmooth} (i): the defect area}\label{subsec:surfacesmooth}

\pass

\noi  Every $x\in\R^2\setminus K$ can be written as $x=s+\lambda^{-\frac{2}{3}}hn_s=s_{\lambda^{-\frac{2}{3}}h}$ with $s\in\partial K$ and $h>0$, the Jacobian of this change of variables being given by
\begin{equation*}\label{eq:jacobcurviligne}
\frac{\dd x}{\dd s \dd h} = \lambda^{-\frac{2}{3}}\big|1+\lambda^{-\frac{2}{3}}hr_s^{-1}\big|.
\end{equation*}

\noi Thus we get from \eqref{eq:airesmooth} that
\begin{align*}
\E(\cA(K_{\lambda}))-\cA(K) = \lambda^{-\frac{2}{3}}\int_{\partial K}\int_0^{\infty} \exp\big(-4\lambda \cA(\Delta\cF_{s_{\lambda^{-2/3}h}})\big)\big|1+\lambda^{-\frac{2}{3}}hr_s^{-1}\big| \dd s \dd h.
\end{align*}

\pass Thanks to Lemma \ref{lem:geom}, we get, for $h>0$ fixed,
\begin{align*}
4\lambda \cA(\Delta\cF_{s_{\lambda^{-2/3}h}}) \underset{\lambda\to\infty}{\sim} 4\big(\lambda^{-2/3}h\big)^{-\frac{3}{2}} \cA(\Delta\cF_{s_{\lambda^{-2/3}h}})  h^{\frac{3}{2}} \underset{\lambda\to\infty}{\sim} C_s h^{\frac{3}{2}}
\end{align*}
where $C_s = 2^{\frac{9}{2}}3^{-1}r_s^{-\frac1{2}}\langle s,n_s\rangle $ and the existence of $C>0$ such that, for all $\lambda>0$ and $s\in \partial K$,
\begin{equation*}
4\lambda \cA(\Delta\cF_{s_{\lambda^{-2/3}h}}) \ge C C_s h^{\frac{3}{2}}.
\end{equation*}

\pass Consequently, we can apply Lebesgue's dominated convergence theorem to obtain,
\begin{align*}
\lambda^{\frac{2}{3}}(\E(\cA(K_{\lambda}))-\cA(K))
   & \underset{\lambda\to\infty}{\sim} \int_{\partial K}\int_0^{\infty} \exp\big(-C_s h^{\frac{3}{2}}\big) \dd s\dd h \\
   & \hspace{0.22cm} = \mfrac{2}{3}\int_{\partial K}C_s^{-\frac{2}{3}}\bigg(\int_0^{\infty} l^{-\frac1{3}}\exp(-l) \dd l\biggl)\dd s \\
   & \hspace{0.22cm} = \mfrac{2}{3}\Gamma\left(\mfrac{2}{3}\right)\int_{\partial K} C_s^{-\frac{2}{3}}\dd s
\end{align*}
which provides assertion (i) of Theorem \ref{theosmooth}. $\square$

\subsection{Proof of Theorem \ref{theosmooth} (ii): support points and defect perimeter}\label{sub:perimetersmooth}

\pass

\noi We start by rewriting \eqref{eq:croftonperimeter} in the special case where $K$ is smooth. Noticing that, for every $s\in \partial K$ such that $n_s=u_\theta$, we get
\begin{equation}\label{eq:chgtsr}
\frac{\dd s}{\dd \theta}=r_s\,\text{ and }\, p_o(K,\theta)=p_o(K,n_s)=\langle s,n_s\rangle,
\end{equation}
we obtain
\begin{equation*}\label{eq:nouvellejolieformuledeyann}
\E(\cU(K_{\lambda}))-\cU(K)=\int_{\partial K} \E(p_o(K_{\lambda},n_s)-p_o(K,n_s))r_s^{-1} \dd s.
\end{equation*}
Using point (ii) of Proposition \ref{prop:supportsmooth} and Lebesgue's dominated convergence theorem, we get Theorem \ref{theosmooth} (ii). $\square$

\pass We only need now to explain how to estimate the mean defect support function in a fixed direction. 
To do so, let us introduce the support point $m_{s,\lambda}$ on $\partial K_{\lambda}$ in direction $n_s$, i.e. the point which satisfies $\langle m_{s,\lambda}, n_s\rangle=p_o(K_{\lambda},n_s)$.
Denoting by $X_{s,\lambda}=\langle m_{s,\lambda},t_s\rangle$ and $Y_{s,\lambda}= p_o(K_{\lambda},n_s)-p_o(K,n_s)$, we can write
\begin{equation*}
m_{s,\lambda}=s+X_{s,\lambda}t_{s}+Y_{s,\lambda} n_s
\end{equation*}
where $(t_s,n_s)$ stands for the Frenet frame at point $s$.

\pass The next proposition investigates the asymptotic distribution of the couple $(X_{s,\lambda},Y_{s,\lambda})$ and provides the required asymptotic estimate for $\E(p_o(K_{\lambda},n_s)-p_o(K,n_s))$.

\begin{prop}\label{prop:supportsmooth}
\noi
\begin{enumerate}
\item[$(i)$] For every $s\in\partial K$, the couple $(\lambda^{\frac1{3}}X_{s,\lambda}, \lambda^{\frac{2}{3}}Y_{s,\lambda})$ converges in distribution when $\lambda\to\infty$
to the distribution with density function $f_s$ given by
\begin{equation*}\label{eq:ftheta}
\hspace{1.5cm} f_s(x,y)=2^{\frac{11}{2}}\langle s,n_s\rangle^2 r_s^{-\frac{3}{2}}\exp\Big(-2^{\frac{9}{2}} 3^{-1}r_s^{-\frac1{2}}\langle s,n_s\rangle
\Big(\frac{x^2}{2r_s}+y\Big)^{\frac{3}{2}}\Big)\Big(\frac{x^2}{2r_s}+y\Big)^{\frac1{2}}y{\emph{\indicat}}_{\{y>0\}}.
\end{equation*}
\item[$(ii)$] There exists $C>0$ such that for every $s\in\partial K$ and $\lambda>0$, $\lambda^{\frac{2}{3}}\E(Y_{s,\lambda})\le C$. Moreover, for every $s\in\partial K$,
\begin{equation*}\label{eq:convmean}
\E(Y_{s,\lambda})=\E(p_o(K_{\lambda},n_s)-p_o(K,n_s))\underset{\lambda\to\infty}{\sim}\lambda^{-\frac{2}{3}}3^{-\frac{4}{3}}\Gamma\Big(\mfrac{2}{3}\Big) r_s^{\frac1{3}}\langle s,n_s\rangle^{-\frac{2}{3}}.
\end{equation*}
\end{enumerate}
\end{prop}

\pass\textbf{Proof of Proposition \ref{prop:supportsmooth}.}

\noi {\it Proof of $(i)$.} We first notice that the point $m_{s,\lambda}$ is necessarily one of the vertices of ${K_{\lambda}}$, i.e. is at the intersection of two bisecting lines between $o$ and two Voronoi neighbors of $o$. For $x_1,x_2\in\cP_{\lambda}\setminus2\cF_o(K)$, we denote by $c_{x_1,x_2}$ the intersection point of the two bisecting lines of the segments $[o,x_1]$ and $[o,x_2]$. In particular,
\begin{equation*}\label{eq:equivalenceplusclaire}
(c_{x_1,x_2}=m_{s,\lambda})\Longleftrightarrow
\left\{ \begin{array}{l}
         c_{x_1,x_2} \mbox{ is extreme in direction $n_s$} \\
         B_{\|c_{x_1,x_2}\|}(c_{x_1,x_2})\cap (\cP_{\lambda}\setminus2\cF_o(K))=\emptyset
         \end{array}\right..
\end{equation*}

\begin{figure}[!h]
\begin{center}
\includegraphics[scale=0.75]{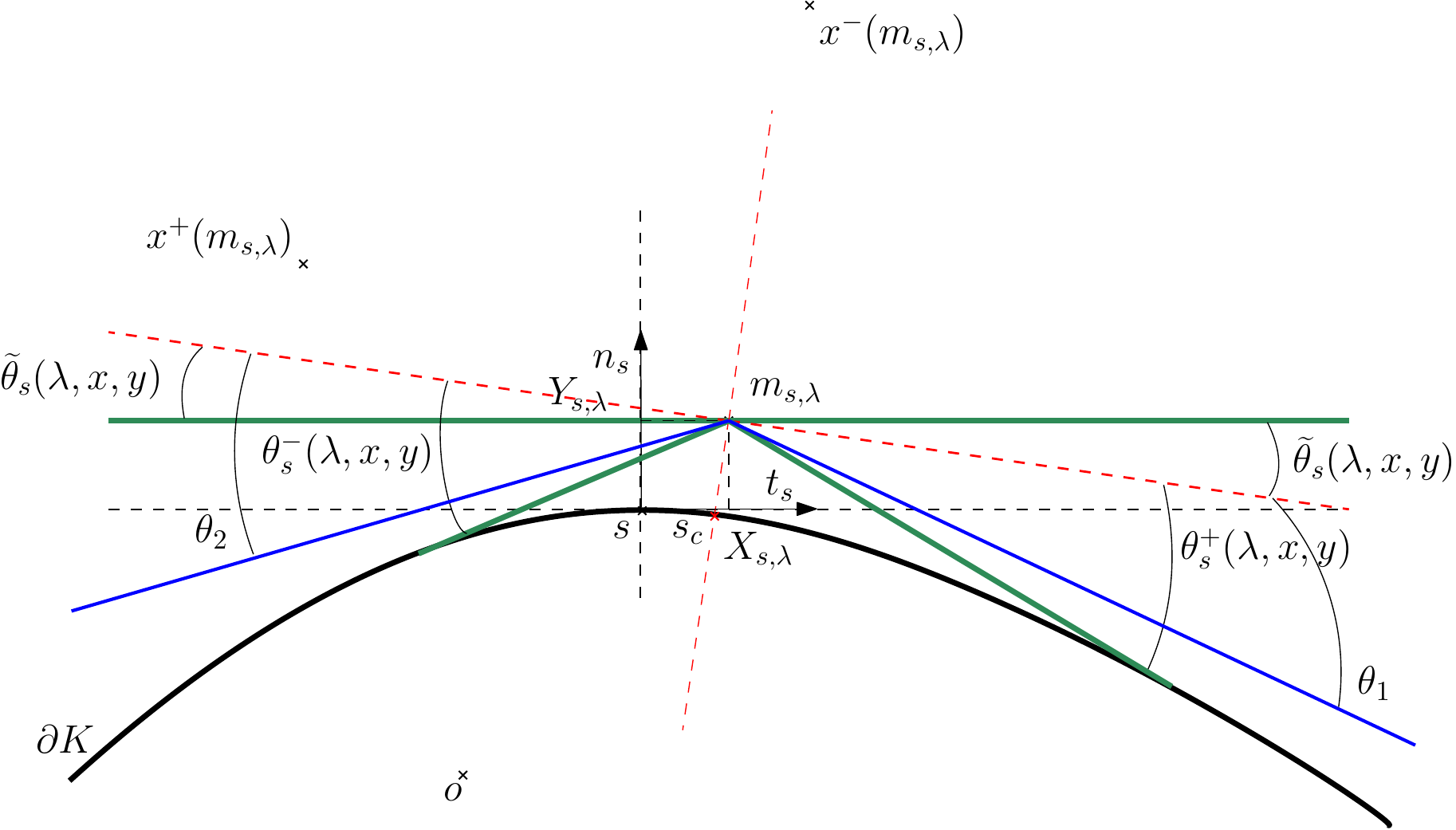}
\vspace{-0.2cm}
\caption{\label{fig:supportsmooth} Realization of the support function.}
\end{center}
\end{figure}

\noi From a given $m_{s,\lambda}$ emanate two segments, one on the left of the half-line $\R_+m_{s,\lambda}$, one on the right. The symmetric points of $o$ with respect to these two segments define the right and the left Poisson-Voronoi neighbors of $o$ with respect to $m_{s,\lambda}$. They will be denoted by $x^+(m_{s,\lambda})$ and $x^-(m_{s,\lambda})$ respectively. Consequently, by Mecke-Slivnyak's formula, for every positive and measurable function $\varphi:\R^2\longrightarrow \R_+$,
\begin{equation*}\label{eq:convloi}
\begin{split}
\E\big(\varphi&(\lambda^{\frac1{3}}X_{s,\lambda},\lambda^{\frac{2}{3}}Y_{s,\lambda})\big) \\
& = \E\bigg(\sum_{(x_1,x_2)\in(\cP_{\lambda}\setminus2\cF_o(K))^2} {\emph{\indicat}}_{\left\{\substack{\tiny c_{x_1,x_2}=m_{s,\lambda}\\x_1=x^+(m_{s,\lambda})\\x_2=x^-(m{s,\lambda})}\right\}}
\varphi(\lambda^{\frac1{3}}\langle c_{x_1,x_2}-s,t_s\rangle,\lambda^{\frac{2}{3}}\langle c_{x_1,x_2}-s,n_s\rangle)\bigg) \\
& = \lambda^2\int_{\R^2\times\R^2}\exp\left(-4\lambda\cA(B_{\|c_{x_1,x_2}\|}(c_{x_1,x_2})\setminus\cF_o(K))\right)
{\emph{\indicat}}_{\left\{
\substack{\mbox{\tiny $c_{x_1,x_2}$ is extreme in direction $n_s$}\\\mbox{{\tiny $x_1=x^+(c_{x_1,x_2})$, $x_2=x^-(c_{x_1,x_2})$}}}
\right\}}
\\
& \hspace*{5cm}\times\varphi(\lambda^{\frac{1}{3}}\langle c_{x_1,x_2}-s,t_s\rangle,\lambda^{\frac{2}{3}}\langle c_{x_1,x_2}-s,n_s\rangle)\dd x_1 \dd x_2.
\end{split}
\end{equation*}

\pass Let $s_c$ be the orthogonal projection of $c_{x_1,x_2}$ onto $K$ and $t_{s_c}=u_{\gamma}$ be the unit outer normal vector of $\partial K$ at $s_c$.
We now apply two consecutive changes of variables in the integral above.

\pass First, we write $c_{x_1,x_2}= ru_{\theta+\gamma}$ and denote by $\theta_1$ and $\theta_2$ the angles between one of the two bisecting lines emanating from $c_{x_1,x_2}$ and $t_{s_c}$. We then use the following lemma providing a change of variables formula which may be understood as a classical formula {\it \`a la} Blaschke-Petkantschin, see e.g. Theorem 7.3.1. from \cite{SW}. It consists essentially in the computation of the Jacobian of a four dimensional transformation.

\begin{lem}\label{lem:jacob}
Let $x=r u_{\theta}$, $r>0$, $\theta\in (0,2\pi)$ and $\theta-\pi<\theta_1<\theta_2<\theta$. Let $x_1=2r\sin(\theta-\theta_1)u_{\theta_1+\frac{\pi}{2}}$ and $x_2=2r\sin(\theta-\theta_2)u_{\theta_2+\frac{\pi}{2}}$ be the symmetric points of the origin $o$ with respect to the lines $x+\R u_{\theta_1}$ and $x+\R u_{\theta_2}$ respectively. Then the Jacobian of the change of variables $(r,\theta,\theta_1,\theta_2)\longmapsto (x_1,x_2)$ is given by
\begin{equation*}
\frac{\dd x_1 \dd x_2}{r \dd r \dd \theta\dd \theta_1\dd \theta_2} = 16 r^2 J(\theta,\theta_1,\theta_2)
\end{equation*}
with
\begin{equation*}
J(\theta,\theta_1,\theta_2) = |\sin(\theta_1-\theta_2)\sin(\theta-\theta_2)\sin(\theta-\theta_1)|.
\end{equation*}
\end{lem}

\noi{\bf Proof of Lemma \ref{lem:jacob}.}
We omit the calculation which is analogous to the proof of the classical Blaschke-Petkantschin's formula, see e.g. Theorem 7.3.1 in \cite{SW}. $\square$

\pass Secondly, we replace the couple $(r,\theta)$ by $(x,y)$ defined by
\begin{equation*}\label{eq:xyh}
x=\lambda^{\frac1{3}}\langle c_{x_1,x_2}-s,t_s\rangle\,\text{ and }\, y=\lambda^{\frac{2}{3}}\langle c_{x_1,x_2}-s,n_s\rangle.
\end{equation*}
We get in particular
\begin{equation}\label{eq:defrho}
r^2=\rho_s^2(\lambda,x,y)=(\langle s,t_s\rangle+\lambda^{-\frac{1}{3}}x)^2+(\langle s,n_s\rangle+\lambda^{-\frac{2}{3}}y)^2
\end{equation}
and a Jacobian given by
\begin{equation*}\label{eq:jacsimple}
\frac{r\dd r\dd \theta}{\dd x \dd y} = \lambda^{-1}.
\end{equation*}

\pass Consequently, we deduce that
\begin{align}\label{eq:formuleexacte}
\E\big(\varphi(\lambda^{\frac1{3}}X_{s,\lambda},\lambda^{\frac{2}{3}}Y_{s,\lambda})\big) = 16\int_{\R\times \R_+}\exp(-\Delta_s(\lambda,x,y))\varphi(x,y)\rho_s^2(\lambda,x,y)J^{\mbox{\scalebox{.5}{supp}}}_s(\lambda,x,y)\dd x \dd y
\end{align}
where
\begin{equation*}
\Delta_s(\lambda,x,y)=4\lambda\cA(B_{\|c_{x_1,x_2}\|}(c_{x_1,x_2})\setminus \cF_o(K))
\end{equation*}
and
\begin{equation*}
J^{\mbox{\scalebox{.5}{supp}}}_s(\lambda,x,y)=\lambda\int_{E^{\mbox{\scalebox{.5}{supp}}}_{s,x,y}} J(\theta_s(\lambda,x,y),\theta_1,\theta_2)\dd \theta_1\dd \theta_2
\end{equation*}
with
\begin{equation*}
\cos(\theta_s(\lambda,x,y))=\bigg\langle \frac{s+\lambda^{-\frac1{3}}xt_s+\lambda^{-\frac{2}{3}}yn_s}{\|s+\lambda^{-\frac{1}{3}}xt_s+\lambda^{-\frac{2}{3}}yn_s\|},t_{s_c}\bigg\rangle
\end{equation*}
and $E^{\mbox{\scalebox{.5}{supp}}}_{s,x,y}$ the set of couples $(\theta_1,\theta_2)$ which satisfy that $c_{x_1,x_2}$ is extremal in the direction of $n_s$ and that the two bisecting lines of $[o,x_1]$ and $[o,x_2]$ do not intersect $K$.

\pass Let us make the set $E^{\mbox{\scalebox{.5}{supp}}}_{s,x,y}$ explicit. Let $\theta_s^+(\lambda,x,y)$, $\theta_s^-(\lambda,x,y)$ and $\widetilde{\theta}_s(\lambda,x,y)$ be, respectively, the angle of aperture at the point $c_{x_1,x_2}=s+\lambda^{-\frac1{3}}x t_s+\lambda^{-\frac{2}{3}}y n_s$ on the right, the angle of aperture at $c_{x_1,x_2}$ on the left and the angle between the vectors $n_s$ and $n_{s_c}$ (see Figure \ref{fig:supportsmooth}).
We obtain
\begin{equation}\label{eq:condanglessmooth}
(\theta_1,\theta_2)\in E^{\mbox{\scalebox{.5}{supp}}}_{s,x,y}\Longleftrightarrow -\theta_s^+(\lambda,x,y)< \theta_1<\widetilde{\theta}_s(\lambda,x,y))<\theta_2<\theta_s^-(\lambda,x,y).
\end{equation}

\pass In order to show the required convergence in distribution, we are going to use Lebesgue's dominated convergence theorem. To do so, we need to prove the convergence of the integrand in \eqref{eq:formuleexacte} and that it is dominated.

\pass $-$ \emph{Convergence and domination of} $\rho_s(\lambda,x,y)$.

\noi We deduce from \eqref{eq:defrho} that
\begin{equation}\label{eq:convrho}
\rho_s(\lambda,x,y)\underset{\lambda\to\infty}{\longrightarrow}\|s\|.
\end{equation}

\pass Moreover, by triangular inequality, we get for all $\lambda\ge 1$,
\begin{equation}\label{eq:dominationrho}
\rho_s(\lambda,x,y)\le \|s\|+(\lambda^{-\frac{2}{3}}x^2+\lambda^{-\frac{4}{3}}y^2)^{\frac1{2}}\le \|s\|+\|(x,y)\|.
\end{equation}

\pass $-$ \emph{Convergence and domination of} $\exp(-\Delta_s(\lambda,x,y))$.

\noi We denote by $h$ the distance from $c_{x_1,x_2}=s+\lambda^{-\frac1{3}}x t_s+\lambda^{-\frac{2}{3}}y n_s$ to $K$. Then the following relation holds,
uniformly in $s$,
\begin{equation}\label{eq:hauteurdevue}
h\underset{\lambda\to\infty}{=}\lambda^{-\frac{2}{3}}\Big(\frac{x^2}{2r_s}+y\Big)+\mathrm{O}(\lambda^{-\frac{4}{3}}).
\end{equation}

\noi Moreover, there exists $C>0$ such that for every $s\in \partial K$ and all $\lambda\ge 1 $,
\begin{equation}\label{eq:dominationh}
C^{-1}\min\Big(\lambda^{-\frac1{3}}\Big(x^2+2yr_s\Big)^{\frac1{2}},\lambda^{-\frac{2}{3}}\Big(\frac{x^2}{r_s}+2y\Big)\Big)\le h\le C \lambda^{-\frac{2}{3}}\Big(\frac{x^2}{2r_s}+y+\frac{y^2}{2r_s}\Big).
\end{equation}

\noi Indeed, let us prove first \eqref{eq:hauteurdevue} and \eqref{eq:dominationh} when $K$ is a disk of radius $r_s$.
We find on the one hand by Pythagora's theorem
\begin{equation*}
\|c_{x_1,x_2}\|^2-r_s^2=(r_s+\lambda^{-\frac{2}{3}}y)^2+\lambda^{-\frac{2}{3}}x^2-r_s^2=\lambda^{-\frac{2}{3}}(x^2+\lambda^{-\frac{2}{3}}y^2+2r_sy)
\end{equation*}
and on the other hand
\begin{equation*}
\|c_{x_1,x_2}\|^2-r_s^2=(\|c_{x_1,x_2}\|-r_s)(\|c_{x_1,x_2}\|+r_s)=h(2r_s+h).
\end{equation*}
Combining these two equalities, we get
\begin{equation*}
h=r_s\Big(1+\lambda^{-\frac{2}{3}}\Big(\frac{x^2}{r_s^2}+\frac{2y}{r_s}+\lambda^{-\frac{2}{3}}\frac{y^2}{r_s^2}\Big)\Big)^{\frac1{2}}-r_s.
\end{equation*}

\noi Using the estimate $\mfrac1{4}\min(u,u^{\frac1{2}})\le (1+u)^{\frac1{2}}-1\le\mfrac{u}{2}$ for every $u>0$ we obtain that the previous equality implies both
\eqref{eq:hauteurdevue} and \eqref{eq:dominationh} when $K$ is a disk.

\pass When $K$ is a smooth convex body, we sandwich its boundary between two disks of radii $r_s+C'\lambda^{-\frac{2}{3}}$ and $r_s-C'\lambda^{-\frac{2}{3}}$ for a fixed positive constant $C'>0$ and we deduce from the previous case both \eqref{eq:hauteurdevue} and \eqref{eq:dominationh} for $K$.

\pass Moreover, using the regularity assumptions on the boundary $\partial K$, we get, uniformly in $s$,
\begin{equation}\label{eq:petitdetailtechnique}
|\langle s_c,n_{s_c}\rangle- \langle s,n_s\rangle|\underset{\lambda\to\infty}{=} \mathrm{O}(\lambda^{-\frac1{3}}).
\end{equation}

\pass Thanks to Lemma \ref{lem:geom}, \eqref{eq:hauteurdevue} and \eqref{eq:petitdetailtechnique}, we get, uniformly in $s$,
\begin{align}\label{eq:limiteaciter}
\Delta_s(\lambda,x,y) \underset{\lambda\to\infty}{\sim} 2^{\frac{9}{2}}3^{-1}r_{s}^{-\frac1{2}}\langle s,n_s\rangle\Big(\frac{x^2}{2r_{s}}+y\Big)^{\frac{3}{2}}.
\end{align}

\noi In particular, thanks to Lemma \ref{lem:geom} (ii) and \eqref{eq:dominationh}, there exists $C>0$, uniform in $s$, such that $\Delta_s(x,y)$ satisfies,
for all $\lambda\ge 1$:
\begin{equation}\label{eq:dominationDelta}
\Delta_s(\lambda,x,y)\ge C \min\Big(\Big(x^2+2r_sy\Big)^{\frac{3}{4}},\Big(\frac{x^2}{2r_s}+y\Big)^{\frac{3}{2}}\Big).
\end{equation}

\pass $-$ \emph{Convergence and domination of} $J^{\mbox{\scalebox{.5}{supp}}}_s(\lambda,x,y)$.

\noi We start by estimating the function $\theta_s(\lambda,x,y)$. We get
\begin{equation}\label{eq:estimthetaxy}
\theta_s(\lambda,x,y)\underset{\lambda\to\infty}{=}\arcsin\bigg(\bigg\langle \frac{s}{\|s\|},n_s\bigg\rangle\bigg)+\mathrm{O}(\lambda^{-\frac1{3}}).
\end{equation}

\noi We now estimate the three angles $\theta_s^+(\lambda,x,y)$, $\theta_s^-(\lambda,x,y)$ and $\widetilde{\theta}_s(\lambda,x,y)$. Using \eqref{eq:angleaperture} and \eqref{eq:hauteurdevue}, we get
\begin{equation}\label{eq:estimtheta+}
\theta_s^+(\lambda,x,y)  \underset{\lambda\to\infty}{=}\theta_s^-(\lambda,x,y)+\mathrm{O}(\lambda^{-1})
  \underset{\lambda\to\infty}{=}\lambda^{-\frac1{3}} 2^{\frac1{2}}r_s^{-\frac1{2}}\Big(\frac{x^2}{2r_s}+y\Big)^{\frac1{2}}+ \mathrm{O}(\lambda^{-1}).
\end{equation}

\noi Thanks to \eqref{eq:dominationangleouverture} and \eqref{eq:dominationh}, we have additionally the inequality, for some $C>0$,
\begin{equation}\label{eq:dominationtheta+}
\theta_s^+(\lambda,x,y)\le 2^{\frac1{2}}r_s^{-\frac1{2}}h^{\frac1{2}}\le C r_s\lambda^{-\frac1{3}}\Big(\frac{x^2}{2r_s}+y+\frac{y^2}{2r_s}\Big)^{\frac1{2}}.
\end{equation}

\noi The same inequality holds for $\theta_s^+(\lambda,x,y)$. We turn now our attention to $\widetilde{\theta}_s(\lambda,x,y)$. When $K$ is a disk we get
\begin{equation*}
\widetilde{\theta}_s(\lambda,x,y)=\arctan\bigg(\frac{\lambda^{-\frac1{3}}x}{r_s+\lambda^{-\frac{2}{3}}y}\bigg)
\underset{\lambda\to\infty}{=}\lambda^{-\frac1{3}}r_s^{-1}x+ \mathrm{O}(\lambda^{-1}).
\end{equation*}

\noi When $K$ is a smooth convex body, we sandwich again its boundary between two disks of radii $r_s+C\lambda^{-\frac{2}{3}}$ and $r_s-C\lambda^{-\frac{2}{3}}$ and we obtain
\begin{equation*}\label{eq:estimtildetheta}
\widetilde{\theta}_s(\lambda,x,y)\underset{\lambda\to\infty}{=}\lambda^{-\frac1{3}}r_s^{-1}x+\mathrm{O}(\lambda^{-1}).
\end{equation*}

\noi Consequently, we deduce from \eqref{eq:estimthetaxy}, \eqref{eq:condanglessmooth} and \eqref{eq:estimtheta+} that
\begin{align}\label{eq:intangles}
J^{\mbox{\scalebox{.5}{supp}}}_{s}(\lambda,x,y)
& \underset{\lambda\to\infty}{\sim} \lambda\int_{E^{\mbox{\scalebox{.5}{supp}}}_{s,x,y}}|\sin(\theta_2-\theta_1)\sin(\theta_s(\lambda,x,y)-\theta_2)\sin(\theta_s(\lambda,x,y)-\theta_1)|
     \dd \theta_1\dd \theta_2 \nonumber \\
& \underset{\lambda\to\infty}{\sim}\lambda\sin^2\Big(\arcsin\Big(\Big\langle\frac{s}{\|s\|},n_s\Big\rangle\Big)\Big)
   \int_{\theta_1=-\theta_s^+(\lambda,x,y)}^{\widetilde{\theta}_s(\lambda,x,y))}
   \int_{\theta_2=\widetilde{\theta}_s(\lambda,x,y))}^{\theta_s^-(\lambda,x,y)} (\theta_2-\theta_1)\dd \theta_1 \dd\theta_2 \nonumber\\
& \underset{\lambda\to\infty}{\sim} \lambda \Big\langle \frac{s}{\|s\|},n_s\Big\rangle^2 \theta_s^+(\lambda,x,y)(\theta_s^+(\lambda,x,y)^2-\widetilde{\theta}_s(\lambda,x,y)^2) \nonumber\\
& \underset{\lambda\to\infty}{\sim} \Big\langle \frac{s}{\|s\|},n_s\Big\rangle^2 2^{\frac{3}{2}}r_s^{-\frac{3}{2}}y\left(\frac{x^2}{2r_s}+y\right)^{\frac1{2}}.
\end{align}

\pass Moreover, thanks to \eqref{eq:dominationtheta+}, we get for some $C>0$,
\begin{equation}\label{eq:dominationJbarre}
J^{\mbox{\scalebox{.5}{supp}}}_{s}(\lambda,x,y)\le  C r_s^3\Big(\frac{x^2}{2r_s}+y+\frac{y^2}{2r_s}\Big)^{\frac{3}{2}}.
\end{equation}

\pass $-$ \emph{Conclusion}.

\noi Combining \eqref{eq:convrho}, \eqref{eq:limiteaciter} and \eqref{eq:intangles}, we obtain that the integrand in \eqref{eq:formuleexacte} converges to
\begin{equation*}
8\|s\|^2\varphi(x,y)\exp\Big(-2^{-\frac{9}{2}}3^{-1}r_{s}^{-\frac1{2}}\big\langle s,n_s\rangle\Big(\frac{x^2}{2r_{s}}+y\Big)^{\frac{3}{2}} \Big)\Big\langle\frac{s}{\|s\|},n_s\Big\rangle^2 2^{\frac{3}{2}}r_s^{-\frac{3}{2}}y\left(\frac{x^2}{2r_s}+y\right)^{\frac1{2}}.
\end{equation*}

\pass Now the estimates \eqref{eq:dominationrho}, \eqref{eq:dominationDelta} and \eqref{eq:dominationJbarre} show that we can apply Lebesgue's dominated convergence theorem for any function $\varphi$ bounded by a polynomial of $x$ and $y$, say. This proves assertion (i).

\pass{\it Proof of $(ii)$.}
We start by rewriting the proof of (i) when $\varphi(x,y)=y$. Since $K$ is a compact convex set with bounded positive curvature and containing the origin in its interior, the non-negative quantities $\|s\|$ and $r_s$ are bounded from above and from below. Consequently, the estimates \eqref{eq:dominationrho}, \eqref{eq:dominationDelta} and \eqref{eq:dominationJbarre} imply that the integral on the right-hand side of \eqref{eq:formuleexacte} is bounded independently of $s$, i.e. that there exists $C>0$ such that $\lambda^{\frac{2}{3}}\E(Y_{s,\lambda})\le C $ for every $s\in\partial K$.
When applying Lebesgue's dominated convergence theorem, we get that
\begin{equation*}\label{eq:premiereformule}
\lambda^{\frac{2}{3}}\E(Y_{s,\lambda}) \underset{\lambda\to\infty}{\longrightarrow}  2^{\frac{11}{2}} r_s^{-\frac{3}{2}}\langle s,n_s\rangle^2I_s
\end{equation*}
where
\begin{equation*}\label{eq:defIs}
I_s =\int_{\R\times\R_+}\exp\Big(-2^{\frac{9}{2}} 3^{-1}r_s^{-\frac1{2}}\langle s,n_s\rangle \Big(\frac{x^2}{2r_s}+y\Big)^{\frac{3}{2}}\Big)\Big(\frac{x^2}{2r_s}+y\Big)^{\frac1{2}}y^2 \dd x \dd y.
\end{equation*}

\noi It remains to make the integral $I_s$ explicit. Recalling that $C_s= 2^{\frac{9}{2}}3^{-1}r_s^{-\frac1{2}}\,\big\langle s,n_s\big\rangle$, we get
\begin{align*}\label{eq:decompIs}
I_s = \frac{4}{3 C_s}\int_0^{\infty}\bigg(\int_0^{\infty} e^{-u}\Big(\frac{u^{\frac{2}{3}}}{C_s^{\frac{2}{3}}}-\frac{x^2}{2r_s}\Big)^2
{\emph{\indicat}}_{\{0<x<2^{\frac1{2}}r_s^{\frac1{2}}C_s^{-\frac1{3}}u^{\frac1{3}}\}}\dd x\bigg)\dd u =I_s^{(1)}+I_s^{(2)}-I_s^{(3)}
\end{align*}
where
\begin{align*}
I_s^{(1)} & = \frac{4}{3C_s}\int_0^{\infty}e^{-u}C_s^{-\frac{4}{3}}u^{\frac{4}{3}}2^{\frac1{2}}r_s^{\frac1{2}}C_s^{-\frac1{3}}u^{\frac1{3}}\dd u
= 2^{\frac{5}{2}}3^{-2}5C_s^{-\frac{8}{3}}r_s^{\frac{1}{2}}\Gamma\Big(\frac{5}{3}\Big) , \\
I_s^{(2)} & = \frac{4}{3C_s}\int_0^{\infty}e^{-u}20^{-1}r_s^{-2}(2^{\frac1{2}}r_s^{\frac1{2}}C_s^{-\frac1{3}}u^{\frac1{3}})^5\dd u
= 2^{\frac{5}{2}}3^{-2}C_s^{-\frac{8}{3}}r_s^{\frac{1}{2}}\Gamma\Big(\frac{5}{3}\Big) , \\
I_s^{(3)} & = \frac{4}{3C_s}\int_0^{\infty}e^{-u}C_s^{-\frac{2}{3}}r_s^{-1}3^{-1}
u^{\frac{2}{3}}(2^{\frac1{2}}r_s^{\frac1{2}} C_s^{-\frac1{3}}u^{\frac1{3}})^3\dd u
= 2^{\frac{7}{2}}3^{-3}5C_s^{-\frac{8}{3}}r_s^{\frac{1}{2}}\Gamma\Big(\frac{5}{3}\Big).
\end{align*}

\noi Finally, inserting these equalities into $I_s$ yields the required result. $\square$

\subsection{Proof of Theorem \ref{theosmooth} (iii): intensity and number of vertices}\label{subsec:verticessmooth}

\pass

\noi Using \eqref{eq:chgtsr}, we rewrite \eqref{eq:efronvertices} when $K$ is smooth and we obtain the following relation, when the intensity $\lambda\to\infty$,
\begin{equation*}\label{eq:celinedion}
\E(\cN(K_{\lambda}))\underset{\lambda\to\infty}{\sim} 4\lambda\int_{\partial K}\langle s,n_s\rangle\E(Y_{s,\lambda})r_s^{-1}\dd s.
\end{equation*}
Theorem \ref{theosmooth} (iii) is then deduced from Proposition \ref{prop:supportsmooth} (ii). $\square$

\pass Actually, we can provide a more precise result on the asymptotic intensity of the point process of vertices of $K_{\lambda}$. This new result that we describe below could alternatively be used to get Theorem \ref{theosmooth} (iii) via an integration of the intensity given in the next proposition.

\pass Let us fix a vertex $s\in\partial K$ and consider the point process $V_\lambda$ of vertices of $K_{\lambda}$. We rewrite a point $v\in V_\lambda$ as $v=s+x_v t_s+y_v n_s$.

\begin{prop}\label{prop:intensitysmooth}
Consider the point process $(x_v,y_v)_{v\in V_\lambda}$ of the vertices of $K_{\lambda}$.
Then the rescaled point process $(\lambda^{\frac{1}{3}}x_v,\lambda^{\frac{2}{3}}y_v)_{v\in V_\lambda}$ has an asymptotic intensity given by
\begin{equation*}\label{eq:intensitysmooth}
\sigma_s(x,y)=2^{\frac{15}{2}}3^{-1}\exp\bigg(-2^{\frac{9}{2}}3^{-1} \bigg(\frac{x^2}{2r_s}+y\bigg)^{\frac{3}{2}} r_{s}^{-\frac1{2}}\langle s,n_s\rangle\bigg)
\langle s,n_s\rangle^2r_s^{-\frac{3}{2}}\bigg(\frac{x^2}{2r_s}+y\bigg)^{\frac{3}{2}}.
\end{equation*}
\end{prop}

\noindent
Proposition \ref{prop:intensitysmooth}
provides an extra valuable information on the point process of vertices which is clearly not of Poisson type. To some extent, this is also reminiscent of the description of the rescaled point process of vertices of random polytopes in the unit-ball or random Gaussian polytopes, as a sub-product of a growth parabolic process based on a Poisson point process, see e.g. Theorem 1.1 in \cite{CY3}.

\pass\textbf{Proof of Proposition \ref{prop:intensitysmooth}.}
Let $X\times Y \subset \R^2\setminus K$ be fixed and denote by $\cN_s(X\times Y)$ the number of points of the process $(\lambda^{\frac1{3}}x_v,\lambda^{\frac{2}{3}}y_v)_{v\in V_\lambda}$ belonging to the set $X\times Y$. We have to show that
\begin{align*}
\E(\cN_s(X\times Y)) \underset{\lambda\to\infty}{\longrightarrow} \int_{X\times Y} \sigma_s(x,y)\dd x \dd y.
\end{align*}

\noi The strategy is very similar to the proof of Proposition \ref{prop:supportsmooth}, i.e. it consists in applying Mecke-Slivnyak's formula, then the change of variables provided by Lemma \ref{lem:jacob} and finally Lebesgue's dominated convergence theorem.

\pass Consequently, we deduce that
\begin{align}\label{eq:calculintermcotes1}
\E(\cN_s(X\times Y)) = 16\int_{\R^2}\exp(-\Delta_s(\lambda,x,y))\rho_s^2(\lambda,x,y)J^{\mbox{\scalebox{.5}{vert}}}_s(\lambda,x,y)\dd x \dd y
\end{align}
where
\begin{equation*}
J^{\mbox{\scalebox{.5}{vert}}}_s(\lambda,x,y)
=\lambda\int_{E^{\mbox{\scalebox{.5}{vert}}}_{s,x,y}} J(\theta_s(\lambda,x,y),\theta_1,\theta_2)\dd \theta_1\dd \theta_2
\end{equation*}
and $E^{\mbox{\scalebox{.5}{vert}}}_{s,x,y}$ is the set of couples $(\theta_1,\theta_2)$ which satisfy that the two bisecting lines of $[o,x_1]$ and $[o,x_2]$ do not intersect $K$.

\pass Let us make the set $E^{\mbox{\scalebox{.5}{vert}}}_{s,x,y}$ noticing that
\begin{equation}\label{eq:condangles}
(\theta_1,\theta_2)\in E^{\mbox{\scalebox{.5}{vert}}}_{s,x,y}\Longleftrightarrow -\theta_s^+(\lambda,x,y)< \theta_1<\theta_2<\theta_s^-(\lambda,x,y).
\end{equation}

\begin{figure}[!h]
\begin{center}
\includegraphics[scale=0.75]{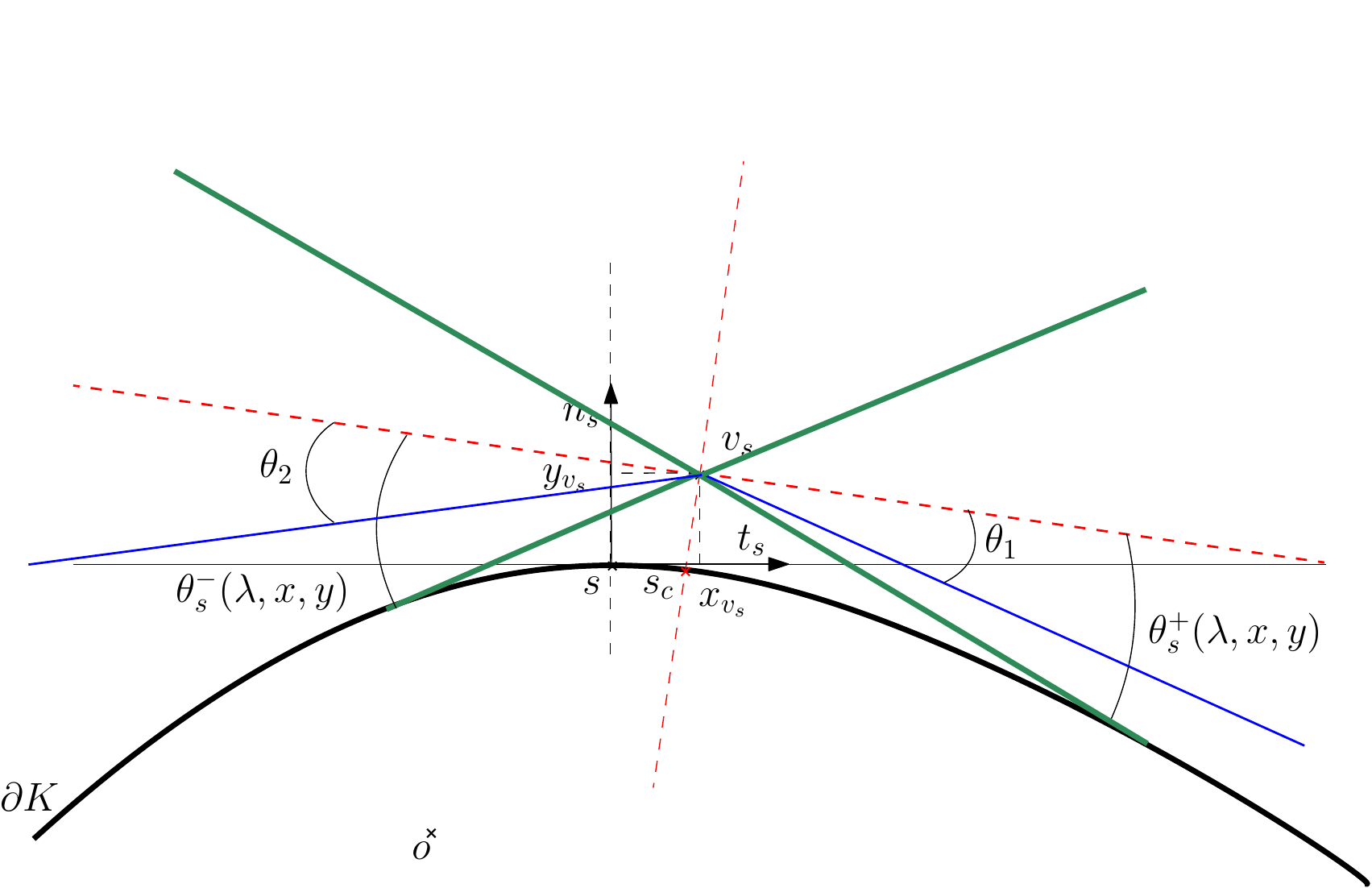}
\vspace{-0.2cm}
\caption{\label{fig:sommetsmooth} Intensity of vertices near a fixed point $s\in\partial K$.}
\end{center}
\end{figure}

\pass The convergence and domination of $\rho_s(\lambda,x,y)$, $\exp(-\Delta_s(\lambda,x,y))$, $\theta_s^+(\lambda,x,y)$ and $\theta_s^-(\lambda,x,y)$ is identical to what has been done in the proof of Proposition \ref{prop:supportsmooth}. We turn our attention to the convergence of $J^{\mbox{\scalebox{.5}{vert}}}_s(\lambda,x,y)$.

\begin{align}\label{eq:convJbarre2}
J^{\mbox{\scalebox{.5}{vert}}}_{s}(\lambda,x,y)
& \underset{\lambda\to\infty}{\sim}\lambda \int_{E^{\mbox{\scalebox{.5}{vert}}}_{s,x,y}}|\sin(\theta_1-\theta_2)\sin(\theta_s(\lambda,x,y)-\theta_2)\sin(\theta_s(\lambda,x,y)-\theta_1)|
\dd \theta_1\dd \theta_2 \nonumber \\
& \underset{\lambda\to\infty}{\sim} \lambda\sin^2\Big(\arcsin\Big(\Big\langle \frac{s}{\|s\|},n_s\Big\rangle\Big)\Big) \int_{\theta_1=-\theta_s^+(\lambda,x,y)}^{{\theta}_s^-(\lambda,x,y))}\int_{\theta_2=\theta_1}^{\theta_s^-(\lambda,x,y)}(\theta_2-\theta_1) \dd \theta_2 \dd\theta_1 \nonumber\\
& \underset{\lambda\to\infty}{\sim} \mfrac{4}{3}\lambda \Big\langle \frac{s}{\|s\|},n_s\Big\rangle^2 \theta_s^+(\lambda,x,y)^3 \nonumber \\
& \underset{\lambda\to\infty}{\sim} 2^{\frac{7}{2}}3^{-1}\Big\langle \frac{s}{\|s\|},n_s\Big\rangle^2 r_s^{-\frac{3}{2}}\bigg(\frac{x^2}{2r_s}+y\bigg)^{\frac{3}{2}}.
\end{align}

\pass Moreover, thanks to \eqref{eq:dominationtheta+}, we get for some positive constant $C>0$,
\begin{equation}\label{eq:dominationJbarre2}
J^{\mbox{\scalebox{.5}{vert}}}_{s}(\lambda,x,y)\le  C \Big\langle \frac{s}{\|s\|},n_s\Big\rangle^2r_s^{-\frac{3}{2}}\bigg(\frac{x^2}{2r_s}+y\bigg)^{\frac{3}{2}}.
\end{equation}

\pass Combining \eqref{eq:convrho}, \eqref{eq:limiteaciter} and \eqref{eq:convJbarre2}, we obtain that the integrand in \eqref{eq:calculintermcotes1} converges to
\begin{equation*}
2^{\frac{15}{2}}3^{-1}\exp\Big(-2^{\frac{9}{2}}3^{-1} \bigg(\frac{x^2}{2r_s}+y\bigg)^{\frac{3}{2}} r_{s}^{-\frac1{2}}\langle s,n_s\rangle\Big)
\langle s,n_s\rangle^2r_s^{-\frac{3}{2}}\bigg(\frac{x^2}{2r_s}+y\bigg)^{\frac{3}{2}}.
\end{equation*}

\pass Now the estimates \eqref{eq:dominationrho}, \eqref{eq:dominationDelta} and \eqref{eq:dominationJbarre2} show that we can apply Lebesgue's dominated convergence theorem. The result follows. $\square$

\subsection{A R\'enyi-Sulanke approach}\label{sec:rsapproach}

\pass

\noi In this section, we explain how points (i) and (iii) of Theorem \ref{theosmooth} could be deduced from Theorem A applied to a disk and from a similar mean estimate from \cite{s02}. We describe below the method which is based on the application of an inversion with respect to the osculating disk at $s$ for any $s\in\partial K$. This is reminiscent of both \cite{cs05} for the idea of transforming a Voronoi cell into a convex hull by an inversion and \cite{cy14} for the rewriting of the expectations as an integral over $\partial K$ of a mean of a so-called {\it score} and the replacement of $K$ by a disk in the calculation of the score. We only sketch below the main steps of the approach, as a thorough proof would involve many more technical details. It is certainly less natural than the method used previously and also specific to the smooth case, as well as to the two functionals $\cA(K_\lambda)$ and $\cN(K_\lambda)$. Nevertheless, we have chosen to present it because it reinforces the parallel with R\'enyi and Sulanke's work while not relying on a new calculation and because it could be extended to the calculation of limiting variances.

\pass The first step consists in associating to any $x\in \cP_\lambda\setminus 2\cF_o(K)$ the point $y$ which is the closest point to $\partial K$ on the bisecting line of the segment $[0,x]$. An easy calculation shows that the new point process of such points $y$ has a local intensity near the boundary of $K$ of $4\lambda \langle s,n_s\rangle r_s^{-1}$ with respect to the coordinates $(s,h)$.

\noi In a second step, we rewrite the expectation $\E(\cN(K_\lambda))$ or $\E(\cA(K_\lambda))$ as a sum over all such points $y$ of the contribution of $y$. For instance, in the case of the functional $\cN(K_\lambda)$, this contribution called {\it score} is equal to $1$ if the line containing $y$ intersects the boundary of $K_\lambda$ and $0$ if not. We then apply Mecke-Slivnyak's theorem to rewrite it as an integral with respect to $(s,h)\in\partial K\times (0,\infty)$. For a fixed $s$, we apply the change of variables in the integral over $h$ provided by the inversion with respect to the osculating disk at $s$. In particular, it preserves locally the intensity of the point process and transforms, up to a negligible term, the score into the indicator function of a point being extreme with respect to a homogeneous Poisson point process inside $B_{r_s}(\omega_s)$ of intensity $4\lambda \langle s,n_s\rangle r_s^{-1}$. The integral of this new score with respect to $h$ is equal in turn, up to the multiplicative term $(2\pi r_s)^{-1}$, to the expected number of extreme points of that Poisson point process. This means we can apply point (iii) of Theorem A to the intensity $4\lambda \langle s,n_s\rangle r_s^{-1}$ and $K=B_{r_s}(o)$ and integrate the result divided by $2\pi r_s$ over $s\in\partial K$.

\pass The method for obtaining the limiting expectation of $\cA(K_\lambda)$ goes along the same lines, save for the fact that the functional obtained by inversion is the defect area of the Voronoi flower of a homogeneous Poisson point process inside $B_{r_s}(\omega_s)$. The asymptotics for the expectation of such area has been obtained by Schreiber in \cite{s02}.

\pass Following the ideas from \cite{cy14}, we would expect the method to provide in a similar way the limiting variance of $\cA(K_\lambda)$ and $\cN(K_\lambda)$. Up to technical justifications, we can claim in particular that up to a multiplicative constant not depending on $K$, they are equal to the respective limiting expectations, i.e. there exist two positive constants $C$ and $C'$ such that
\begin{equation*}
\mbox{Var}(\cA(K_\lambda))\underset{\lambda\to\infty}{\sim}C\lambda^{-\frac{2}{3}}\displaystyle\int_{\partial K} r_s^{\frac1{3}}\left\langle s,n_s\right\rangle^{-\frac{2}{3}}\dd s
\end{equation*}
and
\begin{equation*}
\mbox{Var}(\cN(K_{\lambda}))\underset{\lambda\to\infty}{\sim}C'\lambda^{\frac{1}{3}}\displaystyle\int_{\partial K} r_s^{-\frac{2}{3}}\left\langle s,n_s\right\rangle^{\frac1{3}}\dd s.
\end{equation*}
To the best of our knowledge, there is no easy way to extend the technique to the calculation of the expectation or variance of the perimeter of $K_\lambda$ when $K$ is smooth. Finally, when $K$ is a polygon, there is little hope to find a transformation which would play the role of the inversion as the asymptotic rates for the respective expected functionals in Theorems \ref{theopolygon} and B do not even coincide.

\section{Answer to Question 2: proof for the polygonal case}\label{sec:polygon}

\noi In this section $K$ is a convex polygon with vertices $a_1,a_2,\ldots,a_{n_K}$ containing the origin in its interior. Recall that $\alpha_i$ is
the interior angle at vertex $a_i$ and that $o_i$ is the orthogonal projection of $o$ onto the line $(a_i,a_{i+1})$. A point outside $K$ will be located by its polar coordinates from one vertex $a_i$ (see Figure \ref{fig:polygonsurf}), i.e. we consider a point $s_{a_i,\rho,\alpha}=a_i+\rho u_{\pi-\alpha}$ in the neighborhood of $a_i$, with $\rho>0$ and $\alpha\in(0,\alpha_i)$.

\noi The proof requires to decompose the set $\R^2\setminus K$ into several regions, namely $n_K$ cones above the vertices of $K$ and $n_K$  strips above the edges of $K$. More precisely, for every $1\le i \le n_K$, let us define
\begin{equation*}
\cG_i=\{s_{a_i,\rho,\alpha}:\rho>0 \text{ and } \alpha\in(\tfrac{\pi}{2},\tfrac{3\pi}{2}-\alpha_i)\}
\end{equation*}
and $$\cS_i=\{s_{a_i,\rho,\alpha}:\rho>0 \text{ and } \rho\cos\alpha\in(0,\|a_{i+1}-a_i\|)\}$$ the connected component of $\R^2\setminus (K\cup(\cup_{i=1}^{n_K}\cG_i))$ with $(a_i,a_{i+1})$ on its boundary.

\subsection{Increase of the area of the Voronoi flower}\label{subsec:fleurpoly}

\pass

\noi First, the following geometric lemma provides accurate estimates of the area of the set
\begin{equation*}\label{def:diffvorflowpoly}
\Delta\cF_{a_i,\rho,\alpha}=\cF_o(K \cup \{s_{a_i,\rho,\alpha}\})\setminus\cF_o(K).
\end{equation*}

\begin{lem}\label{lem:geompoly1}
Assume that $K$ is a convex polygon and let $a_i\in\partial K$, $1\leq i \leq n_K$, a fixed vertex of $K$. 
\begin{enumerate}
\item[$(i)$] We get, uniformly in $\rho$,
\begin{equation*}\label{estimeevolfleurpoly1}
  \cA(\Delta\cF_{a_i,\rho,\alpha})\indicat_{\{s_{a_i,\rho,\alpha}\in \cS_i\}}\underset{\alpha\to0}{\sim} \alpha^2\frac{\|o_i\|}{2}\frac{\rho\|a_{i+1}-a_i\|}{\|a_{i+1}-a_i\|-\rho}.
\end{equation*}
\item[$(ii)$] Moreover, there exists $C>0$ such that, for all $(\rho,\alpha)$ such that $s_{a_i,\rho,\alpha}$ belongs to $\cG_i\cup\cS_i$,
\begin{equation}\label{estimeevolfleurpoly2}
\cA(\Delta\cF_{a_i,\rho,\alpha})\ge C\max(1,\rho)\rho\alpha^2.
\end{equation}
\end{enumerate}
\end{lem}

\pass\textbf{Proof of Lemma \ref{lem:geompoly1}.}

\begin{figure}[h!]
\includegraphics[scale=0.8]{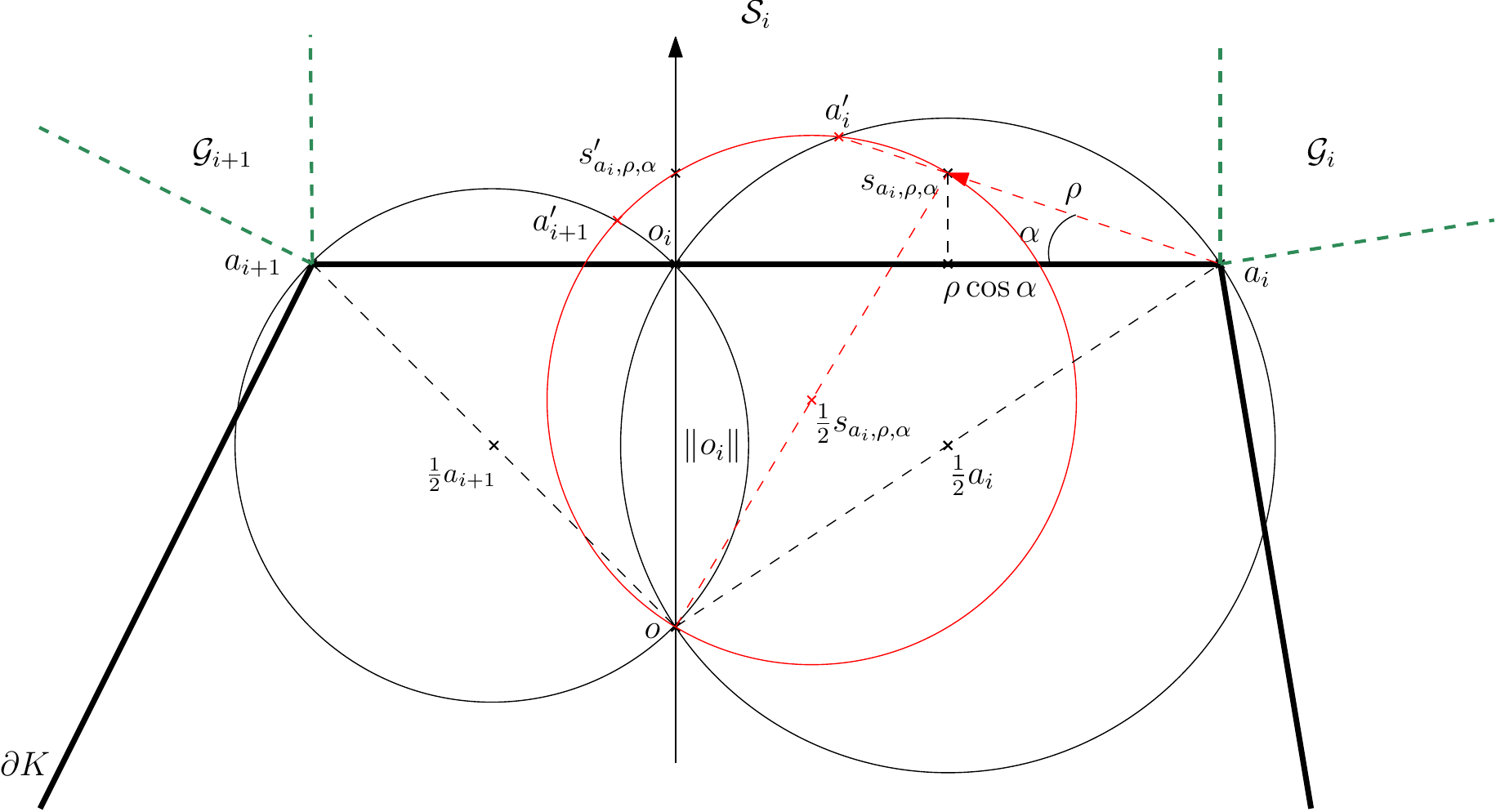}
\caption{\label{fig:polygonsurf} The analogue of Figure \ref{fig:geolem-proof1} in the polygonal case.}
\end{figure}

\noi{\it Proof of $(i)$.}
For $\alpha$ small enough, the  set $\Delta \cF_{a_i,\rho,\alpha}$ is nothing but
\begin{equation*}\label{eq:fleurdecompdisques}
\Delta \cF_{a_i,\rho,\alpha}=B_{\frac1{2}\|s_{a_i,\rho,\alpha}\|}(\tfrac1{2} s_{a_i,\rho,\alpha})\setminus \left(B_{\frac1{2}\|a_i\|}(\tfrac{1}{2}a_i)\cup B_{\frac{1}{2}\|a_{i+1}\|}(\tfrac1{2}a_{i+1})\right),
\end{equation*}
that is a curvilinear triangle with vertices $o_i$, $a_i'$ and $a_{i+1}'$ where $a_i'$ and $a_{i+1}'$ are respectively the intersection
of $\partial B_{\frac1{2}\|s_{a_i,\rho,\alpha}\|}(\tfrac1{2} s_{a_i,\rho,\alpha})$ with $\partial B_{\frac1{2}\|a_i\|}(\tfrac1{2}a_i)$ and $\partial B_{\frac1{2}\|a_{i+1}\|}(\tfrac1{2}a_{i+1})$ (see Figure \ref{fig:polygonsurf}).

\pass We aim at computing estimates for the area $\cA(\Delta\cF_{a_i,\rho,\alpha})$ of this curvilinear triangle. To do this, we split it into the
curvilinear triangles with vertices $o_i$, $s_{a_i,\rho,\alpha}'$, $a_i'$ and $o_i$, $s_{a_i,\rho,\alpha}'$, $a_{i+1}'$ respectively, where $s_{a_i,\rho,\alpha}'$ is the intersection of the line $(o,o_i)$ with the circle $\partial B_{\frac1{2}\|s_{a_i,\rho,\alpha}\|}(\tfrac1{2} s_{a_i,\rho,\alpha})$.

\pass Let us focus on the first curvilinear triangle. As $\alpha\to0$, it tends to a straight triangle whose area is given by
\begin{equation*}
\cA(o_i s_{a_i,\rho,\alpha}' a_i')\underset{\alpha\to 0}{\sim} \frac1{2}\|s'_{a_i,\rho,\alpha}-o_i\|(\|\overset{\frown}{a_i'o_i}\|\sin\beta_i).
\end{equation*}
where $\beta_i$ is the angle between the line $(o_i,s'_{a_i,\rho,\alpha})$ and the tangent line to the disk $B_{\frac{1}{2}a_i}(\tfrac{1}{2}\|a_i\|)$ at $o_i$.

\pass Observe first that we get by symmetry
\begin{equation*}
\|s'_{a_i,\rho,\alpha}-o_i\| = \rho\sin \alpha.
\end{equation*}

\pass Let us now compute the length of the arc $\overset{\frown}{a_i'o_i}$. Observe now that the lines $(s_{a_i,\rho,\alpha},a_i')$ and $(a_i,a_i')$ are both perpendicular to $(o,a_i')$. Therefore the points $a_i$, $s_{a_i,\rho,\alpha}$ and $a_i'$ are aligned. It follows that the angle between the lines $(a_i,o_i)$ and $(a_i,a_i')$ is the same as the angle between $(a_i,o_i)$ and $(a_i,s_{a_i,\rho,\alpha})$ which is nothing but $\alpha$. Thus the central angle between the lines  $(\frac{1}{2}a_i,a_i')$ and $(\frac{1}{2}a_i,o_i)$ is $2\alpha$. Since the arc $\overset{\frown}{a_i'o_i}$ belongs to the circle with center $\tfrac{1}{2}a_i$ and diameter $\|a_i\|$ we get
\begin{equation*}
\|\overset{\frown}{a_i'o_i}\| = \frac1{2}\|a_i\|(2\alpha) = \|a_i\|\alpha.
\end{equation*}

\noi Finally, observing that $\beta_i$ is equal to the angle between $(a_i,o_i)$ and $(a_i,o)$ by the inscribed angle theorem,
we deduce that
\begin{align*}
\cA(o_i s_{a_i,\rho,\alpha}' a_i') \underset{\alpha\to 0}{\sim} \frac1{2}\Big(\rho\sin \alpha\Big)\Big(\|a_i\|\alpha\sin\beta_i \Big)
 = \frac1{2}\Big(\rho\sin \alpha\Big)\Big(\|a_i\|\alpha \frac{\|o_i\|}{\|a_i\|}\Big)
 \underset{\alpha\to0}{\sim} \frac1{2}\rho \alpha^2\|o_i\|.
\end{align*}

\pass Now, computing similarly the area $\cA(o_i s_{a_i,\rho,\alpha}' a_{i+1}')$ of the other curvilinear triangle, we
obtain
\begin{align*}
\cA(o_i s_{a_i,\rho,\alpha}' a_{i+1}') = \frac1{2}\rho' \alpha'^2\|o_i\|
\end{align*}
where,
\begin{equation*}
\Bigg\{
  \begin{array}{lll}
    \rho'   & = \Big((\|a_{i+1}-a_i\|-\rho\cos\alpha)^2+(\rho\sin\alpha)^2\Big)^{\frac1{2}} & \underset{\alpha\to 0}{\sim}  \|a_{i+1}-a_i\|-\rho \\
    \alpha' & = \arctan\Big(\frac{\rho\sin\alpha}{\|a_{i+1}-a_i\|-\rho\cos\alpha}\Big)     & \underset{\alpha\to 0}{\sim} \frac{\rho\alpha}{\|a_{i+1}-a_i\|-\rho}
  \end{array}
\Bigg.
\end{equation*}

\noi Therefore,
\begin{align*}
\cA(o_i s_{a_i,\rho,\alpha}' a_{i+1}') \underset{\alpha\to 0}{\sim} \frac1{2}\frac{\rho^2\alpha^2\|o_i\|}{\|a_{i+1}-a_i\|-\rho}.
\end{align*}

\pass Finally, summing the areas of each triangle, we obtain
\begin{align*}
\cA(\Delta\cF_{a_i,\rho,\alpha})\underset{\alpha\to0}{\sim} \frac1{2}\rho\alpha^2\|o_i\|\bigg(1+\frac{\rho}{\|a_{i+1}-a_i\|-\rho}\bigg)
\underset{\alpha\to0}{\sim} \frac1{2}\rho\alpha^2\|o_i\|\frac{\|a_{i+1}-a_i\|}{\|a_{i+1}-a_i\|-\rho}.
\end{align*}
so that (i) holds.

\pass{\it Proof of $(ii)$}.
We first assume that $s_{a_i,\rho,\alpha}\in \cS_i$. Because of point (i), there exists $\alpha_0\in (0,\tfrac{\pi}{2})$ such that for every $\alpha<\alpha_0$,
$\cA(\Delta\cF_{a_i,\rho,\alpha})\ge \mfrac{\|o_i\|}{4}\rho\alpha^2$ and $\rho\le \mfrac{\|a_{i+1}-a_i\|}{\cos\alpha_0}$.
This shows \eqref{estimeevolfleurpoly2} as soon as $\alpha<\alpha_0$. When $\alpha>\alpha_0$, it is enough to show that $ \cA(\Delta\cF_{a_i,\rho,\alpha})\ge C \max (\rho,\rho^2)$
for some positive constant $C$. This last inequality comes now from the fact that $\Delta\cF_{a_i,\rho,\alpha}$ contains both a disk of radius proportional to $\rho$ and an angular sector with thickness $\rho$ and constant angular width. Finally, the exact same argument holds when $s_{a_i,\rho,\alpha}\in \cG_i$ so this completes the proof. $\square$

\subsection{Proof of Theorem \ref{theopolygon} (i): the defect area}\label{subsec:surfacepolygon}
\pass

\pass
We can write, recalling the notation of Section \ref{subsec:fleurpoly},
\begin{align*}
\E(\cA(K_{\lambda})) & -\cA(K) \\
& = \sum_{i=1}^{n_K}\int_{\cS_i}\exp\left(-4\lambda\cA(\Delta\cF_{a_i,\rho,\alpha})\right)\rho\dd \rho \dd\alpha
       + \sum_{i=1}^{n_K}\int_{\cG_i}\exp\left(-4\lambda\cA(\Delta\cF_{a_i,\rho,\alpha})\right) \rho\dd \rho \dd \alpha .
\end{align*}

\pass It is then enough to show that, for every $1\le i\le n_K$,
\begin{equation}\label{eq:madonna}
\lambda^{\frac1{2}}\int_{\cS_i}\exp\left(-4\lambda\cA(\Delta\cF_{a_i,\rho,\alpha})\right)\rho\dd \rho \dd\alpha
\underset{\lambda\to\infty}{\longrightarrow} 2^{-\frac{9}{2}}\pi^{\frac{3}{2}}\|o_i\|^{-\frac1{2}}\|a_{i+1}-a_i\|^{\frac{3}{2}}
\end{equation}
and
\begin{equation}\label{eq:rayoflight}
  \lambda^{\frac{1}{2}}
  \int_{\cG_i}\exp\left(-4\lambda\cA(\Delta\cF_{a_i,\rho,\alpha})\right) \rho\dd \rho \dd \alpha
\underset{\lambda\to\infty}{\longrightarrow} 0.
\end{equation}

\pass{\it Proof of \eqref{eq:madonna}}.
Let us fix $1\le i\le n_K$. The change of variables $\beta=\lambda^{\frac1{2}}\alpha$ yields
\begin{align*}
\lambda^{\frac1{2}}\int_{\cS_i}\exp & \left(-4\lambda \cA(\Delta\cF_{a_i,\rho,\alpha})\right)\rho\dd \rho \dd\alpha \\
& = \int_{0}^{\rho_i(\lambda^{-1/2}\beta)}\int_{0}^{\lambda^{1/2}\frac{\pi}{2}}
\exp\left(-4\lambda\cA(\Delta\cF_{a_i,\rho,\lambda^{-1/2}\beta})\right)\rho\dd \rho \dd\beta
\end{align*}
where $\rho_i(\cdot)$ denotes the equation of the line containing $a_{i+1}$ and orthogonal to $(a_{i},a_{i+1})$ with respect to the polar coordinates $(\rho,\alpha)$.

\pass Thanks to Lemma \ref{lem:geompoly1}, we have
\begin{align*}
4\lambda\cA(\Delta\cF_{a_i,\rho,\lambda^{-1/2}\beta}) \underset{\lambda\to\infty}{\longrightarrow}
2\|o_i\|\frac{\rho\|a_{i+1}-a_i\|}{\|a_{i+1}-a_i\|-\rho}\beta^2 .
\end{align*}
and, for all $\lambda>0$,
\begin{equation*}
\lambda\cA(\Delta\cF_{a_i,\rho,\lambda^{-1/2}\beta}) \ge C \max(1,\rho)\rho \beta^2
\end{equation*}
where $C$ is a positive constant.

\pass Consequently, we can apply Lebesgue's dominated convergence theorem to obtain
\begin{align*}
 & \int_0^{\rho_i(\lambda^{-1/2}\beta)}\int_{0}^{\lambda^{1/2}\frac{\pi}{2}}
\exp\left(-4\lambda\cA(\Delta\cF_{a_i,\rho,\lambda^{-1/2}\beta})\right)\rho\dd \rho \dd\beta \\
 & \hspace{3cm} \underset{\lambda\to\infty}{\longrightarrow} \int_0^{\|a_{i+1}-a_i\|}\int_0^{\infty}
\exp\left(-2\|o_i\|\frac{\rho\|a_{i+1}-a_i\|}{\|a_{i+1}-a_i\|-\rho}\beta^2 \right)\rho\dd \rho \dd\beta \\
 & \hspace{3.22cm}  = \|a_{i+1}-a_i\|^2\int_0^1\bigg(\int_0^{\infty} \exp\left(-2\|o_i\|\frac{l\|a_{i+1}-a_i\|}{1-l}\beta^2 \right)\dd\beta \bigg) l\dd \l \\
 & \hspace{3.22cm} = \|a_{i+1}-a_i\|^{\frac{3}{2}} \pi^{\frac1{2}}2^{-\frac{3}{2}}\|o_i\|^{-\frac1{2}} \int_0^1 (l(1-l))^{\frac1{2}}\dd l \\
 & \hspace{3.22cm} = \|a_{i+1}-a_i\|^{\frac{3}{2}} \pi^{\frac{3}{2}}2^{-\frac{9}{2}}\|o_i\|^{-\frac1{2}}.
\end{align*}

\pass{\it Proof of \eqref{eq:rayoflight}}.
Using the second part of Lemma \ref{lem:geompoly1} we have successively, for all $\lambda>0$,
\begin{align*}
& \int_{\cG_i}\exp\left(-4\lambda\cA(\Delta\cF_{a_i,\rho,\alpha})\right) \rho\dd \rho \dd \alpha \\
& \hspace{1cm} \leq \int_0^{\infty}\int_{\frac{\pi}{2}}^{\frac{3\pi}{2}-\alpha_i}\exp\left(-C\lambda\max(1,\rho)\rho\alpha^2\right) \rho\dd \rho \dd \alpha  \\
& \hspace{1cm} = \int_{\frac{\pi}{2}}^{\frac{3\pi}{2}-\alpha_i}\bigg(\int_0^1\exp\left(-C\lambda\rho\alpha^2\right) \rho\dd \rho\bigg) \dd \alpha
     + \int_{\frac{\pi}{2}}^{\frac{3\pi}{2}-\alpha_i}\bigg(\int_1^{\infty}\exp\left(-C\lambda\rho^2\alpha^2\right) \rho\dd \rho\bigg) \dd \alpha  \\
& \hspace{1cm} \leq \int_{\frac{\pi}{2}}^{\frac{3\pi}{2}-\alpha_i} \lambda^{-2}\alpha^{-4} \dd \alpha
     + \int_{\frac{\pi}{2}}^{\frac{3\pi}{2}-\alpha_i} \lambda^{-1}\alpha^{-2} \dd \alpha  \\
& \hspace{1cm} \leq C \lambda^{-1}.
\end{align*}

\pass That implies \eqref{eq:rayoflight} and completes the proof. $\square$

\subsection{Proof of Theorem \ref{theopolygon} (ii): support points and defect perimeter}\label{sub:perimeterpolygon}

\pass

\noi The proof relies on \eqref{eq:croftonperimeter}, i.e. the rewriting of the mean defect perimeter of $K_\lambda$ as the integral in all directions of the mean defect support function of $K_\lambda$. Let us denote by $u_{\delta_i}$, $\delta_i\in (0,2\pi)$, the external unit normal vector to the line $(a_{i},a_{i+1})$ (with the convention $\delta_0=\delta_{n_K}$). We expect $\E(p_o(K_\lambda,\theta)-p_o(K,\theta))$ to be maximal in directions {\it close} to $\delta_i$ for every $i$ while the remaining directions should have a negligible contribution inside the integral on the right-hand side of \eqref{eq:croftonperimeter}. Let us first rewrite \eqref{eq:croftonperimeter} as
\begin{align*}
\E( &\cU(K_{\lambda})) -\cU(K) \\
& = \sum_{i=1}^{n_K} \bigg(\int_{\delta_{i-1}}^{\frac1{2}(\delta_{i-1}+\delta_{i})} \E(p_o(K_{\lambda},\theta)-p_o(K,\theta))\dd\theta
   + \int_{\frac1{2}(\delta_{i-1}+\delta_{i})}^{\delta_{i}} \E(p_o(K_{\lambda},\theta)-p_o(K,\theta))\dd\theta\bigg).
\end{align*}

\pass We estimate only one integral $\int_{\frac1{2}(\delta_{i-1}+\delta_i)}^{\delta_i} \E(p_o(K_{\lambda},\theta)-p_o(K,\theta))\dd\theta$
as the others will be treated analogously. Using the change of variable $\theta=\delta_i-\lambda^{-\gamma}$, we get
\begin{align}\label{eq:perimpolyginterm0}
\frac{\lambda}{\log \lambda}&\int_{\frac1{2}(\delta_{i-1}+\delta_i)}^{\delta_i} \E(p_o(K_{\lambda},\theta)-p_o(K,\theta))\dd\theta \nonumber \\
& \underset{\lambda\to\infty}{\sim} \frac{\lambda}{\log \lambda}\int_0^{\infty} \E(p_o(K_{\lambda},\theta_i-\lambda^{-\gamma})-p_o(K,\delta_i-\lambda^{-\gamma}))(\lambda^{-\gamma}\log \lambda)\dd\gamma.
\end{align}

\pass It is now a consequence of the next two Propositions that only directions $\delta_i-\lambda^{-\gamma}$ up to the critical value $\gamma=\tfrac1{2}$ will contribute. Precisely, we deduce from point (ii) of Proposition \ref{prop:supportpoly} below combined with Lebesgue's dominated convergence theorem that
\begin{equation}\label{eq:perimpolyginterm1}
\frac{\lambda}{\log \lambda}\int_0^{\frac1{2}} \E(p_o(K_{\lambda},\theta_i-\lambda^{-\gamma})-p_o(K,\delta_i-\lambda^{-\gamma}))(\lambda^{-\gamma}\log \lambda)\dd\gamma \underset{\lambda\to\infty}{\longrightarrow} \frac1{12\|o_i\|}
\end{equation}
and from point (ii) of Proposition \ref{prop:supportpoly2} that
\begin{equation}\label{eq:perimpolyginterm2}
  \lambda\int_{\frac1{2}}^{\infty} \E(p_o(K_{\lambda},\theta_i-\lambda^{-\gamma})-p_o(K,\delta_i-\lambda^{-\gamma}))\lambda^{-\gamma}\dd\gamma\le \lambda\int_{\frac1{2}}^{\infty} C\lambda^{-\frac1{2}} \lambda^{-\gamma}\dd \gamma \underset{\lambda\to\infty}{\longrightarrow} 0.
\end{equation}

\pass Inserting \eqref{eq:perimpolyginterm1} and \eqref{eq:perimpolyginterm2} into \eqref{eq:perimpolyginterm0} shows that
\begin{align*}
\int_{\frac1{2}(\delta_{i-1}+\delta_i)}^{\delta_i} \E(p_o(K_{\lambda},\theta)-p_o(K,\theta))\dd\theta
& \underset{\lambda\to\infty}{\sim} \int_{\delta_{i-1}}^{\frac1{2}(\delta_{i-1}+\delta_{i})} \E(p_o(K_{\lambda},\theta)-p_o(K,\theta))\dd\theta \\
& \underset{\lambda\to\infty}{\sim}\frac{\log \lambda}{\lambda}\frac1{12\|o_i\|}.
\end{align*}
Summing these estimates over the vertices of $K$ provides then the required result. $\square$

\pass As in the smooth case, we need now to explain how to estimate the defect support function of $K_\lambda$ in a fixed direction $\delta_i-\lambda^{-\gamma}$. As emphasized in the proof of point (ii) of Theorem \ref{theopolygon}, we will treat separately the cases $\gamma\in(0,\tfrac1{2})$ and $\gamma\ge \tfrac1{2}$. More precisely, let us introduce the support point $m_{\lambda^{-\gamma}}$ on $\partial K_{\lambda}$ which satisfies
\begin{equation*}
\langle m_{\lambda^{-\gamma}}, u_{\delta_i-\lambda^{-\gamma}}\rangle=p_o(K_{\lambda},u_{\delta_i-\lambda^{-\gamma}})
\end{equation*}
and let us denote by $(R_{\lambda^{-\gamma}},A_{\lambda^{-\gamma}})$ the polar coordinates of $m_{\lambda^{-\gamma}}$ with respect to the coordinate system with origin $a_i$ and first axis $(a_i,a_{i+1})$. In particular, we notice that $A_{\lambda^{-\gamma}}\ge \lambda^{-\gamma}$ almost surely since $p_o(K_{\lambda},u_{\delta_i-\lambda^{-\gamma}})\ge p_o(K,u_{\delta_i-\lambda^{-\gamma}})\ge \langle a_i,u_{\delta_i-\lambda^{-\gamma}}\rangle$.

\pass The next proposition investigates the asymptotic distribution of the couple $(R_{\lambda^{-\gamma}},A_{\lambda^{-\gamma}})$ for $\gamma\in (0,\tfrac12)$.

\begin{prop}\label{prop:supportpoly}
\noi
\begin{enumerate}
\item[$(i)$] For every $\gamma\in (0,\tfrac{1}{2})$, the couple $(\lambda^{1-2\gamma}R_{\lambda^{-\gamma}}, \lambda^{\gamma}A_{\lambda^{-\gamma}})$
converges in distribution when $\lambda\to\infty$ to the distribution with density function $f_i$ given by
\begin{equation*}\label{eq:fthetapoly}
f_i(\rho,\alpha)= 8\|o_i\|^2\exp\left(-2\|o_i\|\rho\alpha^2\right)\alpha(\alpha-1)\rho \indicat_{\{\rho>0\}}\indicat_{\{\alpha>1\}}.
\end{equation*}
\item[$(ii)$] There exists $C>0$ such that for every $\gamma\in (0,\tfrac{1}{2})$ and $\lambda>0$,
\begin{equation}\label{eq:estimsupppolyg}
  \lambda^{1-\gamma}\E(R_{\lambda^{-\gamma}}\sin(A_{\lambda^{-\gamma}}-\lambda^{-\gamma}))\le C.
\end{equation}
Moreover, for every $\gamma\in (0,\tfrac{1}{2})$,
\begin{equation*}
\hspace{1.5cm} \E(R_{\lambda^{-\gamma}}\sin(A_{\lambda^{-\gamma}}-\lambda^{-\gamma}))
= \E(p_o(K_{\lambda},\delta_i -\lambda^{-\gamma})-p_o(K,\delta_i-\lambda^{-\gamma}))\underset{\lambda\to\infty}{\sim}\lambda^{\gamma-1}\frac1{6\|o_i\|}.
\end{equation*}
\end{enumerate}
\end{prop}

\pass\textbf{Proof of Proposition \ref{prop:supportpoly}.}

\begin{figure}[h!]
\includegraphics[scale=0.8]{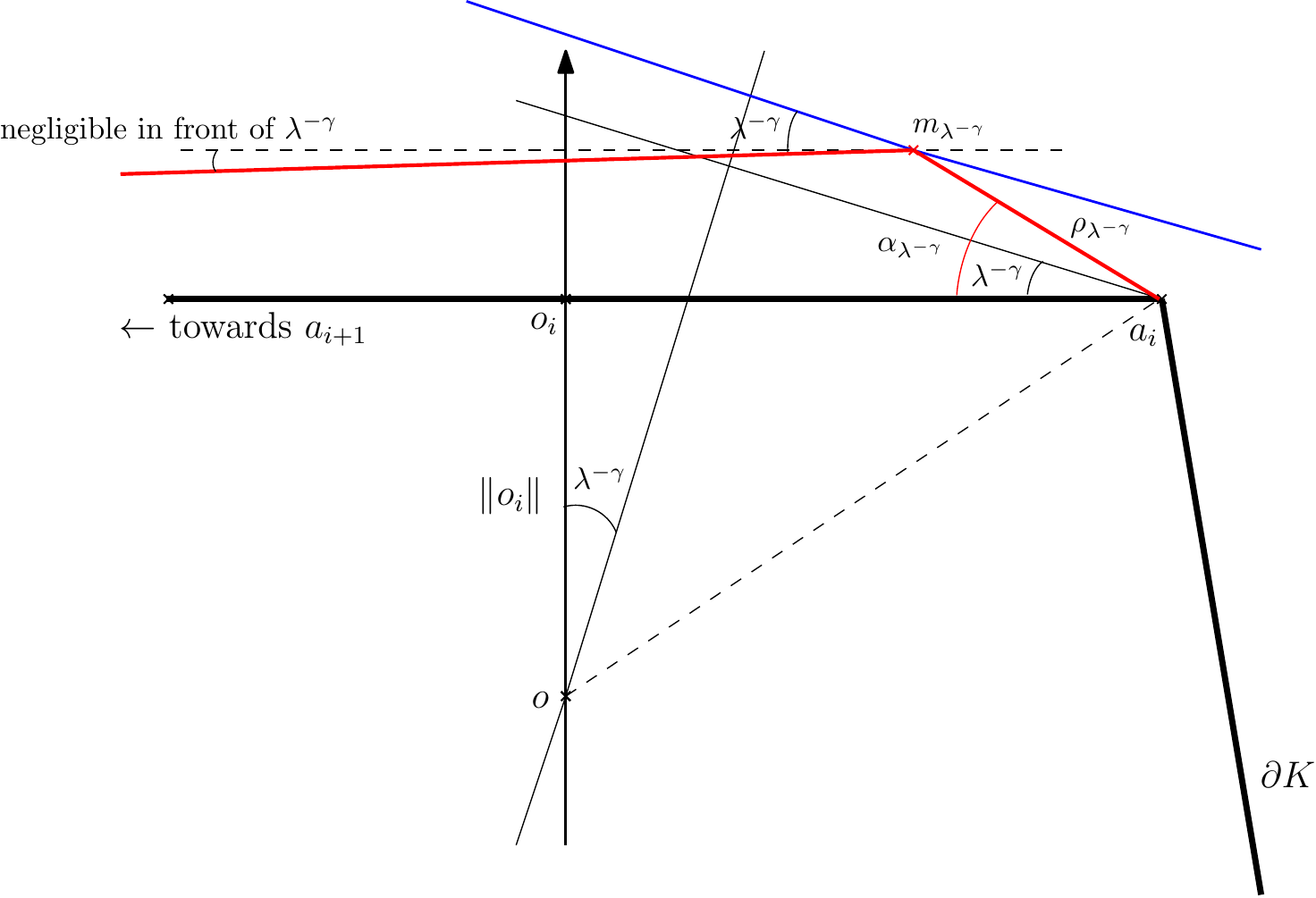}
\caption{\label{fig:supportpolygon} The analogue of Figure \ref{fig:supportsmooth} in the polygonal case: $m_{\lambda^{-\gamma}}$ denotes the support point in the direction $\mfrac{\pi}{2}-\lambda^{-\gamma}$.}
\end{figure}

\noi {\it Proof of $(i)$.} Without loss of generality, we can assume in the proof that $\delta_i=\tfrac{\pi}{2}$. Once again, the strategy of the proof consists in going along the same lines as for the smooth case. We start with the same identity but written in polar coordinates, that is
$c_{x_1,x_2}=m_{\lambda^{-\gamma}}= a_i+ru_{\theta}$ (see Figure \ref{fig:supportpolygon}). Notice that $(r,\theta)\in\cS_i$ when $\lambda\to\infty$.

\pass We then proceed with two consecutive changes of variables:
\begin{enumerate}[label=-,leftmargin=1\parindent]
\item denoting by $\theta_1$ and $\theta_2$ the angles between the two bisecting lines emanating from $c_{x_1,x_2}$ corresponding to the right and left neighbor of $o$ respectively, we use Lemma \ref{lem:jacob}.
\item secondly, we replace the couple $(r,\theta)$ by $(\rho,\alpha)$ defined by
\begin{align*}
\rho = \lambda^{1-2\gamma}r \,\text{ and }\, \alpha=\lambda^{\gamma}\theta.
\end{align*}

\noi We get in particular
\begin{equation*}
r=r_i(\lambda,\rho,\alpha)=\|s_{a_i,\lambda^{2\gamma-1}\rho,\lambda^{-\gamma}\alpha}\|
\end{equation*}
and
\begin{equation*}
\theta=\theta_i(\lambda,\rho,\alpha)
= \arcsin\bigg(\bigg\langle \frac{s_{a_i,\lambda^{2\gamma-1}\rho,\lambda^{-\gamma}\alpha}}{\|s_{a_i,\lambda^{2\gamma-1}\rho,\lambda^{-\gamma}\alpha}\|}\bigg\rangle\bigg)
\end{equation*}
whose Jacobian is given by
\begin{equation*}
\frac{r\dd r \dd \theta }{\rho\dd\rho \dd\alpha} = \lambda^{3\gamma-2}
\end{equation*}
\end{enumerate}

\pass Consequently, as in the proof of point (i) of Proposition \ref{prop:supportsmooth}, we deduce that  for every positive and measurable function $\varphi:\R^2\longrightarrow \R_+$
\begin{align*}
\E(\varphi(\lambda^{1-2\gamma}R_{\lambda^{-\gamma}} &, \lambda^{\gamma}A_{\lambda^{-\gamma}})) \\
& = 16\int_1^{\lambda^{\gamma}\pi+1}\int_0^{\lambda^{1-2\gamma}\rho_i(\alpha)}\exp(-\Delta_i(\lambda,\rho,\alpha))r_i^2(\lambda,\rho,\alpha)
  J^{\mbox{\scalebox{.5}{supp}}}_i(\lambda,\rho,\alpha)\varphi(\rho,\alpha)\rho\dd\alpha\dd \rho  \\
\end{align*}
where
\begin{equation*}
  \Delta_i(\lambda,\rho,\alpha)= 4\lambda\cA\big(B_{\|s_{a_i,\lambda^{2\gamma-1}\rho,\lambda^{-\gamma}\alpha}\|}(s_{a_i,\lambda^{2\gamma-1}\rho,\lambda^{-\gamma}\alpha})\setminus \cF_o(K)\big)
  = 4\lambda\cA\big(\Delta \cF_{a_i,\lambda^{2\gamma-1}\rho,\lambda^{-\gamma}\alpha}\big)
\end{equation*}
and
\begin{equation*}
J^{\mbox{\scalebox{.5}{supp}}}_i(\lambda,\rho,\alpha)=\lambda^{3\gamma}\int_{E_i(\lambda,\rho,\alpha)} J(\theta_i(\lambda,\rho,\alpha),\theta_1,\theta_2)\dd \theta_1\dd \theta_2
\end{equation*}
where $E_i^{\mbox{\scalebox{.5}{supp}}}(\lambda,\rho,\alpha)$ stands for the set of couples $(\theta_1,\theta_2)$ which satisfy that the two bisecting lines of $[o,x_1]$ and $[o,x_2]$ do not intersect $K$.

\pass Thanks to Lemma \ref{lem:geompoly1}, we have
\begin{align*}
\Delta_i(\lambda,\rho,\alpha) \underset{\lambda\to\infty}{\longrightarrow} 2\|o_i\|\rho\alpha^2
\end{align*}
as soon as $\lambda^{2\gamma-1}\rho\cos(\lambda^{-\gamma}\alpha)\le \|a_{i+1}-a_i\|$ and
\begin{equation}\label{eq:dominDeltapolyg}
  \Delta_i(\lambda,\rho,\alpha)\ge C \rho \alpha^2
\end{equation}
for some constant $C>0$. Moreover,
\begin{align*}
r_i(\lambda,\rho,\alpha) \underset{\lambda\to\infty}{\longrightarrow} \|a_i\|
\end{align*}
and for $\lambda\ge 1$,
\begin{equation}\label{eq:dominripolyg}
r_i(\lambda,\rho,\alpha)\le \|a_i\|+\lambda^{2\gamma-1}\rho\le \|a_i\|+\rho.
\end{equation}

\pass Let us now turn on the term $J^{\mbox{\scalebox{.5}{supp}}}_i(\lambda,\rho,\alpha)$. Using the convergence

\begin{equation}
\sin(\theta_i(\lambda,\rho,\alpha))\underset{\lambda\to\infty}{\longrightarrow} \frac{\|o_i\|}{\|a_i\|},
\end{equation}
we get successively
\begin{align*}
& J^{\mbox{\scalebox{.5}{supp}}}_i(\lambda,\rho,\alpha) \\
  & \underset{\lambda\to\infty}{\sim} \lambda^{3\gamma}\int_{E_i(\lambda,\rho,\alpha)}J(\theta_i(\lambda,\rho,\alpha),\theta_1,\theta_2)\dd \theta_1\dd \theta_2 \\
  & \underset{\lambda\to\infty}{\sim} \lambda^{3\gamma}\int_{(0,(\alpha-1)\lambda^{-\gamma})\times(0,\lambda^{-\gamma})} |\sin(\theta_i(\lambda,\rho,\alpha)-\theta_1)\sin(\theta_i(\lambda,\rho,\alpha)-\theta_2)\sin(\theta_2-\theta_1)|\dd \theta_1\dd \theta_2 \\
  & \underset{\lambda\to\infty}{\sim} \lambda^{3\gamma}\Big(\frac{\|o_i\|}{\|a_i\|}\Big)^2 \Big(\mfrac1{2}\lambda^{-3\gamma}\alpha(\alpha-1)\Big) \\
  &  \hspace{0.22cm} = \frac1{2}\Big(\frac{\|o_i\|}{\|a_i\|}\Big)^2\alpha(\alpha-1).
\end{align*}

\noi Finally, we notice that any couple $(\theta_1,\theta_2)\in E_i(\lambda,\rho,\alpha)$ satisfies that one of the two angles is at most equal to $\alpha$ and the other to $\rho\alpha$ up to a multiplicative constant. Consequently,  we can show that for some constant $C>0$,
\begin{equation}\label{eq:dominJpolyg}
  J^{\mbox{\scalebox{.5}{supp}}}_i(\lambda,\rho,\alpha)\le C\rho\alpha^2.
\end{equation}

\pass A method based on Lebesgue's dominated convergence theorem and analogous to the proof of Proposition \ref{prop:supportsmooth} will show that
\begin{align*}
& \E(\varphi(\lambda^{1-2\gamma}R_{\lambda^{-\gamma}}, \lambda^{\gamma}A_{\lambda^{-\gamma}}))\\
&\hspace*{1cm}\underset{\lambda\to\infty}{\sim} 16\int_{(0,\infty)^2}\exp(2\|o_i\|\rho\alpha^2)\|a_i\|^2(\mfrac1{2}\big(\mfrac{\|o_i\|}{\|a_i\|}\big)^2\alpha(\alpha-1)\indicat_{\{\alpha>1\}})
\varphi(\rho,\alpha)\rho\dd \rho \dd\alpha.
\end{align*}
This implies the required result.

\pass{\it Proof of $(ii)$.}
This is a direct consequence of the convergence in distribution proved in (i) and of the equality
\begin{equation*}
p_o(K_{\lambda},\delta_i -\lambda^{-\gamma})-p_o(K,\delta_i-\lambda^{-\gamma})=R_{\lambda^{-\gamma}}\sin(A_{\lambda^{-\gamma}}-\lambda^{-\gamma}).
\end{equation*}

\noi Indeed, applying the method used in (i) to $\varphi(\rho,\alpha)=\rho(\alpha-1)$, we get \eqref{eq:estimsupppolyg} from \eqref{eq:dominDeltapolyg}, \eqref{eq:dominripolyg}, \eqref{eq:dominJpolyg} and the inequality $\sin(x)\le x$ for $x>0$. Moreover,
\begin{align*}
\lambda^{1-\gamma}(\E(p_o(K_{\lambda},\tfrac{\pi}{2} -\lambda^{-\gamma})-p_o(K,\tfrac{\pi}{2}-\lambda^{-\gamma})))
 & \underset{\lambda\to\infty}{\sim} \E(\lambda^{1-2\gamma}R_{\lambda^{-\gamma}}(\lambda^{\gamma}A_{\lambda^{-\gamma}}-1)) \\
 & \hspace{0.22cm} = \int_{\R^2} \rho(\alpha-1) f_i(\rho,\alpha)\dd \rho\dd \alpha \\
 & \hspace{0.22cm} = \frac1{6\|o_i\|}.\quad\square
\end{align*}

\pass

\noi We will need an analogous result for the support function of $K$ with respect to $o$ in a direction of the form $\tau\lambda^{-\frac1{2}}$, for $\tau\ge0$. Let us introduce the point $m_{\tau\lambda^{-1/2}}$ on $\partial K_{\lambda}$ which satisfies
\begin{equation*}
\langle m_{\tau\lambda^{-1/2}}, u_{\frac{\pi}{2}-\tau\lambda^{-1/2}}\rangle=p_o(K_{\lambda},u_{\frac{\pi}{2}-\tau\lambda^{-1/2}})
\end{equation*}
and denote by $R_{\tau\lambda^{-1/2}}$ and $A_{\tau\lambda^{-1/2}}$ the polar coordinates of $m_{\tau\lambda^{-1/2}}$ with respect to $a_i$.

\pass The next proposition provides the asymptotic distribution of the couple $(R_{\tau\lambda^{-1/2}},A_{\tau\lambda^{-1/2}})$. Since it is very similar to Proposition \ref{prop:supportpoly}, the proof is omitted.

\newpage

\begin{prop}\label{prop:supportpoly2}
\noi
\begin{enumerate}
\item[$(i)$] Let $\tau\ge 0$. The couple $(R_{\tau\lambda^{-1/2}}, \lambda^{-\frac1{2}}A_{\tau\lambda^{-1/2}})$ converges in distribution when $\lambda\to\infty$ to the distribution with density function $g_i$ given by
\begin{align*}
\hspace{1.5cm} g_i(\rho,\alpha)
& = 8\|o_i\|^2\frac{\rho\|a_{i+1}-a_i\|}{\|a_{i+1}-a_i\|-\rho}\exp\left(-2\|o_i\|\frac{\rho\|a_{i+1}-a_i\|}{\|a_{i+1}-a_i\|-\rho}\alpha^2\right) \\
& \hspace{2cm} \times (\alpha-\tau)\Big(\frac{\rho\alpha}{\|a_{i+1}-a_i\|-\rho}+\tau\Big)\indicat_{\{\rho\in(0,\|a_{i+1}-a_i\|)\}}\indicat_{\{\alpha>\tau\}}.
\end{align*}
\item[$(ii)$] There exists a positive constant $\tau\ge 0$ such that, for every $\gamma>\tfrac1{2}$,
\begin{align*}
\E(p_o(K_{\lambda},\delta_i-\lambda^{-\gamma})-p_o(K,\delta_i-\lambda^{-\gamma})) \leq \tau\lambda^{-\frac1{2}}.
\end{align*}
\end{enumerate}
\end{prop}

\pass Let us notice that the special case $\tau=0$ provides the asymptotic distribution of the highest point of $K_{\lambda}$ above the edge $(a_i,a_{i+1})$. Straightforward computation show that the asymptotic distribution of $R_0$ admits the simple density function
\begin{align*}
\rho\longmapsto \frac1{\|a_{i+1}-a_i\|}\indicat_{\{\rho\in(0,\|a_{i+1}-a_i\|)\}}
\end{align*}
that is the highest point is asymptotically uniformly distributed along the edge $(a_i,a_{i+1})$.

\subsection{Proof of Theorem \ref{theopolygon} (iii): intensity and number of vertices}\label{subsec:verticespolygon}

\pass

\noi Let us rewrite \eqref{eq:efronvertices} when $K$ is a polygon. We proceed in the same way as for the the proof of point (ii), i.e. we decompose the integral in \eqref{eq:efronvertices} into $2n_K$ integrals over intervals $(\tfrac1{2}(\delta_{i-1}+\delta_i),\delta_i)$ and $(\delta_i,\tfrac1{2}(\delta_i+\delta_{i+1}))$.
Similar changes of variables yield
\begin{align*}
\E(\cN(K_{\lambda}))
& \underset{\lambda\to\infty}{\sim} 8\lambda\sum_{i=1}^{n_K} \int_{0}^{\infty} p_o(K,\delta_i-\lambda^{-\gamma})\E(p_o(K_{\lambda},\delta_i-\lambda^{-\gamma})-p_o(K,\delta_i-\lambda^{-\gamma}))(\lambda^{-\gamma}\log\lambda)\dd\gamma\\
& \underset{\lambda\to\infty}{\sim} 8\lambda\sum_{i=1}^{n_K} (I_i^{(1)}(\lambda)+I_i^{(2)}(\lambda))
\end{align*}

\noi where
\begin{align*}
I_i^{(1)}(\lambda) & = \int_{0}^{\frac12} p_o(K,\delta_i-\lambda^{-\gamma})\E(p_o(K_{\lambda},\delta_i-\lambda^{-\gamma})-p_o(K,\delta_i-\lambda^{-\gamma}))(\lambda^{-\gamma}\log\lambda)\dd\gamma \\
& \hspace{-2cm}\text{and} \\
I_i^{(2)}(\lambda) & =\int_{\frac12}^{\infty} p_o(K,\delta_i-\lambda^{-\gamma})\E(p_o(K_{\lambda},\delta_i-\lambda^{-\gamma})-p_o(K,\delta_i-\lambda^{-\gamma}))(\lambda^{-\gamma}\log\lambda)\dd\gamma.
\end{align*}

\pass Because of point (ii) of Proposition \ref{prop:supportpoly}, the integrand of $I_i^{(1)}(\lambda)$ satisfies
\begin{equation*}
\lambda^{1-\gamma}p_o(K,\tfrac{\pi}{2}-\lambda^{-\gamma})\E(p_o(K_{\lambda},\tfrac{\pi}{2}-\lambda^{-\gamma})-p_o(K,\tfrac{\pi}{2}-\lambda^{-\gamma}))
(\lambda^{-\gamma}\log\lambda)\underset{\lambda\to\infty}{\sim}\|o_i\|\times\frac{1}{6\|o_i\|}=\frac1{6}
\end{equation*}
and we may apply Lebesgue's dominated convergence theorem. Now point (ii) of Proposition \ref{prop:supportpoly2} shows that $I_i^{(2)}(\lambda)$ is negligible, which in turn implies Theorem \ref{theopolygon} (iii). $\square$

\pass We now aim at being more specific on the localization of the vertices of $K_\lambda$. The following statement shows a striking self-similarity of the limiting intensity of the point process of vertices around a fixed vertex of $K$.

\begin{prop}\label{prop:intensitypolygon}
Consider the point process $(\rho_v,\alpha_v)_{v\in V_i}$ of the vertices of $K_{\lambda}$ belonging to $\cS_i$. Then, for all $\gamma\in(0,\tfrac1{2})$, the rescaled point process $(\lambda^{1-2\gamma}\rho_v,\lambda^{\gamma}\alpha_v)_{v\in V_i}$ has an asymptotic intensity given by
\begin{equation}\label{eq:intensitypoly}
\sigma_i(\rho,\alpha)=\frac{8}{3}\|o_i\|^2\rho\alpha^3\exp(-2\|o_i\|\rho\alpha^2)\indicat_{\{\rho>0\}}\indicat_{\{\alpha>0\}}.
\end{equation}
\end{prop}

\pass\textbf{Proof of Proposition \ref{prop:intensitypolygon}.}  Let $R\times A \subset \cS_i$ be fixed and denote by $\cN_i(R\times A)$ the number of points of the process $(\lambda^{1-2\gamma}\rho_v,\lambda^{\gamma}\alpha_v)_{v\in V_i}$ belonging to the set $R\times A$. We have to show that
\begin{align*}
\E(\cN_i(R\times A)) \underset{\lambda\to\infty}{\longrightarrow} \int_{R\times A} \sigma_i(\rho,\alpha)\dd\rho \dd\alpha.
\end{align*}

\begin{figure}[h!]
\includegraphics[scale=0.8]{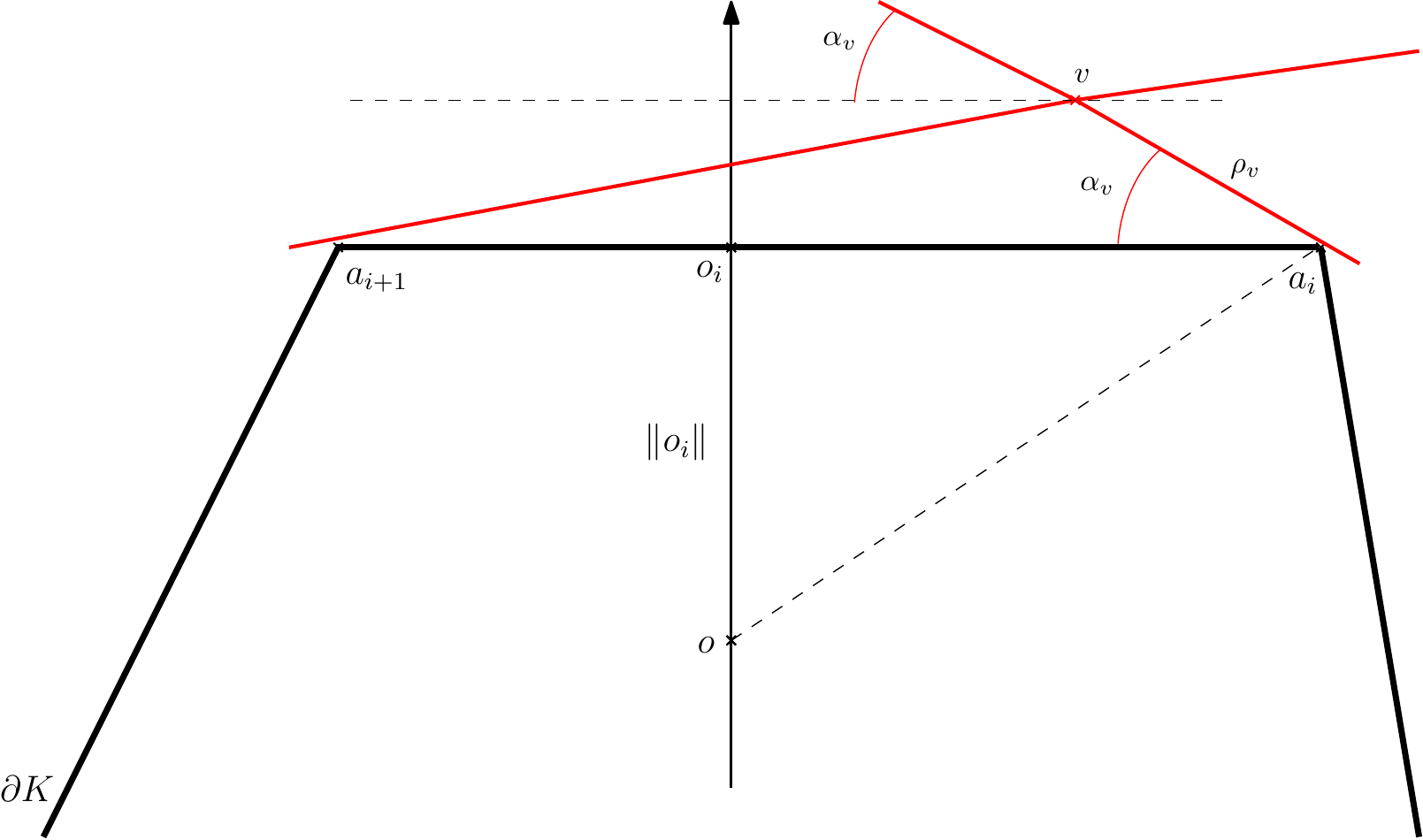}
\caption{\label{fig:sommetpolygon} The analogue of Figure \ref{fig:sommetsmooth} in the polygonal case.}
\end{figure}

\noi The strategy of the proof consists in going along the same lines as for the smooth case and proceeding exactly like for the proof of Proposition \ref{prop:supportpoly}. Precisely, we obtain that
\begin{align*}
\E(\cN_i(R\times A)) = 8\int_{R\times A} \exp\left(-\Delta_i(\lambda,\rho,\alpha)\right)r_i^2(\lambda,\rho,\alpha)J^{\mbox{\scalebox{.5}{supp}}}_i(\lambda,\rho,\alpha)\rho\dd \rho \dd\alpha
\end{align*}
where
\begin{equation*}
\Delta_i(\lambda,\rho,\alpha)= 4\lambda\cA\big(B_{\|s_{a_i,\lambda^{2\gamma-1}\rho,\lambda^{-\gamma}\alpha}\|}(s_{a_i,\lambda^{2\gamma-1}\rho,\lambda^{-\gamma}\alpha})\setminus \cF_o(K)\big)
  = 4\lambda\cA\big(\Delta \cF_{a_i,\lambda^{2\gamma-1}\rho,\lambda^{-\gamma}\alpha}\big)
\end{equation*}
and
\begin{equation*}
J^{\mbox{\scalebox{.5}{supp}}}_i(\lambda,\rho,\alpha)=\lambda^{3\gamma}\int_{E_i(\lambda,\rho,\alpha)} J(\theta_i(\lambda,\rho,\alpha),\theta_1,\theta_2)\dd \theta_1\dd \theta_2
\end{equation*}
where $E_i(\lambda,\rho,\alpha)$ stands for the set of couples $(\theta_1,\theta_2)$ which satisfy that the two bisecting lines of $[o,x_1]$ and $[o,x_2]$ do not intersect $K$ (see Figure \ref{fig:sommetpolygon}).

\pass Thanks to Lemma \ref{lem:geompoly1}, we get
\begin{align*}
\Delta_i(\lambda,\rho,\alpha) \underset{\lambda\to\infty}{\longrightarrow} 2\|o_i\|\rho\alpha^2 .
\end{align*}
Moreover,
\begin{align*}
r_i(\lambda,\rho,\alpha) \underset{\lambda\to\infty}{\longrightarrow} \|a_i\|\,\text{ and }\,\sin(\theta_i(\lambda,\rho,\alpha)) \underset{\lambda\to\infty}{\longrightarrow} \frac{\|o_i\|}{\|a_i\|}.
\end{align*}

\pass Let us now turn on the term $J^{\mbox{\scalebox{.5}{supp}}}_i(\lambda,\rho,\alpha)$. We get successively
\begin{align*}
J^{\mbox{\scalebox{.5}{supp}}}_i(\lambda,\rho,\alpha)
 & \underset{\lambda\to\infty}{\sim} \lambda^{3\gamma}\int_{(0,\lambda^{-\gamma}\alpha)^2} |\sin(\theta_i(\lambda,\rho,\alpha)-\theta_1)\sin(\theta_i(\lambda,\rho,\alpha)-\theta_2)\sin(\theta_2-\theta_1)|\dd \theta_1\dd \theta_2 \\
  & \underset{\lambda\to\infty}{\sim} \lambda^{3\gamma}\Big(\frac{\|o_i\|}{\|a_i\|}\Big)^2 \Big(\frac1{3}(\lambda^{-\gamma}\alpha)^3\Big)\\
  & \hspace{0.22cm} = \frac1{3}\Big(\frac{\|o_i\|}{\|a_i\|}\Big)^2\alpha^3.
\end{align*}

\pass We apply now again Lebesgue's dominated convergence theorem, omitting the domination step which is very similar to what we did in the proof of Proposition \ref{prop:supportpoly}. It follows that
\begin{align*}
\E(\cN_i(R\times A)) \underset{\lambda\to\infty}{\longrightarrow} 8\int_{R\times A} \exp(-2\|o_i\|\rho\alpha^2)\|a_i\|^2\times\frac1{3}\Big(\frac{\|o_i\|}{\|a_i\|}\Big)^2\alpha^3 \rho\dd\rho\dd\alpha
\end{align*}
which implies the required result. $\square$

\section{Answer to Question $3$: the role of the Steiner point}\label{sec:steiner}

\noi In the two previous sections, the cell that we considered is associated with a deterministic nucleus at the origin which is added to the Poisson point process. In particular, the asymptotic shape of the cell depends on both the choice of the convex body $K$ and the position of the origin $o$. In this section, we investigate a modified question which is intrinsic in $K$, i.e. we ask for the behavior of the cell $\widehat{K}_{\lambda}$ containing $K$ when the Poisson point process is conditioned on its associated Voronoi tessellation to not intersect $K$. Since the problem is invariant under translation, we are allowed to assume that the Steiner point of $K$ coincides with the origin, without it being a nucleus of the tessellation. More precisely, the Steiner point of $K$ denoted by  $\mbox{st}(K)$ is defined by the equality
\begin{equation*}
  \mbox{st}(K)=\frac{1}{\pi}\int_0^{2\pi} p_o(K,\theta)u_\theta\dd\theta.
\end{equation*}

\noi When $K$ is smooth, $\mbox{st}(K)$ can be rewritten as
\begin{equation}
\mbox{st}(K)=\frac{1}{\pi}\int_{\partial K}r_s^{-1}\langle s,n_s\rangle n_s\dd s.
\end{equation}
In particular, $\mbox{st}(K)$ is included in the relative interior of $K$, see e.g. Section 1.7 in \cite{Schneider} for the definition and the general properties of $\mbox{st}(K)$. Note in particular that its definition is intrinsic to $K$, i.e. is independent of the choice of the origin $o$. We show as a byproduct of the proof of Proposition \ref{steinerpoint} below that alternatively, $\mbox{st}(K)$ is the unique point $x$ which minimises the function $x\longmapsto \cA(\cF_x(K))$.

\pass Let ${\mathscr S}_{\lambda}$ be the event such that there is one cell $\widehat{K}_\lambda$ of the Poisson-Voronoi tessellation associated with $\cP_{\lambda}$ which contains $K$. We are interested in showing that conditional on ${\mathscr S}_{\lambda}$, the nucleus of $\widehat{K}_\lambda$, denoted by $Z_{\lambda}$, is close to the Steiner point $\mbox{st}(K)$.

\begin{prop}\label{steinerpoint}
Conditional on ${\mathscr S}_{\lambda}$, the rescaled nucleus $\lambda^{\frac{1}{2}}Z_{\lambda}$ converges in distribution as $\lambda\to\infty$ to the centered Gaussian distribution with covariance matrix $(4\pi)^{-1}$ times the identity matrix.
\end{prop}

\pass\textbf{Proof of Proposition \ref{steinerpoint}.}
We start by calculating both $\P({\mathscr S}_{\lambda})$ and the density of $Z_{\lambda}$ conditional on ${\mathscr S}_{\lambda}$. In the sequel, we denote by ${\mathcal C}_x$ the Voronoi cell associated with $x\in\R^2$. 
For any bounded and measurable function $\varphi:\R^2\longrightarrow \R$, we deduce from Mecke-Slivnyak's formula that
\begin{align*}
\E\bigg(\sum_{x\in\cP_{\lambda}}\indic_{\{K\subset \cC_x\}}\varphi(x)\bigg)
=\lambda\int_{\R^2} \P(K \subset \cC_x) \varphi(x) \dd x =\lambda\int_{\R^2} \exp(-4\lambda\cA(\cF_x(K))) \varphi(x) \dd x.
\end{align*}

\noi Taking $\varphi=1$ in the last equality above, we obtain that
\begin{equation*}
\P({\mathscr S}_{\lambda})=\lambda\int_{\R^2}\exp(-4\lambda\cA(\cF_x(K)))\dd x
\end{equation*}
and that the conditional density $f_{{\lambda}}$
of $\lambda^{\frac12}Z_{\lambda}$ is proportional to $x\longmapsto\exp(-4\lambda\cA(\cF_{\lambda^{-1/2}x}(K)))$.

\pass We now turn our attention to the calculation of $\cA(\cF_x(K))$.
For any $x\in \R^2$, we denote by $\cE(x)\subset [0,2\pi]$ the set of all directions such that $p_x(K,\theta)> 0$. Denoting by
\begin{equation*}\label{eq:resteairefleur}
\cR(x)=\frac1{2}\int_{[0,2\pi]\setminus\cE(x)}p_x^2(K,\theta)\dd \theta,
\end{equation*}
we rewrite $\cA(\cF_x(K))$ as
\begin{equation*}
\cA(\cF_x(K)) =\frac1{2}\int_0^{2\pi}p_x^2(K,\theta)\dd\theta-\cR(x).
\end{equation*}

\noi Using now \eqref{eq:suppchangeref}, we get
\begin{align}\label{eq:functionflower}
\cA(\cF_x(K))
& =\frac1{2}\int_{\cE(x)} (p_o(K,\theta)-\langle x,u_{\theta}\rangle)^2\dd\theta-\cR(x) \nonumber \\
& =\cA(\cF_o(K))-\Big\langle x, \int_0^{2\pi} p_o(K,\theta)u_{\theta}\dd\theta\Big\rangle +\frac1{2}\int_0^{2\pi}\langle x,u_{\theta}\rangle^2\dd\theta-\cR(x) \nonumber \\
& =\cA(\cF_o(K))+\frac{\pi}{2}\|x\|^2-\cR(x)
\end{align}
where we have used both the fact that $o$ is the Steiner point of $K$ and the equality
\begin{equation*}
\int_0^{2\pi}\langle x,u_{\theta}\rangle^2\dd\theta= \int_0^{2\pi}(\cos\theta)^2\|x\|^2\dd\theta =\pi \|x\|^2.
\end{equation*}

\pass Let us show two basic properties of the rest ${\cR}(x)$.
\begin{enumerate}[label=-,leftmargin=1\parindent]
\item When $x$ is in the interior of $K$, $\cE(x)=[0,2\pi]$ and $\cR(x)=0$.
\item When $x$ is not in the interior of $K$, because of \eqref{eq:suppchangeref}, $0<p_o(K,\theta)\le \langle x,u_\theta\rangle$ as soon as $\theta\in [0,2\pi]\setminus\cE(x)$ and consequently, for any $\theta\in [0,2\pi]\setminus\cE(x)$,
    \begin{equation*}
    p_x(K,\theta)^2=\langle x,u_\theta\rangle^2-p_o(K,\theta)(2\langle x,u_\theta\rangle-p_o(K,\theta))\le \langle x,u_\theta\rangle^2.
    \end{equation*}
\end{enumerate}

\pass Combining this inequality with the fact that $[0,2\pi]\setminus \cE(x)$ is an interval of length at most $\pi$ implies in turn that
\begin{equation*}\label{eq:boundR(x)}
0\le \cR(x)\le \frac12\int_{[0,2\pi]\setminus\cE(x)} \langle x,u_\theta\rangle^2\dd\theta\le \frac{\pi}{4}\|x\|^2.
\end{equation*}
In view of \eqref{eq:functionflower}, this means in particular that $o$ is the unique minimum of the function $x\longmapsto \cA(\cF_x(K))$.

\pass Now, inserting \eqref{eq:functionflower} into the conditional density function $f_{\lambda}$ of $\lambda^{\frac12}Z_\lambda$, we obtain that
$f_{\lambda}(x)$ is proportional to $\exp(-2\pi\|x\|^2+4\lambda{\cR}(\lambda^{-\frac12}x))$. Using the properties of $\cR$ detailed above, we get that for every $x\in\R^2$, $f_\lambda(x)$ converges to $2\exp(-2\pi\|x\|^2)$ and
\begin{equation*}
\exp(-2\pi\|x\|^2+4\lambda{\cR}(\lambda^{-\frac12}x))\le \exp(-\pi\|x\|^2).
\end{equation*}

\noi Consequently, an application of Lebesgue's dominated convergence theorem shows that for any measurable function $g:\R^2\longrightarrow \R$ which is bounded by a polynomial of $\|x\|$,
\begin{equation*}
\int_{\R^2} g(x)f_\lambda(x)\dd x  \underset{\lambda\to\infty}{\longrightarrow} 2\int_{\R^2} g(x)\exp(-2\pi\|x\|^2)\dd x
\end{equation*}
which completes the proof of Proposition \ref{steinerpoint}. $\square$

\pass The next proposition provides a precise description of the conditional distribution of $\cP_{\lambda}$ given ${\mathscr S}_{\lambda}$.

\begin{prop}\label{steinerpoint2}
The conditional distribution of $\cP_\lambda$ given ${\mathscr S}_\lambda$ is equal in distribution to $\{Z_{\lambda}\}\cup \cP^{(Z_\lambda)}_\lambda$ where $Z_\lambda$ is a random variable distributed according to a density function proportional to $x\longmapsto\exp(-4\lambda\cA(\cF_x(K)))$ and, given $\{Z_{\lambda}=x\}$, $\cP^{(Z_\lambda)}_\lambda$ is a Poisson point process of intensity $\lambda{\emph{\indicat}}_{\R^2\setminus 2\cF_x(K)}$.
\end{prop}

\pass\textbf{Proof of Proposition \ref{steinerpoint2}.}
Let $L$ be a fixed compact set in $\R^2$. Using Mecke-Slivnyak's formula, we get successively
\begin{align*}
\E\Big(\indic_{{\mathscr S}_{\lambda}}\indic_{\{\cP_{\lambda}\cap L=\emptyset\}}\Big)
   & = \E\bigg(\sum_{x\in \cP_{\lambda}}\indic_{\{K\subset \cC_x\}}\indic_{\{\cP_{\lambda}\cap L=\emptyset\}}\bigg)\\
   & =\int_{\R^2} \P(K\subset \cC_x,(\cP_{\lambda}\cup\{x\})\cap L=\emptyset)\dd x\\
   & =\int_{\R^2\setminus L}\P(\cP_{\lambda}\cap (L\cup 2\cF_x(K))=\emptyset)\dd x\\
   & =\int_{\R^2\setminus L}\exp(-\lambda\cA(L\setminus 2\cF_x(K)))\exp(-4\lambda\cA(\cF_x(K)))\dd x.
\end{align*}
Dividing the last equality by $\P({\mathscr S}_{\lambda})$, we get the required result. $\square$

\pass We are now ready to get asymptotic expectations for the area, perimeter and number of vertices of $\widehat{K}_\lambda$.

\begin{theo}\label{theosteiner}
Let $K$ be a convex body with its Steiner point at the origin. The asymptotics of $\E(\cA(\widehat{K}_\lambda))-\cA(K)$, $\E(\cU(\widehat{K}_\lambda))-\cU(K)$ and $\E(\cN(\widehat{K}_\lambda))$ are then provided by Theorem \ref{theosmooth} when $K$ is smooth and by Theorem \ref{theopolygon} when $K$ is a polygon.
\end{theo}

\pass\textbf{Proof of Theorem \ref{theosteiner}.}
We prove the result for $\E(\cA(\widehat{K}_\lambda))-\cA(K)$ and explain at the end how to adapt the arguments for $\E((\cU(\widehat{K}_\lambda))-\cU(K)$ and $\E(\cN(\widehat{K}_\lambda))$.

\pass Let $\cC_o$ be the Voronoi cell associated with the origin when $o$ is added to the set of nuclei $\cP_\lambda$. The cell $\widehat{K}_{\lambda}$ containing $K$ is distributed, up to a translation, as $\cC_o$ conditional on $(K+Z_{\lambda})\subset \cC_o$, where $Z_{\lambda}$ is distributed as in Proposition \ref{steinerpoint2}. Recalling that $f_{\lambda}$ is the density function of $\lambda^{\frac12}Z_\lambda$, we obtain
\begin{align*}
\E(\cA(\widehat{K}_\lambda)-\cA(K)) = \int_{\R^2}\E(\cA((K+\lambda^{-\frac{1}{2}}x)_\lambda)-\cA(K))f_{\lambda}(x)\dd x = I_1(\lambda)+I_2(\lambda)
\end{align*}
where
\begin{align*}
I_1(\lambda) & =\int_{\lambda^{\frac14}K}\E(\cA((K+\lambda^{-\frac{1}{2}}x)_\lambda)-\cA(K))f_{\lambda}(x)\dd x \\
&\hspace{-4cm}\text{and} \\
I_2(\lambda) & =\int_{\R^2\setminus\lambda^{\frac14}K}\E(\cA((K+\lambda^{-\frac{1}{2}}x)_\lambda)-\cA(K))f_{\lambda}(x)\dd x.
\end{align*}

\pass We start by showing that
\begin{equation}\label{eq:convI1}
I_1(\lambda)\underset{\lambda\to\infty}{\sim}\E(\cA({K}_\lambda)-\cA(K)).
\end{equation}
Indeed, a method similar to what has been done in Sections \ref{subsec:surfacesmooth} and \ref{subsec:surfacepolygon} shows that, uniformly in $x\in\R^2$,
\begin{equation*}
\E(\cA((K+\lambda^{-\frac1{2}}x)_\lambda)-\cA(K))\indicat_{\{x\,\in\,\lambda^{\frac14}K\}}\underset{\lambda\to\infty}{\sim}\E(\cA({K}_\lambda)-\cA(K)).
\end{equation*}
Combining this with the convergence and domination of the function $f_\lambda$ showed in the proof of Proposition \ref{steinerpoint}, we get \eqref{eq:convI1}.

\pass Let us show now that the integral $I_2(\lambda)$ is negligible in front of $I_1(\lambda)$. To do so, we denote by $R_x$ the maximal distance from $o$ to the farthest point in $(K+\lambda^{-\frac12}x)_\lambda$. We notice in particular that $\cA((K+\lambda^{-\frac1{2}}x)_\lambda)\le \pi R_x^2$. Moreover, by a method similar to Lemma $1$ in \cite{fs96}, we obtain that, for any $r>0$,
\begin{equation*}
\P(R_x\ge r)\le C'\exp(-Cr^2+C'\lambda^{-1}\|x\|^2),
\end{equation*}
for some positive constants $C,C'>0$.

\noi Consequently, $\E(\cA((K+\lambda^{-\frac1{2}}x)_\lambda)-\cA(K))$ is bounded by $1+\lambda^{-1}\|x\|^2$ up to a multiplicative constant. Using the domination of $f_\lambda$ showed in Proposition \ref{steinerpoint}, we get that $I_2(\lambda)\to 0$ exponentially fast as $\lambda\to\infty$. Combining this last result with \eqref{eq:convI1}, we obtain the required convergence for the mean defect area of $\widehat{K}_\lambda$.

\pass Finally, the estimates for $\E(\cU(\widehat{K}_\lambda)-\cU(K))$ and $\E(\cN(\widehat{K}_\lambda))$ follow from similar arguments, as soon as we are able to get bounds for $\E(\cU((K+\lambda^{-\frac1{2}}x)_\lambda)-\cU(K))$ and $\E(\cN((K+\lambda^{-\frac1{2}}x)_\lambda))$. Using the inclusion $(K+\lambda^{-\frac12}x)_\lambda\subset B_o(R_x)$, we get that, up to a multiplicative constant, $\lambda^{-\frac12}\|x\|$ is an upper-bound of $\E(\cU((K+\lambda^{-\frac{1}{2}}x)_\lambda))-\cU(K))$. A use of Proposition \ref{prop:efron} (i) combined with the same inclusion finally shows that, up to a multiplicative constant, $\lambda^{-1}\|x\|^2$ is an upper-bound of $\E(\cN((K+\lambda^{-\frac1{2}}x)_\lambda))$. $\square$

\section{Answer to Question $1$}\label{sec:tele}

\noi Given a bounded closet set $\cD\subset\R^2$ containing $o$ in its interior, we draw a Poisson point process with intensity $\lambda$ outside $\cD$ and look at the limit shape of the cell $\cC_\lambda(\cD)$ containing $o$ associated with the Voronoi tessellation when $\lambda\to\infty$. We would like to apply previous results which describe the shape of such a cell but when the Poisson point process is drawn outside the Voronoi flower of a convex body. Actually, this is not possible directly. However, it will turn out that the change of $\cD$ by the maximal Voronoi flower $\cF_\cD$ with respect to $o$ included in $\cD$  will not affect the first-order asymptotics of the characteristics of the cell.

\pass In this section, all the Voronoi flowers will be considered with respect to $o$.

\subsection{Voronoi flower and antiorthotomic curve of a set}\label{subsec:isotel}

\pass

\noi This section is devoted to the description of the set $\cF_\cD$. Namely, we will prove that this set is the homothetical image with ratio $2$ of the flower of the so-called \emph{antiorthotomic curve} $\Gamma_\cD$ of $\cD$, that is the curve $\Gamma_\cD$ made of points which are equidistant from $o$ and $\partial \cD$. Equivalently, $\partial \cF_\cD$ is  the locus of reflexions of $o$ in the tangent lines to $\Gamma_\cD$, see e.g. Exercise 2 page 132 in \cite{bg84}. We will use in several places that $\Gamma_\cD$ is the boundary of a convex body $K$, namely the set of points which are closer to $o$ than to $\partial \cD$. Indeed, if $x$ and $y$ are in $K$, $\cF_o(\{x,y\})$ does not meet $\partial \cD$ and since $\cF_o(\{x,y\})=\cF_o([x,y])$, $[x,y]$ is also included in $K$.

\pass Let us start before with a general result characterizing sets which are a Voronoi flower. In the sequel, we denote by $i$ the inversion with pole $o$ and with respect to the unit-circle.

\begin{lem}\label{lem:isaflower1}
The set $\cD$ is the Voronoi flower of a convex set of $\R^2$ if and only if the complementary of its image by $i$ is convex.
Moreover, when such a convex set exists it is unique. 
\end{lem}

\pass\textbf{Proof of Lemma \ref{lem:isaflower1}.}
Assume first that $\cD$ is the Voronoi flower of a convex set $K$. Then the set $\cD$ is equal to
\begin{equation*}
\cD=\Big\{t p_o(K,u_\theta)u_{\theta} : (t,\theta)\in[0,1]\times[0,2\pi]\Big\}.
\end{equation*}

\noi Hence, the complementary of its image by $i$ writes
\begin{equation*}
\R^2\setminus i(\cD)= \left\{\frac{tz}{p_o(K,z)} : (t,z)\in [0,1)\times\R^2\setminus\{0\} \right\}.
\end{equation*}

\noi We use now the second proof of Theorem $1.7.1$ of \cite{Schneider} which claims that the set $\{ (z,h) : z\in\R^2 \mbox{ and } h\ge p_o(K,z)\}$ of $\R^3$ is a convex cone. Its section $\{(\tfrac{tz}{p_o(K,z)},1): (t,z)\in[0,1]\times\R^2\}$ is therefore also a convex set. The interior of the projection on the first two coordinates is nothing but $\R^2\setminus i(\cD)$ which is then convex too. Conversely, assume that $\R^2\setminus i(\cD)$ is a convex set. Notice that since it contains $o$ it is starlike with respect to $o$. It is also bounded since $i(\cD) $ contains the image by $i$ of a small disk centered at $o$. Hence, there exists a function $g: \R^2 \to (0,\infty)$ such that $\tfrac1{g}$ is homogeneous of degree $1$ and such that $\R^2\backslash i(\cD)= \{ tg(z)z : (t,z) \in [0,1)\times\R^2\setminus\{0\}\}$.

\noi Consequently, the cone of $\R^3$ with apex $o$ and generated by $(\R^2\setminus i(\cD))\times\{1\}$ is the set $\{(z,h): z\in\R^2 \mbox{ and } h>\tfrac1{g(z)}\}.$ Since the convexity of the epigraph of a homogenous function of degree $1$ is equivalent to the sublinearity of the function, the converse of Theorem 1.7.1 of \cite{Schneider} implies that $\tfrac1{g}$ is the support function of some convex set $K$. Since $\cD=\{ \tfrac{t}{g(u_{\theta})}u_{\theta} : (t,\theta)\in[0,1]\times[0,2\pi] \}$ it follows that $\cD$ is the Voronoi flower of $K$. $\square$

\pass The following result provides a geometric description of $\cF_\cD$ which will lead to an analytical expression of $\cF_\cD$.

\begin{prop}\label{prop:isaflower3}
Let us denote by $\cD^*$ the maximal starlike set with respect to $o$ included in $\cD$ and by $\cF_\cD$ the set defined by
\begin{equation}\label{eq:fleurmax}
\cF_\cD=\R^2\setminus i(\conv(\R^2\setminus i(\cD^*))).
\end{equation}
Then ${\cF_\cD}$ is the the unique maximal Voronoi flower included in $\cD$. Moreover, the set $\tfrac1{2}\partial\cF_\cD$ is the flower of the antiorthotomic curve $\Gamma_\cD$ with respect to $o$.
\end{prop}

\pass\textbf{Proof of Proposition \ref{prop:isaflower3}.}
Let us prove first that the set defined by \eqref{eq:fleurmax} is the unique maximal Voronoi flower included in $\cD$. Relation \eqref{eq:fleurmax} implies straightforwardly that $\cF_\cD\subset\cD^*$. Moreover $\cF_\cD$ is a Voronoi flower by application of Lemma \ref{lem:isaflower1}. Now $\cF_\cD$ is maximal in $\cD^*$. Indeed, for any Voronoi flower $\cF$ included in $\cD^*$, the set $\R^2\setminus i(\cF)$ is a convex set containing $\R^2\setminus i(\cD^*)$, thus containing $\conv(\R^2\setminus i(\cD))$. This inclusion $\conv(\R^2\setminus i(\cD^*))\subset \R^2\setminus i(\cF)$ implies in turn that
\begin{equation*}
\cF_\cD=\R^2\setminus i(\conv(\R^2\setminus i(\cD^*))) \supset \R^2\setminus i(\R^2\setminus i(\cF))=\cF.
\end{equation*}
We conclude by noticing that a maximal Voronoi flower included in $\cD$ is starlike with respect to $o$ and therefore is included in $\cD^*$. Hence $\cF_\cD$ is also maximal in $\cD$.

\pass Let us determine now the Voronoi flower of $\Gamma_\cD$. For all $x\in\Gamma_\cD$, the disk $ B_{\|x\|}(x)$ is included in $\cD$. Hence by maximality of $\cF_\cD$, the flower of $\Gamma_\cD$ is included in $\tfrac1{2}\cF_\cD$. Conversely, $\tfrac1{2}\cF_\cD$ is the flower of some set $K$. If this set $K$ is not included in the interior of $\Gamma_\cD$ then it exists $x\in K$ but not in the interior of the set with boundary $\Gamma_\cD$ such that $B_{\|x\|}(x)$ has a part outside $\cD$. This is not possible since $\cF_\cD\subset\cD$. $\square$

\subsection{The main result}\label{subsec:telemain}

\pass

\begin{figure}[!h]
\begin{center}
\includegraphics[scale=0.75]{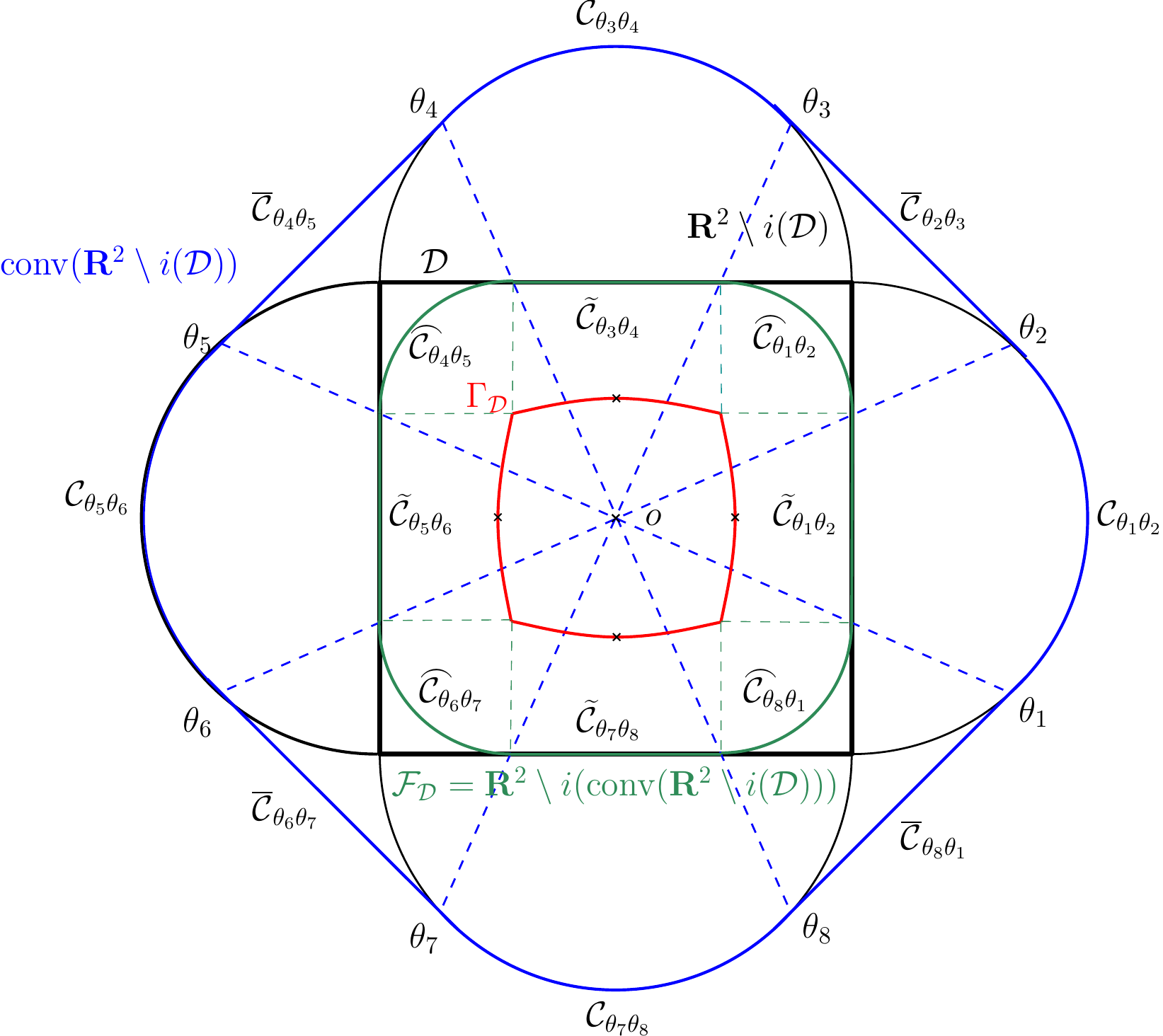}
\caption{\label{fig:television} Construction of the unique maximal flower $\cF_\cD$ included in the bounded closed set $\cD$. Actually $\cF_\cD$ is the Voronoi flower of $2\Gamma_\cD$.}
\end{center}
\end{figure}

\noi Equality \eqref{eq:fleurmax} combined with the fact that the flower of $\Gamma_\cD$ is $\tfrac1{2}\partial \cF_\cD$ implies that $\Gamma_\cD$ only depends on $\cD^*$, i.e. the maximal starlike set with respect to $o$ included in $\cD$. In the next lemma, we fix a polar equation of $\cD^*$ with respect to $o$ and deduce from it a parametric equation of $\Gamma_\cD$ (see Figure \ref{fig:television}).

\begin{lem}\label{lem:isaflower2}
Assume that $\partial \cD^* =\{d(\theta):\theta\in[0,2\pi]\}$ with $d$ a function which is piecewise of class $\cC^3$.
\begin{enumerate}
\item[$(i)$] If $\cD^*$ is a Voronoi flower then
\begin{equation*}\label{eq:polfleur1}
\Gamma_\cD = \dfrac1{2}\Big\{(d(\theta)\cos\theta - d'(\theta)\sin\theta, d(\theta)\sin\theta+d'(\theta)\cos\theta) : \theta\in[0,2\pi] \Big\}
\end{equation*}
and the Voronoi flower of $\Gamma_\cD$ satisfies $2\cF(\Gamma_\cD)=\cD$. \\
\item[$(ii)$] If $\cD^*$ is not a Voronoi flower but is such that it exists angles
$0\leq\theta_1\leq\theta_2\leq\cdots\leq\theta_{2n}\leq 2\pi$ such that $\conv(\R^2\setminus i(\cD))$ is strictly convex in the cone of apex $o$ and directions in $\Theta=[\theta_1,\theta_2]\cup [\theta_3,\theta_4]\cup\cdots \cup [\theta_{2n-1},\theta_{2n}]$ then
\begin{equation*}\label{eq:polfleur2}
\hspace{1.5cm} \Gamma_\cD = \frac1{2}\bigcup_{i=0}^{n-1}\Big\{(d(\theta)\cos\theta - d'(\theta)\sin\theta, d(\theta)\sin\theta+d'(\theta)\cos\theta) : \theta\in[\theta_{2i+1},\theta_{2i+2}]\Big\}
\end{equation*}
and the Voronoi flower of $\Gamma_\cD$ satisfies $2\cF(\Gamma_\cD)\subsetneq \cD^*$ with
\begin{equation*}\label{eq:polfleur3}
\cD^*\backslash 2\cF(\Gamma_\cD)= \bigcup_{i=1}^{n-1} \cD_{\theta_{2i}\theta_{2i+1}}
\end{equation*}
where $\cD_{\theta_{2i}\theta_{2i+1}}$ is the set delimited on one hand by $\partial \cD^*$ and on the other hand by a circular arc $\overset{\frown}{\cC}_{\theta_{2i}\theta_{2i+1}}$ which is a part of a circle containing $o$ and tangent to $\partial \cD^*$.
\end{enumerate}
\end{lem}

\pass\textbf{Proof of Lemma \ref{lem:isaflower2}.}

\noi {\it Proof of $(i)$.} The set $\cD$ is the set of the orthogonal projections of $o$ onto the tangent lines to $2\Gamma_\cD$. Conversely, since $\cD$ is a Voronoi flower, the set $2\Gamma_\cD$ is the envelope of all lines containing $s$ and orthogonal to $s$ for $s\in\partial\cD$. This envelope is usually called the reciprocal pedal curve of $\partial \cD$. Using Exercise 2, page 132 in \cite{bg84}, straightforward computations show that $2\Gamma_\cD$ may be parametrized as
\begin{equation*}
\left\{
  \begin{array}{ll}
    x(\theta)= d(\theta)\cos\theta - d'(\theta)\sin\theta\\
    y(\theta)= d(\theta)\sin\theta + d'(\theta)\cos\theta
   \end{array}
\right..
\end{equation*}

\pass{\it Proof of $(ii)$.}
Keeping the notations of Proposition \ref{prop:isaflower3} and denoting by $\cF(\Gamma_\cD)$ the Voronoi flower of $\Gamma_\cD$ we have
\begin{equation*}
2\cF(\Gamma_\cD)=\cF_\cD = \R^2\setminus i(\conv(\R^2\setminus i(\cD))).
\end{equation*}

\noi It follows from hypotheses that $\partial(\conv(\R^2\setminus i(\cD)))$ is a union of a finite number $n$ of strictly convex and regular curves $\cC_{\theta_{2i+1}\theta_{2i+1}}$ and a finite number $n-1$ of straight lines $\overline{\cC}_{\theta_{2i}\theta_{2i+1}}$ alternated
\begin{equation*}
\partial(\conv(\R^2\setminus i(\cD))) =  \cC_{\theta_1\theta_2} \cup  \overline{\cC}_{\theta_2\theta_3} \cup \cC_{\theta_3\theta_4}\cup \cdots \cup \overline{\cC}_{\theta_{2n-2}\theta_{2n-1}}\cup\cC_{\theta_{2n-1}\theta_{2n}}.
\end{equation*}

\noi Therefore, by properties of inversion, \eqref{eq:fleurmax} implies that $\partial\cF_\cD$ is the corresponding union of a finite number $n$ of regular curves $\widetilde{\cC}_{\theta_{2i+1}\theta_{2i+2}}$ and a finite number $n-1$ of circular arcs $\overset{\frown}{\cC}_{\theta_{2i}\theta_{2i+1}}$ consecutively
\begin{equation*}
  \partial\cF_\cD =  \widetilde{\cC}_{\theta_1\theta_2} \cup  \overset{\frown}{\cC}_{\theta_2\theta_3} \cup \widetilde{\cC}_{\theta_3\theta_4}\cup \cdots \cup \overset{\frown}{\cC}_{\theta_{2n-2}\theta_{2n-1}}\cup\widetilde{\cC}_{\theta_{2n-1}\theta_{2n}}.
\end{equation*}

\noi Let us denote by $\{\widetilde{d}(\theta):\theta\in[0,2\pi]\}$ the polar representation of $\cF_\cD$. Noticing that \eqref{eq:fleurmax} implies that the convex parts of $\partial\cD$ are unchanged under the transformation giving $\cF_\cD$, we obtain more precisely
\begin{equation*}
\partial\cD\cap \partial \cF_\cD
= \bigcup_{i=0}^{n-1} \widetilde{\cC}_{\theta_{2i+1}\theta_{2i+2}}
= \bigcup_{i=0}^{n-1}\Big\{(d(\theta)\cos\theta, d(\theta)\sin\theta) : \theta\in[\theta_{2i+1},\theta_{2i+2}]\Big\}
\end{equation*}
and
\begin{equation*}
\partial\cF_\cD\backslash \partial(\cD\cap\cF_\cD) = \bigcup_{i=1}^{n-1} \overset{\frown}{\cC}_{\theta_{2i}\theta_{2i+1}}
= \bigcup_{i=1}^{n-1}\Big\{(\widetilde{d}(\theta)\cos\theta, \widetilde{d}(\theta)\sin\theta) : \theta\in[\theta_{2i},\theta_{2i+1}]\Big\}.
\end{equation*}

\noi Noticing moreover that each arc $\overset{\frown}{\cC}_{\theta_{2i}\theta_{2i+1}}$ is a part of a circle containing $o$ and joining tangentially the two neighboring curves $\widetilde{\cC}_{\theta_{2i-1}\theta_{2i}}$ and $\widetilde{\cC}_{\theta_{2i+1}\theta_{2i+2}}$, this proves the second equality in (ii).

\pass Furthermore, this property of arcs $\overset{\frown}{\cC}_{\theta_{2i}\theta_{2i+1}}$ yields
\begin{equation*}
\widetilde{d}(\theta)=\left\{
  \begin{array}{ll}
   d(\theta) & \hbox{if $\theta\in\Theta$} \\
   d(\theta_{2i}) & \hbox{if $\theta\in[\theta_{2i},\theta_{2i+1}]$}
  \end{array}
\right..
\end{equation*}
We then obtain the first equality in (ii) by applying (i) with $\widetilde{d}$. $\square$

\pass

\pass We are now able to state the main result of this section which answers to Question 1. We recall that  $\cC_\lambda(\cD)$ is the Voronoi cell containing the origin $o$ when $\cP_\lambda$ is conditioned on not intersecting $\cD$.

\begin{theo}\label{theo:tele}
Assume that $\partial \cD^*=\{d(\theta):\theta\in[0,2\pi]\}$ with $d$ a function which is piecewise of class $\cC^3$. Assume moreover the existence of angles
$0\leq\theta_1\leq\theta_2\leq\cdots\leq\theta_{2n}\leq 2\pi$ such that $\conv(\R^2\setminus i(\cD^*))$ is strictly convex in the cone of apex $o$ and directions in $\Theta=[\theta_1,\theta_2]\cup[\theta_3,\theta_4
]\cup \cdots \cup [\theta_{2n-1},\theta_{2n}]$. Under these conditions,

\vspace{0.15cm}
\begin{enumerate}
\item[$(i)$] The cell $\cC_{\lambda}(\cD)$ satisfies, when the intensity $\lambda\to\infty$,
\begin{equation*}
\cC_{\lambda}(\cD) \underset{\lambda\to\infty}{\longrightarrow} K
\end{equation*}
almost surely in the Hausdorff metric, where $K$ is the compact and convex set containing $o$ and whose boundary is the antiorthotomic curve $\Gamma_\cD$ of $\cD$ given by
\begin{equation*}\label{eq:isotele}
\hspace{1.5cm} \Gamma_\cD = \frac1{2}\bigcup_{i=0}^{n-1}\Big\{(d(\theta)\cos\theta - d'(\theta)\sin\theta, d(\theta)\sin\theta+d'(\theta)\cos\theta) : \theta\in[\theta_{2i+1},\theta_{2i+2}]\Big\}.
\end{equation*}
\item[$(ii)$] The defect area, defect perimeter and number of vertices of $\cC_{\lambda}(\cD)$ have respectively the following asymptotics when the intensity $\lambda\to\infty:$
\begin{align*}
& \hspace{.5cm} \E(\cA(\cC_{\lambda}(\cD)))-\cA(K) \underset{\lambda\to\infty}{\sim} \lambda^{-\frac{2}{3}}2^{-\frac{8}{3}}3^{-\frac1{3}}\Gamma\Big(\mfrac{2}{3}\Big)
\int_\Theta (d(\theta)+d''(\theta))^{\frac{4}{3}}d(\theta)^{-\frac{2}{3}}\,\dd \theta \\
& \hspace{.5cm} \E(\cU(\cC_{\lambda}(\cD)))-\cU(K) \underset{\lambda\to\infty}{\sim} \lambda^{-\frac{2}{3}}2^{-\frac{2}{3}}3^{-\frac{4}{3}}\Gamma\Big(\mfrac{2}{3}\Big)
\int_\Theta (d(\theta)+d''(\theta))^{\frac1{3}}d(\theta)^{-\frac{2}{3}}\,\dd \theta \\
& \hspace{.5cm} \E(\cN(\cC_{\lambda}(\cD))) \underset{\lambda\to\infty}{\sim} \lambda^{\frac1{3}}2^{-\frac{8}{3}} 3^{-\frac{4}{3}}\Gamma\Big(\mfrac{2}{3}\Big)\int_\Theta (d(\theta)+d''(\theta))^{\frac1{3}}d(\theta)^{\frac1{3}}\,\dd \theta.
\end{align*}
\end{enumerate}
\end{theo}

\pass\textbf{Proof of Theorem \ref{theo:tele}.}

\noi{\it Proof of $(i)$.}
This assertion is a consequence of the area asymptotics stated in $(ii)$.
Indeed, fixing $h>0$, the inequality $\dd_H(\cC_{\lambda}(\cD), K)>h$ combined with the fact that both $K$ and $\cC_{\lambda}(\cD)$ are convex bodies, implies that
\begin{equation*}
\cA(\cC_{\lambda}(\cD)\setminus K)\ge \inf_{s\in \partial K} \cA(\conv(K\cup \{s+hn_s\})\setminus K)=\eta_h >0.
\end{equation*}
Consequently,
\begin{equation*}
\P(\dd_H(\cC_{\lambda}(\cD),K)>h)\le \P(\cA(\cC_{\lambda}(\cD)\setminus K)\ge\eta_h)\le \eta_h^{-1}\E(\cA(\cC_{\lambda}(\cD)\setminus K))\underset{\lambda\to\infty}{\longrightarrow} 0
\end{equation*}
because of the area asymptotics in $(ii)$. This shows that the convergence of $\cC_\lambda(\cD)$ to $K$ occurs in probability. Since $\cC_\lambda(\cD)$ is decreasing for the inclusion, its almost sure limit can only be $K$. Finally, the representation of $\Gamma_\cD$ follows from Lemma \ref{lem:isaflower2}.

\pass{\it Proof of $(ii)$.}
The strategy is to apply the results of Section \ref{sec:smooth} and Section \ref{sec:polygon} with the set $K$ which is, up to a finite number of isolated angular points in $\partial K$, analogous to the smooth case. Indeed, since $d$ is piecewise $\cC^3$, the boundary of $K$, i.e. the curve $\Gamma_{\cD}$, is piecewise $\cC^2$ and with bounded and positive curvature. The influence of the angular points will be treated by using the estimates proved for the polygonal case in Theorem \ref{theopolygon}.

\pass When $\cD$ is a Voronoi flower, it follows from point (ii) of Lemma \ref{lem:isaflower2} that $\partial K$ is a smooth curve and we may apply Theorem \ref{theosmooth}.

\begin{figure}[!h]
\begin{center}
\includegraphics[scale=0.7]{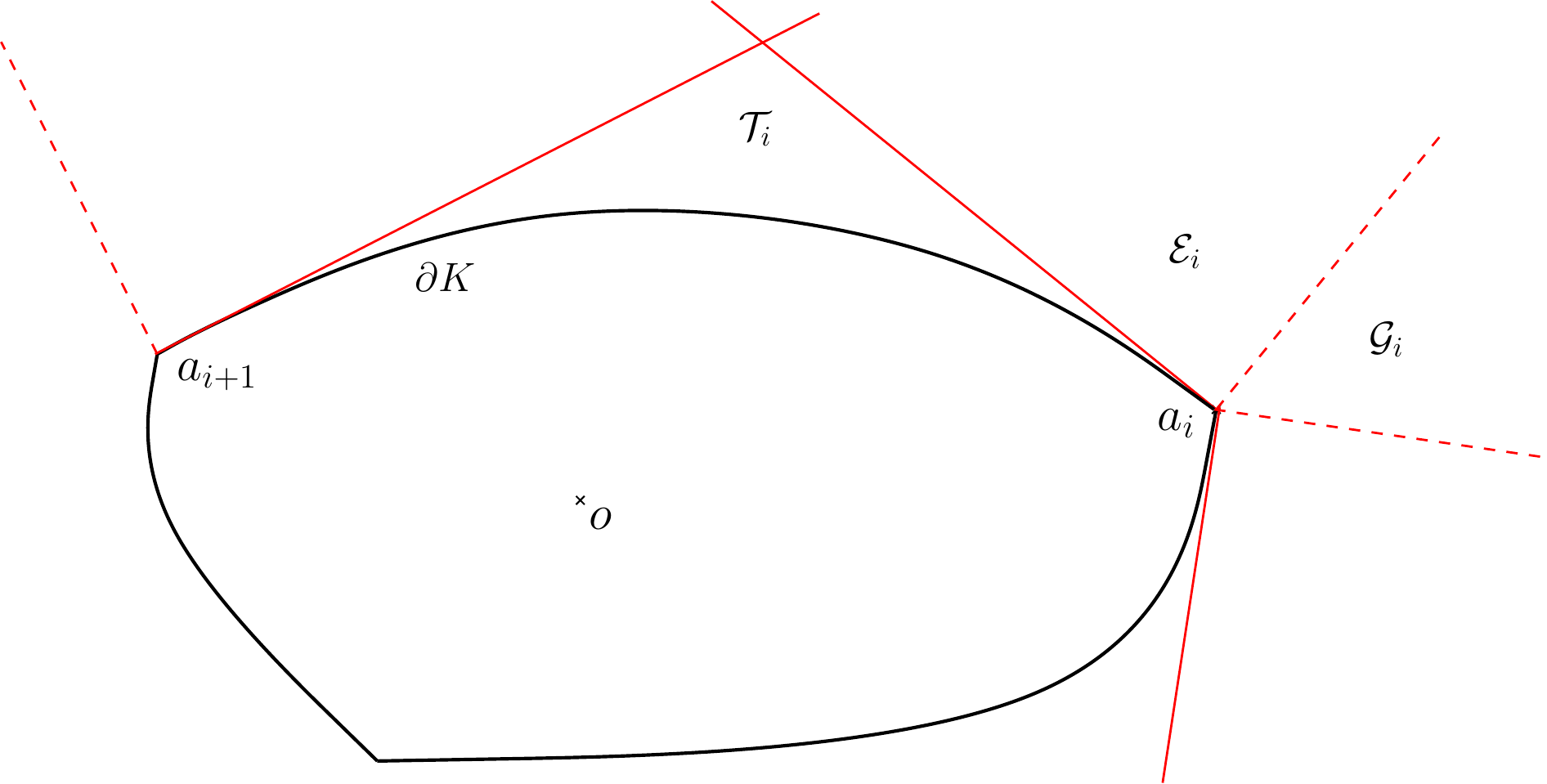}
\vspace{0.2cm}
\caption{\label{fig:mixedcase} The decomposition of space outside of $K$ near an angular point}
\end{center}
\end{figure}

\pass When $\cD$ is a not a Voronoi flower, but is still a starlike set, the error with the previous case is due to the influence of points of the Poisson point process inside the set $\cD\setminus \cF_\cD$. This set is exactly described by the second point of Lemma \ref{lem:isaflower2} (ii). Let us then consider another Poisson point process $\widetilde\cP_{\lambda}$ drawn into $\cD\cup \cup_{i=1}^{n-1} \cD_{\theta_{2i}\theta_{2i+1}}$ independently of the original process $\cP_{\lambda}\cap(\R^2\setminus\cD)$. Let us denote by $\widetilde\cC_{\lambda}(\cD)$ the corresponding cell of the associated Voronoi tessellation containing $o$. In particular, $\widetilde\cC_{\lambda}(\cD)\subset\cC_{\lambda}(\cD)$. We investigate the difference $\cC_{\lambda}(\cD)\setminus \widetilde\cC_{\lambda}(\cD)$ by dividing it into three parts. Precisely, at each angular point $a_i$, we denote by $\cT_i$, $\cG_i$ and $\cE_i$ respectively the region delimited by $\Gamma_\cD$, the semi-tangent half-line at $a_i$ and the semi-tangent half-line at $a_{i+1}$, the region delimited by the two semi-normal half-lines at $a_i$ and the region delimited by the semi-tangent half-line at $a_i$ and the corresponding semi-normal half-line at $a_i$ (see Figure \ref{fig:mixedcase}). We describe below the contribution of each of these regions to the asymptotics of $\E(\cA(\cC_{\lambda}(\cD))-\cA(\widetilde\cC_{\lambda}(\cD)))$, $\E(\cU(\cC_{\lambda}(\cD))-\cU(\widetilde\cC_{\lambda}(\cD)))$ and $\E(\cN(\cC_{\lambda}(\cD))-\cN(\widetilde\cC_{\lambda}(\cD)))$.

\pass $-$ \emph{Contribution of $\cT_i$.} We notice that $ \widetilde\cC_{\lambda}(\cD)\cap \cT_i =\cC_{\lambda}(\cD)\cap \cT_i$ so there is no contribution from $\cT_i$.

\pass $-$ \emph{Contribution of $\cG_i$.} The error coming from $\cG_i$ is of the same order as in the polygonal case, namely $\mathrm{O}(\lambda^{-1})$ for the area by \eqref{eq:rayoflight}, $\mathrm{O}(\lambda^{-1}\log\lambda)$ for the perimeter and $\mathrm{O}(\log\lambda)$ for the number of vertices. All three of them are respectively negligible in front of $\E(\cA(\widetilde\cC_{\lambda}(\cD)))$, $\E(\cU(\widetilde\cC_{\lambda}(\cD)))$ and $\E(\cN(\widetilde\cC_{\lambda}(\cD)))$ provided by Theorem \ref{theosmooth}.

\pass $-$ \emph{Contribution of $\cE_i$.} The height of $\cC_{\lambda}(\cD)$ above the smooth part of $\partial K$ between $a_i$ and $a_{i+1}$ is of order $\mathrm{O}(\lambda^{-\frac23})$. This strip intersects $\cE_i$ at distance of order $\lambda^{-\frac13}$ from $a_i$. Consequently, the contribution to the area is of order $\mathrm{O}(\lambda^{-\frac23}\lambda^{-\frac13})=\mathrm{O}(\lambda^{-1})$, which is negligible in front of $\E(\cA(\widetilde\cC_{\lambda}(\cD)))$ by Theorem \ref{theosmooth} applied to $K$. Regarding perimeter and number of vertices, the contributions of $\cE_i$ are deduced from Theorem \ref{theopolygon}, i.e. are of order $\mathrm{O}(\lambda^{-1}\log\lambda)$ and $\mathrm{O}(\log\lambda)$ respectively. They are consequently negligible in front of $\E(\cU(\widetilde\cC_{\lambda}(\cD)))$ and $\E(\cN(\widetilde\cC_{\lambda}(\cD)))$ respectively.

\noi Finally, when $\cD$ is a not a Voronoi flower and not starlike, we compare $\cC_\lambda(\cD)$ to $\cC_\lambda(\cD^*)$. The error between the two comes from the points located in $\cD\setminus \cD^*$. We can show in a very similar way that this error does not affect the expectation asymptotics.

\pass To conclude the proof, it remains to rewrite the asymptotics given by Theorem \ref{theosmooth} using the change of variables $s=s(\theta)$. Straightforward computations give $r_s=r(\theta)=d(\theta)+d''(\theta)$, $\langle s,n_s\rangle=d(\theta)$, $\dd s = r(\theta)\dd\theta$ and then the desired results. $\square$

\section{Further results and questions}\label{sec:conclusion}

\noi In this section, we discuss extensions of our work and describe open questions.

\subsection{Higher dimension}

\pass

\noi The statements of Proposition \ref{prop:efron} can be extended to dimension $d\ge 3$. Indeed, we can show in the exact same way that
\begin{align*}
\E(\cN_{d-1}(K_\lambda)) & = 2^d\lambda(\E(\cA(\cF_o(K_{\lambda}))-\cA(\cF_o(K))) \\
& \hspace{-0.22cm} \underset{\lambda\to\infty}{\sim} 2^d\lambda\int_{{\mathbb S}^{d-1}}p_o^{d-1}(K,u)\E(p_o(K_\lambda,u)-p_o(K,u))\dd \sigma_d(u)
\end{align*}
where $\cN_{d-1}(K_\lambda)$ denotes the number of hyperfaces of $K_\lambda$, $\sigma_d$ is the uniform measure on the unit-sphere ${\mathbb S}^{d-1}$ and,
for every $u\in{\mathbb S}^{d-1}$, $p_o(K,u)$ is the support function of $K$ in direction $u$.

\pass Consequently, the asymptotics of the mean defect volume of $K_\lambda$ and its number of hyperfaces are to be deduced from the asymptotics
of, respectively, the increase of the volume of the Voronoi flower when a point is added near the boundary of $K$ and the mean defect support
function of $K_\lambda$. The problem now is that we have calculated these two quantities in dimension $2$ with a local approximation of
$\partial K$ by an osculating circle in the smooth case and by a broken line in the polygonal case. The extension of these methods to higher
dimension is unclear. Indeed, a smooth convex body is locally approximated by an ellipsoid, which makes the calculation more delicate.
Besides, we cannot use an affine transformation as for the convex hull model since we need to preserve the Euclidean structure.
A similar problem occurs when the R\'enyi-Sulanke approach of Section \ref{sec:rsapproach} is put to the test. In the polyhedral case,
the situation seems even more tricky.

\pass We expect Proposition \ref{steinerpoint} to have an analogue for $d\ge 3$. Only, the Steiner point should be replaced by the point
\begin{equation*}
\mbox{st}_d(K)=\frac1{d\kappa_d}\int_{{\mathbb S}^{d-1}}p_o^{d-1}(K,u)u\dd \sigma_d(u)
\end{equation*}
where $\kappa_d$ is the volume of the $d$-dimensional unit-ball.

\pass Finally, Lemma \ref{lem:isaflower1} and Proposition \ref{prop:isaflower3} can be extended without difficulty. The limit shape of $\cC_\lambda(\cD)$
should be an analogue of the antiorthotomic curve, i.e. the set of points which are equidistant from $o$ and $\partial \cD$. Only, there is little hope
to get an explicit decomposition and explicit spherical equations in the spirit of Lemma \ref{lem:isaflower2} (ii).

\subsection{The Crofton cell}

\pass

\noi The three main questions of the paper prove to be equally appealing when the Poisson-Voronoi tessellation is replaced by any random line tessellation in the plane and in particular by the stationary and isotropic Poisson line tessellation, see e.g. \cite{cal10}. This tessellation is obtained by taking a Poisson point process $\cQ_\lambda$ of intensity measure $\lambda \|x\|^{-1}\dd x$ in $\R^2$ and constructing for every $x$ in the point process, the line $L_{x}$ containing $x$ and with normal vector $x$. The cell containing the origin is the so-called {\it Crofton cell} and is defined as the intersection of all closed half-planes containing the origin and delimited by lines $L_{x}$. We denote by $\overline{K}_\lambda$ a cell distributed as the Crofton cell conditional on the event that no line crosses $K$, which is equivalent to say that no point from $\cQ_\lambda$  meets $\cF_o(K)$. We recall that this event
has probability
\begin{equation}
\exp\bigg(-\lambda\int_{x\in \cF_o(K)}\|x\|^{-1}\dd x\bigg)=\exp\bigg(-\lambda\int_{0}^{2\pi}p_o(K,\theta)\dd \theta\bigg)=\exp(-\lambda\cU(K)).
\end{equation}
This new random polygon $\overline{K}_{\lambda}$ satisfies \eqref{eq:croftonperimeter} whereas \eqref{eq:airesmooth} is replaced by
\begin{align*}
  \E(\cA(\overline{K}_\lambda)-\cA(K)) & =\int_{\R^2\setminus K} \exp\bigg(-\lambda (\cU(\conv(K\cup\{x\}))-\cU(K))\bigg)\dd x.
\end{align*}

\noi Moreover, as in the Voronoi case, we can establish an exact Efron-type identity which connects in a very simple way the mean defect perimeter of $\overline{K}_\lambda$ to its mean number of vertices. Indeed, denoting by $\cC_{x}$ the Crofton cell associated with the set $(\cQ_\lambda\setminus \cF_o(K))\setminus\{x\}$, we get
\begin{equation*}
\E(\cN(\overline{K}_\lambda))=\E\bigg(\sum_{x\in \cQ_{\lambda}\setminus \cF_o(K)}{\emph{\indicat}}_{\{x\in \cF_o(\cC_x)\}}\bigg) = \lambda\int_{\R^2\setminus \cF_o(K)} \P(x\in (\cF_o(\overline{K}_\lambda)\setminus\cF_o(K)))\|x\|^{-1}\dd x.
\end{equation*}
Thanks to Fubini theorem, this yields
\begin{equation}\label{eq:efroncrofton}
  \E(\cN(\overline{K}_\lambda))=\lambda(\E(\cU(\overline{K}_\lambda))-\cU(K)).
\end{equation}

\pass In the smooth case, we can use arguments similar to Lemma \ref{lem:geom} and Proposition \ref{prop:supportsmooth} to get
\begin{align*}
& \E(\cA(\overline{K}_{\lambda}))-\cA(K)
\underset{\lambda\to\infty}{\sim}\lambda^{-\frac{2}{3}}2^{-\frac{2}{3}}3^{-\frac1{3}}\Gamma\Big(\mfrac{2}{3}\Big)
\displaystyle \int_{\partial K}r_s^{\frac1{3}} \dd s, \\
& \E(\cU(\overline{K}_{\lambda}))-\cU(K)
\underset{\lambda\to\infty}{\sim}\lambda^{-\frac{2}{3}}2^{\frac43}3^{-\frac{4}{3}}\Gamma\Big(\mfrac{2}{3}\Big)
\displaystyle\int_{\partial K} r_s^{-\frac{2}{3}}\dd s \\
& \hspace{-3.37cm}\text{and} \\
& \E(\cN(\overline{K}_\lambda)) \underset{\lambda\to\infty}{\sim}\lambda^{\frac1{3}}2^{\frac43}3^{-\frac{4}{3}}\Gamma\Big(\mfrac{2}{3}\Big)
\displaystyle\int_{\partial K} r_s^{-\frac{2}{3}}\dd s
\end{align*}
This extends to any smooth convex body the results for the defect area and number of vertices obtained in \cite{cs05}, see Theorem 2 therein, when $K$ is a disk.

\pass In the polygonal case, we adapt Lemma \ref{lem:geompoly1}, Propositions \ref{prop:supportpoly} and \ref{prop:supportpoly2} to get
\begin{align*}
& \E(\cA(\overline{K}_\lambda))-\cA(K)
\underset{\lambda\to\infty}{\sim}\lambda^{-\frac1{2}}\displaystyle 2^{-\frac{5}{2}}\pi^{\frac{3}{2}} \sum_{i=1}^{n_K} \|a_{i+1}-a_i\|^{\frac{3}{2}}  \\
& \E(\cU(\overline{K}_{\lambda}))-\cU(K) \underset{\lambda\to\infty}{\sim}(\lambda^{-1}\log \lambda)\cdot2\cdot 3^{-1}n_K, \\
& \hspace{-3.38cm}\text{and} \\
& \E(\cN(\overline{K}_\lambda)) \underset{\lambda\to\infty}{\sim} (\log\lambda)\cdot2\cdot 3^{-1}n_K.
\end{align*}

\pass In 1968, R\'enyi and Sulanke investigated a model close to the Crofton cell, save for the fact that they did not use the notion of point process in the whole plane. Instead, they fixed a domain $B$ which includes $K$ and they considered the polygon containing $K$ and delimited by $n$ {\it random lines} which intersect $B$ without crossing $K$. This is on a par with the actual Crofton cell when the number of lines is Poissonized and the set $B$ goes to $\R^2$. In this context, they obtained the mean number of vertices in the smooth and polygonal cases, see S\"atze 4 and 5 in \cite{resu68}. Replacing $\tfrac{n}{b-l}=\tfrac{n}{\cU(B)-\cU(K)}$ by $\lambda$ in their formulas provides the exact same results as ours. To the best of our knowledge, they did not cover the calculations for the asymptotic mean area and mean perimeter, nor did they establish an Efron-type relation.

\pass It comes as no surprise that contrary to the Voronoi case, the limiting expectations do not depend on the position of $K$ with respect to $o$. Indeed, by stationarity, the origin has no privileged status among the points of the Crofton cell. As a consequence, Question $3$ is the same as Question $2$ in the context of the Crofton cell. Finally, an analogue of Question $1$ would be: what is the asymptotic shape of the Crofton cell conditional on the event that all the projections of the origin onto the lines fall outside of a fixed domain $\cD$? This problem is dependent on the position of the origin and an almost straightforward adaptation of Section \ref{sec:tele} shows that the limit shape of the Crofton cell is twice the antiorthotomic curve associated with the largest Voronoi flower with respect to $o$ included in $\cD$.

\subsection{Inlets}

\pass

\noi Regarding Question $1$, we are interested in the local structure of the cells intersecting $\partial \cD$, which look like inlets. Indeed, we can observe, for instance on Figure \ref{fig:question1}, the variability of the {\it depth} of such cells inside $\cD$ where by depth, we mean the distance to $\partial \cD$. In particular, we expect the normalized empirical distribution of the depth of all bifurcation points to have a limit. Conversely, there is a hidden branching structure, i.e. a random geometric tree, that can be observed from the origin, for instance in a fixed direction. Being able to describe this tree seems quite stimulating since it should contain information on the domain $\cD$ itself. For example, the connections occur near the medial axis of $\cD$. Finally, this question leads to the following natural issue: what is a Voronoi tessellation of a Poisson point process \textit{without any point} inside a domain $\cD$ looking like?

\bibliographystyle{plain}

\end{document}